%% file: main.tex
\documentclass[11pt,margin=3pt]{article}

\usepackage[hypertexnames=false]{hyperref}
\usepackage{color}
\usepackage{latexsym}
\usepackage{dsfont}
\usepackage{appendix}
\usepackage{amssymb}
\usepackage{subcaption}
\usepackage[margin=.9in]{geometry}
\usepackage{float}
\usepackage{algorithm,algorithmic}
\usepackage[numbers]{natbib}
\usepackage{amsmath, amsfonts,amssymb,euscript,array,enumerate,amsfonts,mathrsfs}
\usepackage[dvipsnames]{xcolor}
\usepackage{amsthm}
\newtheorem{Theorem}{Theorem}[section]

\newtheorem{Proposition}[Theorem]{Proposition}
\newtheorem{Assumption}[Theorem]{Assumption}
\newtheorem{Lemma}[Theorem]{Lemma}

\newtheorem{Remark}[Theorem]{Remark}

\usepackage[T1]{fontenc}
\usepackage[latin1]{inputenc}
\usepackage{hhline}
\usepackage{graphicx}
\usepackage{verbatim}
\usepackage{eurosym}
\usepackage{dsfont}
\allowdisplaybreaks

\DeclareMathOperator{\Tr}{Tr}

\DeclareMathOperator{\sx}{\underline{\sigma}_{\pmb{X}}}
\DeclareMathOperator{\sqmin}{\underline{\sigma}_{\pmb{Q}}}

\DeclareMathOperator{\sqi}{\underline{\sigma}_{Q,\textit{i}}}

\DeclareMathOperator{\sqfir}{\underline{\sigma}_{\pmb{Q},1}}
\DeclareMathOperator{\sqsec}{\underline{\sigma}_{\pmb{Q},2}}

\DeclareMathOperator{\sri}{\underline{\sigma}_{R,\textit{i}}}

\DeclareMathOperator{\srfir}{\underline{\sigma}_{\pmb{R},1}}
\DeclareMathOperator{\srsec}{\underline{\sigma}_{\pmb{R},2}}
\DeclareMathOperator{\srmin}{\underline{\sigma}_{\pmb{R}}}

\DeclareMathOperator{\HH}{\mathcal{H}}

\makeatletter
\@addtoreset{equation}{section}

\newcommand*{\rom}[1]{\expandafter\@slowromancap\romannumeral #1@}
\makeatother
\newcommand{\vertiii}[1]{{\left\vert\kern-0.25ex\left\vert\kern-0.25ex\left\vert #1 
    \right\vert\kern-0.25ex\right\vert\kern-0.25ex\right\vert}}

\begin{document}
\title{Policy Gradient Methods Find the Nash Equilibrium in N-player General-sum Linear-quadratic Games}

\author{Ben Hambly
\thanks{Mathematical Institute, University of Oxford. \textbf{Email:} \{hambly,  xur,  yang\}@maths.ox.ac.uk }
\and
Renyuan Xu \footnotemark[1]
\and
Huining Yang\thanks{ 
Supported by the EPSRC Centre for Doctoral Training in Industrially Focused Mathematical Modelling (EP/L015803/1) in collaboration with BP plc.}
  \footnotemark[1]
}
\maketitle
\begin{abstract}
We consider a general-sum N-player linear-quadratic game with stochastic dynamics over a finite horizon and prove the global convergence of the natural policy gradient method to the Nash equilibrium. In order to prove convergence of the method we require a certain amount of noise in the system. We give a condition, essentially a lower bound on the covariance of the noise in terms of the model parameters, in order to guarantee convergence. We illustrate our results with numerical experiments to show that even in situations where the policy gradient method may not converge in the deterministic setting, the addition of noise leads to convergence.
\end{abstract}
%\tableofcontents
\section{Introduction}
Policy optimization algorithms have achieved substantial empirical successes in addressing a variety of non-cooperative multi-agent problems, including self-driving vehicles \cite{SSS2016}, real-time bidding games \cite{JSLGWZ2018}, and optimal execution in financial markets \cite{hambly2020policy}. However, there have been few results from a theoretical perspective showing why such a class of reinforcement learning algorithms performs well with the presence of competition among agents. In the literature, the convergence of such algorithms is guaranteed only for specific classes of games including normal-form games,  differentiable games, and linear-quadratic games. For normal-form games (in which there are no state dynamics), the policy gradient method does not converge in a general set-up \cite{singh2000nash} and theoretical guarantees for convergence have been established only for some special cases such as policy prediction in two-player two-action bimatrix games \cite{zhang2010multi,song2019convergence} and two-player two-action games \cite{bowling2002multiagent}. For  differentiable games (where the cost function is assumed to be differentiable and, in most cases, the gradient is Lipschitz continuous with respect to the agent's policy parameters), there is a line of recent work \cite{balduzzi2018mechanics,letcher2019differentiable,foerster2017learning,fiez2020implicit} where the convergence guarantees for these algorithms are mainly developed for zero-sum games and cooperative games, but, in most cases, are limited to a subset of local Nash equilibrium points. However, the smoothness properties required in the differentiable game are very restrictive in general and they even fail to hold for LQ games \cite{zhang2019policy,mazumdar2020gradient,bu2019global}. 

As a starting point to tackle this challenging problem, we investigate linear-quadratic (LQ) games which can be seen as a generalization of the linear-quadratic regulator (LQR) from a single agent to multiple agents. In an LQ game, all agents jointly control a linear state process, which may be in high dimensions, where the control (or action) from each individual agent has a linear impact on the state process. Each agent optimizes a quadratic cost function which depends on the state process, the control from this agent and/or the controls from the  opponents.
%Multi-agent reinforcement learning (MARL) has enjoyed substantial successes for analyzing the otherwise challenging games, including two-agent or two-team computer games \cite{SHMGSVSAPL2016,VBCMJCDHGP2019}, self-driving vehicles \cite{SSS2016}, real-time bidding games \cite{JSLGWZ2018}, ride-sharing \cite{LQJYWWWY2019}, and traffic routing \cite{EAA2013}.

LQ games are a relatively simple setting in which to analyze the behavior of multi-agent reinforcement learning (MARL) algorithms in continuous action and state
spaces since they admit global Nash equilibria in the space of linear feedback policies. Moreover, these equilibria can be found by solving a coupled set of Ricatti equations when system parameters are given. As such LQ games are a natural benchmark problem on which to test policy gradient algorithms in multi-agent settings when system parameters are unknown. Furthermore, policy gradient methods open up the possibility to develop new scalable approaches to finding or learning solutions to control problems even with constraints.
Finally, the empirical results presented in \cite{mazumdar2019policy} imply that, even in this relatively straightforward LQ case (with linear dynamics, linear feedback policies, and quadratic costs), policy gradient MARL would be unable to find the local Nash equilibrium in a non-negligible subset of problems. This further demonstrates the necessity of understanding under what circumstances policy gradient methods work for LQ games.

For LQ games with policy gradient algorithms, most of the existing literature has focused on zero-sum games with two players \cite{bu2019global,zhang2019policy,zhang2021derivative}. In the setting of deterministic dynamics and infinite time horizon,  \cite{zhang2019policy}  proposed an alternating policy update scheme with a projection step and showed sublinear convergence for the algorithm. For a similar setting, \cite{bu2019global} proposed a leader-follower type of policy gradient algorithm  which is projection-free and  enjoys a global sublinear convergence rate and asymptotically linear convergence rate.  For the case of stochastic dynamics and finite time horizon, \cite{zhang2021derivative}  provided the first sample complexity result for the global convergence of the policy gradient method with an alternating policy update scheme.

However, little theory has been developed  for the more general class of LQ games with  N players and general-sum cost functions. It has been documented  that the policy gradient method may {\it fail to converge} in such a setting with deterministic dynamics due to the lack of theoretical guidance on how to properly choose the step size and the exploration scheme \cite{mazumdar2019policy}.  For a special class of N-player LQ games  with homogeneous agents where agents interact with each other through a  mean-field type of ``deep state'' which is the aggregated state position of all agents,
\cite{roudneshin2020reinforcement} showed that the policy gradient method converges to the global Nash equilibrium. Up to now, as far as we are aware, providing theoretical guarantees for the convergence of the policy gradient method remains an open problem for general-sum LQ games with more than two players
%. This has also been recognized as a open problem in 
\cite{mazumdar2020gradient}.

%{\color{red}[Discussion: review on other RL methods for zero-sum LQGs? review on policy gradient method for zero-sum MDP games?]}

%In this work, we focus on a more general  but less understood class of games which is the  N-player general-sum LQGs. \cite{mazumdar2019policy} \cite{roudneshin2020reinforcement} \cite{mazumdar2020gradient}

\paragraph{Our Contributions.} In this work, we explore the natural policy gradient method, which can be viewed as a normalized version of the vanilla policy gradient descent method, for a class of N-player general-sum linear-quadratic games. Our main result is Theorem~\ref{thm:local_conv_NPG} in which we provide a global linear convergence guarantee %and a polynomial sample complexity guarantee {\color{blue} Not sure I really understand what is meant by this}
for this approach in the setting of a finite time horizon and stochastic dynamics provided there is a certain level of noise in the system. The noise can either come from the underlying dynamics or carefully designed explorations from the agents.  Intuitively speaking, the noise can help the agents escape from %unpleasant situations 
traps such as those observed in \cite{mazumdar2019policy} when the  dynamics are deterministic. In addition, the noise can help agents find the correct descent direction for convergence to the desired Nash equilibrium point. To the best of our knowledge, this is the first result of its kind in the  MARL literature showing a provable convergence result for N-player general-sum LQ games. From a technical perspective, the main difficulty is to quantify the perturbation of the individual's gradient term (see Step 3 in the proof of Lemma \ref{lemma:one_step_loc_conv}) and to control the descent of the individual's cost function (see Step 4 in the proof of Lemma \ref{lemma:one_step_loc_conv}) in the presence of competition from the other $N-1$ agents. We also note that our main result for the natural policy gradient method can be extended to the case of the vanilla policy gradient method, see Theorem \ref{thm:conv_PG_vani} and Lemma \ref{lemma:one_step_cont_vani}.

We illustrate the performance of our algorithm with three examples. We first perform the natural policy gradient algorithm under the experimental set-up in \cite{mazumdar2019policy} over a finite time horizon. The algorithm converges to the Nash equilibrium for appropriate initial policies and step sizes. The second example is a toy LQ game example with synthetic data. The empirical results suggest that, in practice, the natural policy gradient algorithm can find the Nash equilibrium even if the level of system noise is lower than that required for our theoretical analysis. The third example is a three-player general-sum game, where we show the convergence of natural policy gradient methods with known and unknown parameters.

\paragraph{Comparison to the Literature on General-sum LQ games.} With deterministic dynamics and infinite horizon, \cite{mazumdar2019policy} provided some empirical examples where the policy gradient method fails to converge to the set of Nash equilibria (which may not be unique).  %under such setting with deterministic dynamics due to the lack of theoretical guidance on how to  \cite{mazumdar2019policy}
These empirical examples motivate an examination of the possibility of applying policy gradient methods in the multi-agent environment. Here we explain the difference between our framework and the set-up in \cite{mazumdar2019policy}. In addition, we offer some explanations for why the policy gradient method works in our framework whereas it fails to converge in \cite{mazumdar2019policy}. 
\begin{itemize}
    \item {\bf Well-definedness of LQ games:} The existence and uniqueness of a Nash equilibrium is a prerequisite for the convergence of learning algorithms. In the setting of stochastic dynamics and finite horizon, the general-sum LQ game has a unique Nash equilibrium solution under mild conditions \cite{BasarOlsder1999}. For the setting with an infinite horizon, the existence of Nash equilibria can be proved under some stabilizability and detectability properties. However obtaining explicit and verifiable model conditions for stabilizability and detectability seems to be a quite challenging task let alone finding conditions for the uniqueness of the equilibrium \cite{BasarOlsder1999}.  It is not clear if there exists a unique Nash equilibrium for the setting considered in \cite{mazumdar2019policy}.
    \item  {\bf Self-exploration property of time-dependent policies:} For single-agent LQR problems, \cite{basei2021logarithmic} highlighted that the time-dependent optimal feedback policy enjoys a self-exploration property in the finite time horizon setting. Namely, the time-dependent optimal feedback matrices ensure that the optimal state and control processes span the entire parameter space, which enables the design of efficient exploration-free learning algorithms. By contrast the optimal feedback policy is time-invariant in the infinite time horizon setting and learning algorithms tend to have difficulty converging without efficient exploration schemes \cite{mania2019certainty}. We believe a similar analogy holds for the game setting as well.
    \item {\bf System noise:} In our setting we need Assumption \ref{ass:noise_cond} which implies that a certain level of system noise is essential for the convergence of the policy gradient method. We show via numerical experiments in Section \ref{sec:numerical_experiments} that the circulating and divergence phenomenon described in \cite{mazumdar2019policy} can be avoided with the addition of system noise. %{\color{red} already said} Intuitively the system noise can help the agents avoid weird situations in the competition and find the correct descent direction to reduce the cost. 
 The noise can either come from the original system when the dynamics are stochastic (as suggested in Assumption \ref{ass:noise_cond}) or agents can apply Gaussian exploration (as suggested in \cite{mania2019certainty} for single-agent LQR problems with infinite time horizon).
\end{itemize}
%\paragraph{Comparison to \cite{roudneshin2020reinforcement} .}
In addition, \cite{roudneshin2020reinforcement} showed the global convergence of the policy gradient method for a mean-field type of LQ game in the setting of infinite horizon and stochastic dynamics. In particular, agents are assumed to be  homogeneous and are only able to interact 
through an aggregated state and action pair. With this special formulation, the uniqueness of the Nash equilibrium could be established (see \cite[Theorem 1]{roudneshin2020reinforcement}) and the proof of convergence could be reduced to the single agent case  (see \cite[Theorem 2]{roudneshin2020reinforcement}). In this paper, we focus on a more general LQ game with no homogeneity assumption nor any restriction on the interactions.

%In this work, we show that natural policy gradient, which can be viewed as a normalized version of the vanilla policy gradient descent method, has provably local convergence guarantee to the global Nash equilibrium with stochastic dynamics and finite time horizon. In addition, we discover that 

%\paragraph{Policy gradient for non-cooperative games.}
%\begin{itemize}
%    \item Markov games: \cite{zhao2021provably}, \cite{zhang2020model},\cite{mguni2021learning},\cite{wei2021last}
%\end{itemize}

\paragraph{Organization and Notation} For any matrix $Z =(Z_1,\cdots,Z_d) \in \mathbb{R}^{m\times d}$ with $Z_j \in \mathbb{R}^m$ ($j=1,2,\cdots, d$), we let $Z^{\top} \in \mathbb{R}^{d\times m}$ denote the transpose of $Z$, $\|Z\|$ denotes the spectral norm of the matrix $Z$; $\Tr(Z)$ denotes the trace of a square matrix $Z$; and $\sigma_{\min}(Z)$ denotes the minimal singular value of a square matrix $Z$. For a sequence of matrices $\pmb{D}=(D_0,\cdots,D_T)$, we define a new norm $\vertiii{\pmb{D}}$ as $\vertiii{\pmb{D}} = \sum_{t=0}^T \|D_t\|$, where $D_t \in \mathbb{R}^{m\times d}$; $\gamma_D=\max_{t=0,\cdots,T}\|D_t\|$ denotes the maximum over all $\|D_t\|$. Furthermore we denote by $\mathcal{N}(\mu,\Sigma)$ the Gaussian distribution with mean $\mu\in \mathbb{R}^d$ and covariance matrix $\Sigma\in\mathbb{R}^{d\times d}$.

The rest of the paper is organized as follows. We introduce the mathematical framework and problem set-up in Section \ref{sec:setup}. The convergence analysis of the natural policy gradient method for the case of known model parameters is provided in Section \ref{sec:exact_NPG}. When parameters are unknown,  the sample-based natural policy gradient method is discussed in Section \ref{sec:mf_NPG}. Finally, the algorithm is applied to three numerical examples in Section \ref{sec:numerical_experiments}.

\section{N-player General-Sum Linear-quadratic Games}\label{sec:setup}

\subsection{Problem Set-up}
We consider the following N-player general-sum linear-quadratic (LQ) game over a finite time horizon $T$. The state process evolves as 
\begin{equation}\label{eqn:nonzero_State}
    x_{t+1} = A_tx_t + \sum_{i=1}^N B_t^iu_t^i + w_t,\ t = 0,1,\cdots,T-1,
\end{equation}
where $x_t\in\mathbb{R}^d$ is the state of the system with the initial state $x_0$ drawn from a Gaussian distribution, $u_t^i\in\mathbb{R}^{k_i}$ is the control of player $i$ at time $t$ and $\{w_t\}_{t=0}^{T-1}$ are zero-mean IID Gaussian random variables which are independent of $x_0$. In this paper we focus on the class of Markovian controls. The system parameters
$A_t\in\mathbb{R}^{d\times d}$, $B_t^i\in\mathbb{R}^{d\times k_i}$, for $t=0,1,\cdots,T-1$ are referred to as system (transition) matrices. The objective of player $i$ ($i=1,\cdots,N$) is to minimise their finite time horizon value function:
\begin{equation}\label{eqn:nonzero_cost}
    \inf_{\{u_t^i\}_{t=0}^{T-1}}\mathbb{E}\left[\sum_{t=0}^{T}c_t^i(x_t,u_t^i)\right],
\end{equation}
where the cost function 
\begin{equation}\label{eqn:nonzero_cost_at_t}
    c_t^i(x_t,u_t^i)= x_t^\top Q_t^ix_t+(u_t^i)^\top R_t^{i} u_t^i,\quad t = 0,1,\cdots, T-1,
\end{equation}
with $c_T^i(x_T)=x_T^\top Q_T^ix_T$, where $Q_t^i\in\mathbb{R}^{d\times d}$ and $R_t^{i}\in\mathbb{R}^{k_i\times k_i}$ ($i=1,\cdots,N$) are matrices that parameterize the quadratic costs. Note that the randomness in the LQ game comes from both the initial state and the noise process in the state equation, therefore throughout the paper, unless specified otherwise, the expectation (for example in \eqref{eqn:nonzero_cost}) is taken with respect to both the initial state $x_0$ and the noise $\{w_t\}_{t=0}^{T-1}$. We also denote by $\pmb{u}^i := (u_0^i,\cdots,u_{T-1}^i)$, $\pmb{x} := (x_0,\cdots,x_{T})$, $\pmb{Q}^i := (Q_0^i,\cdots,Q_{T}^i)$, and $\pmb{R}^{i} := (R_0^{i},\cdots,R_{T-1}^{i})$, for $i=1,\cdots,N$.

\begin{Assumption}[Cost Parameter]\label{ass:nonzero_cost}
Assume for $i=1,\cdots,N$, $Q_t^i\in\mathbb{R}^{d\times d}$, for $t=0,1,\cdots,T$, and $R_t^{i}\in\mathbb{R}^{k_i\times k_i}$, for $t=0,1,\cdots,T-1$ are symmetric positive definite matrices. 
\end{Assumption}

\begin{Assumption}[Initial State and Noise Process]\label{ass:nonzero_initial_noise}
Assume 
\begin{enumerate}
    \item Initial state: $x_0$ is Gaussian such that $\mathbb{E}[x_0x_0^\top]$ is positive definite.
    \item Noise: $\{w_t\}_{t=0}^{T-1}$ are IID Gaussian and independent from $x_0$ such that $\mathbb{E}[w_t] = 0$, and $W = \mathbb{E}[w_tw_t^\top]$ is
positive definite, $\forall t=0,1,\cdots,T-1$.
\end{enumerate}
\end{Assumption}

\begin{Assumption}[Existence and Uniqueness of Solution]\label{ass:nonzero_exist_sol}
Assume there exists a unique solution set $\{K_t^{i*}\}_{t=0}^{T-1}$, for $i=1,\cdots,N$ to the following set of linear matrix equations:
\begin{equation}\label{eqn:optimal_K_ti}
     K_t^{i*} = \left(R_t^{i}+(B_t^i)^\top P_{t+1}^{i*}B_t^i\right)^{-1}(B_t^i)^\top P_{t+1}^{i*}\left(A_t-\sum_{j=1,j\neq i}^N B_t^jK_t^{j*}\right),
\end{equation}
where $\{P_t^{i*}\}_{t=0}^T$ are obtained recursively backwards from
\begin{equation}\label{eqn:lyapunov_opti}
    P_{t}^{i*} = Q_t^i + (K_t^{i*})^\top R_t^{i} K_t^{i*} + \left(A_t - \sum_{j=1}^N B_t^jK_t^{j*} \right)^\top P_{t+1}^{i*}\left(A_t -  \sum_{j=1}^N B_t^jK_t^{j*}\right),
\end{equation}
with terminal condition $P_{T}^{i*} = Q_T^i$.
\end{Assumption}

A similar assumption is adopted in \cite{zhang2019policy} for a two-player zero-sum LQ game and in \cite{roudneshin2020reinforcement} for a homogeneous N-player game with mean-field interaction.
\begin{Remark}\label{remark:suff_unique}
%[A sufficient condition for Assumption \ref{ass:nonzero_exist_sol}]
{\rm
A sufficient condition for the unique solvability of \eqref{eqn:optimal_K_ti} is the invertibility of the block matrix $\Phi_t$, $t=0,1,\cdots,T-1$, with the $ii$-th block given by $R_t^i+(B_t^i)^\top P_{t+1}^{i*}B_t^i$ and the $ij$-th block given by $(B_t^i)^\top P_{t+1}^{i*} B_t^j$, where $i,j=1,\cdots,N$ and $j\neq i$. See Remark 6.5 in \cite{BasarOlsder1999}.}
\end{Remark}

\begin{Lemma}[Nash Equilibrium and Equilibrium Cost {\cite[Corollary 6.4]{BasarOlsder1999}}]\label{lemma:nash_equilibrium}
Assume Assumptions \ref{ass:nonzero_cost}, \ref{ass:nonzero_initial_noise}, and \ref{ass:nonzero_exist_sol} hold. Then for $i=1,2,\cdots,N$,
\begin{enumerate}
    \item The Nash equilibrium strategy for player $i$ is given by
    \begin{equation}\label{eqn:nash_linear_feedback}
        u_t^{i*} = - K_t^{i*} x_t,\quad t=0,1,\cdots T-1,
    \end{equation}
    where $K_t^{i*}$ is defined in \eqref{eqn:optimal_K_ti}.
    \item The Nash equilibrium cost for player $i$ is  
\begin{equation}\label{eqn:nonzero_equil_cost}
 \mathbb{E}[x_0^\top P_0^{i*}x_0+N_0^{i*}] = \inf_{\{u_t^i\}_{t=0}^{T-1}}\mathbb{E}\left[\sum_{t=0}^T c_t^i(x_t,u_t^i)\right],
\end{equation}
where $\{P_t^{i*}\}_{t=0}^{T}$ are defined in \eqref{eqn:lyapunov_opti} and
\begin{equation}\label{eqn:nonzero_Nt}
    N_t^{i*} = N_{t+1}^{i*} + \mathbb{E}[w_t^\top P_{t+1}^{i*}w_t] = N_{t+1}^{i*} + \Tr(WP_{t+1}^{i*}),\quad t=0,1,\cdots,T-1
\end{equation}
with terminal condition $N_T^{i*} = 0$.
\end{enumerate}
\end{Lemma}

To find the Nash equilibrium strategy in the linear feedback form \eqref{eqn:nash_linear_feedback}, we only need to focus on the following class of linear admissible policies in feedback form
\begin{equation*}
    u_t^i = -K_t^ix_t, \qquad t=0,1,\cdots,T-1,\quad i=1,\cdots,N
\end{equation*}
which can be fully characterized by $\pmb{K}^i := (K_0^i,K_1^i,\cdots,K_{T-1}^i)$ with $K_t^i\in\mathbb{R}^{k_i\times d}$. We write $\pmb{K}=(\pmb{K}^1,\cdots,\pmb{K}^N)$ for a collection of policies and 
\[
\pmb{K}^{*}=(\pmb{K}^{1*},\cdots,\pmb{K}^{N*}).
\]
for the collection of optimal policies. We will use the notation 
\begin{equation*}
    (\pmb{K}^i,\pmb{K}^{-i*})=(\pmb{K}^{1*},\cdots,\pmb{K}^{(i-1)*},\pmb{K}^{i},\pmb{K}^{(i+1)*},\cdots,\pmb{K}^{N*}),
\end{equation*}
for $i=1,\cdots,N$ for player $i$'s policy, when all other players use their optimal policies.

\section{The Natural Policy Gradient Method with Known Parameters} 
\label{sec:exact_NPG}

In this section, we provide a global linear convergence guarantee %and a polynomial sample complexity guarantee 
for the natural policy gradient method applied to the LQ game \eqref{eqn:nonzero_State} - \eqref{eqn:nonzero_cost}. Throughout this section we assume all the parameters in the LQ game, %\eqref{eqn:nonzero_State} - \eqref{eqn:nonzero_cost}, 
$\{A_t\}_{t=0}^{T-1}$, $\{B_t^i\}_{t=0}^{T-1}$, $\{Q_t^i\}_{t=0}^{T}$, and $\{R_t^{i }\}_{t=0}^{T-1}$ ($i=1,\cdots,N$), are known. The analysis with known parameters paves the way for learning LQ games with unknown parameters which is discussed in Section \ref{sec:mf_NPG}.

For LQ games, any (admissible) feedback policy can be fully characterized by a set of parameters $\pmb{K}=(\pmb{K}^1,\cdots,\pmb{K}^N)$. Therefore, we can correspondingly define player $i$'s cost induced by the joint policy $\pmb{K}$ as
\begin{equation*}
    C^i(\pmb{K}) = \mathbb{E}\left[\sum_{t=0}^{T-1}\left(x_t^\top Q_t^i x_t+(K_t^i x_t)^\top R_t^{i}(K_t^ix_t)\right)+x_T^\top Q_T^ix_T\right],
\end{equation*}
where $\{x_t\}_{t=0}^T$ is the random path from the dynamics \eqref{eqn:nonzero_State} induced by $\pmb{K}$ starting with $x_0$. 

We start by introducing some notation which will be used throughout the analysis.
We define the state covariance matrix $\Sigma_t^{\pmb{K}}$, and let $\Sigma_{\pmb{K}}$ be the sum of $\Sigma_t^{\pmb{K}}$:
\begin{equation}\label{eqn:defn_Sigma_t_K1K2}
    \Sigma_t^{\pmb{K}} = \mathbb{E}[x_t^{\pmb{K}} (x_t^{\pmb{K}})^\top], \quad \Sigma_{\pmb{K}} = \sum_{t=0}^T \Sigma_t^{\pmb{K}},
\end{equation}
where $\{x_t^{\pmb{K}}\}_{t=0}^T$ is a state trajectory generated by following a set of policies $\pmb{K}$. We will write $x_t=x_t^{\pmb{K}}$ when no confusion may occur. Define $\sigma_{\pmb{X}}^{\,\pmb{K}}$ to be the lower bound over all the minimum singular values of $\Sigma_t^{\pmb{K}}$:
\[
\sigma_{\pmb{X}}^{\,\pmb{K}}:=\min_t \sigma_{\min}(\Sigma_t^{\pmb{K}}),
\]
%where $\{x_t^{\pmb{K}}\}_{t=0}^T$ is a state trajectory generated by following a set of policies $\pmb{K}$. 
We also define $\sx$ as
\begin{equation}\label{defn_sx_nonzero}
   \sx  := \min\{\sigma_{\min}(\mathbb{E}[x_0x_0^\top]),\sigma_{\min}(W)\}.
\end{equation}

Similarly, we define 
\[ \sri=\min_t \sigma_{\min}(R_t^{i}), \;\; \sqi=\min_t\sigma_{\min}(Q_t^i),\]
and 
\begin{equation}\label{eqn:defn_srmin_sqmin}
 \srmin = \min_i\{\sri\}, \quad \sqmin = \min_i\{\sqi\}.
\end{equation}
We further define $\gamma_A$, $\gamma_B$, and $\gamma_R$ as
\begin{equation}\label{eqn:defn_gammaAgammaBgammaR}
     \gamma_A = \max_{t=0,\cdots,T-1} \|A_t\|,\quad \gamma_B =\max_i\{ \max_{t=0,\cdots,T-1} \|B_t^i\|\}, \quad\gamma_R =\max_i\{ \max_{t=0,\cdots,T-1} \|R_t^i\|\}.
\end{equation}
Under Assumption \ref{ass:nonzero_cost}, we have $\sri\geq \srmin>0$ and $\sqi\geq \sqmin>0$, for $i=1,\dots,N$. For the well-definedness of the state covariance matrix, we have the following result.

\begin{Lemma}
%[Well-definedness of the State Covariance Matrix]
Assume Assumption \ref{ass:nonzero_initial_noise} holds. Then $\mathbb{E}[x_tx_t^\top]$ is positive definite for $t=0,1,\cdots,T$ under any set of policies $\pmb{K}=(\pmb{K}^1,\cdots,\pmb{K}^N)$ and we have $\sigma_{\pmb{X}}^{\,\pmb{K}}\geq \sx>0$.
\end{Lemma}

\begin{proof}
Let $\{x_t\}_{t=0}^T$ be the state trajectory induced by an arbitrary policy set $\pmb{K}$. By Assumption \ref{ass:nonzero_initial_noise} the matrix $\mathbb{E}[x_0 x_0^\top ]$ is positive definite. For $t \ge 1$, we have from \eqref{eqn:nonzero_State} and taking expectations
$$\mathbb{E}[x_t x_t^\top ] =\left(A_{t-1}-\sum_{i=1}^N B_{t-1}^iK_{t-1}^i\right) \mathbb{E}[x_{t-1}x_{t-1}^\top]\left(A_{t-1}-\sum_{i=1}^N B_{t-1}^iK_{t-1}^i\right)^\top + \mathbb{E}[w_{t-1} w_{t-1}^\top].$$ 
Now as $(A_{t-1}-\sum_{i=1}^N B_{t-1}^iK_{t-1}^i) \mathbb{E}[x_{t-1}x_{t-1}^\top](A_{t-1}-\sum_{i=1}^N B_{t-1}^iK_{t-1}^i)^\top$ is positive semi-definite and by assumption $\mathbb{E}[w_{t-1} w_{t-1}^\top]$ is positive definite, we have $\mathbb{E}[x_t x_t^\top ]$ is positive definite and as a result the statement holds. 
% In this case, we can simply take  $\sx = \min\{\sigma_{\min}(\mathbb{E}[x_0x_0^\top]),\sigma_{\min}(W)\}$.
\end{proof}

We write  $\HH =\left\{ h\,|\, h \,\,\text{are polynomials in the model parameters}\right\}$ and $\HH(.)$ when there are other dependencies. The model parameters are expressed in terms of 
 $\frac{1}{1+\sum_i k_i}$, $\frac{1}{d+1}$, $d$, $\frac{1}{N+1}$, $\frac{1}{T+1}$, $\frac{1}{\|W\|+1}$,$\frac{1}{\|\Sigma_0\|+1}$, $\frac{1}{\gamma_A+1}$, $\frac{1}{\gamma_B+1}$, $\frac{1}{\gamma_R+1}$, $\frac{1}{\sx+1}$, $\sx$,  $\frac{1}{\srmin+1}$, $\srmin$, $\frac{1}{\sqmin+1}$, $\sqmin$, and $\frac{1}{\vertiii{\pmb{K}^*}+1}$.
 
\paragraph{Natural Policy Gradient Method.} We consider the following natural policy gradient updating rule for each player $i$ ($i=1,\cdots,N$), in a sequence of games for $m=1,\dots, M$,
\begin{equation}\label{eqn:nonzero_grad_update_rule}
    K_{t}^{i,(m)} = K_t^{i,(m-1)} - \eta\nabla_{K_t^i}C^i(\pmb{K}^{(m-1)})(\Sigma_t^{\pmb{K}^{(m-1)}})^{-1},\quad\forall\, 0\leq t\leq T-1,
\end{equation}
where %$m$ ($1\leq m\leq M$) is the number of iterations, 
$\nabla_{K_t^i}C^i(\pmb{K}^{(m-1)})=\frac{\partial C^i(\pmb{K}^{(m-1)})}{\partial K_t^i}$ is the gradient of $C^i(\pmb{K}^{(m-1)})$ with respect to $K_t^i$, and $\eta$ is the step size.

Natural policy gradient methods \cite{kakade2001natural} -- and related algorithms such as trust region policy optimization \cite{schulman2015trust} and the natural actor critic \cite{peters2008natural} -- are some of the most popular and effective variants of the vanilla policy gradient methods in single-agent reinforcement learning.  The natural policy gradient method performs better than the vanilla policy gradient method in practice since it takes the information  geometry (\cite{bagnell2003covariant,kakade2001natural}) into consideration by normalizing the gradient term by $(\Sigma_t^{\pmb{K}^{(m)}})^{-1}$ in \eqref{eqn:nonzero_grad_update_rule}. By doing so, the natural gradient method represents the steepest descent direction based on the underlying structure of the parameter space \cite{kakade2001natural}. It is worth mentioning that \cite{FGKM2018} provided the first theoretical guarantee that the natural gradient method has an  improved constant in the convergence rate compared to the vanilla version.

%{\color{red}[Huining: we may want to mention the advantages of NPG compared with vanilla PG, e.g. Fazel 2018 (NPG helps simplify the proof and leads to smaller $N$ and simple $\eta$ but I didn't see any order improvement.)]}

 With access to the model parameters, players can directly calculate their own policy gradient and perform the natural policy gradient steps iteratively. See Algorithm \ref{alg:NPG_known} for the details. Unlike the nested-loop updates in some zero-sum LQ game settings \cite{zhang2020sample,zhang2019policy} in which the inner-loop updates are sample inefficient, each player updates their policies simultaneously at each iteration in Algorithm \ref{alg:NPG_known}. This 
simultaneous updating framework is a more realistic set-up for practical examples such as online auction bidding.

\paragraph{Equivalent Form of the Natural Policy Gradient.}  Following the explanation of the natural  gradient method in the single-agent setting \cite{FGKM2018}, we provide an equivalent form of the  natural gradient method with Fisher information in the game setting. In our N-player game case, we assume player $i$'s updating rule follows:
\begin{eqnarray}
K_t^{i\prime} \longleftarrow K_t^{i} - \eta\, G_{K_t^i}^{-1}\nabla _{K_t^{i}}C^i(\pmb{K}),\label{eq:natural_fisher_version}
\end{eqnarray}
where  $G_{K_t^i}$ is the Fisher information matrix:
\begin{eqnarray}
G_{K_t^i} = \mathbb{E}\Big[\nabla_{K_t^i}\log \pi^i_{K_t^i}(u_t|x_t)\nabla_{K_t^i}\log \pi^i_{K_t^i}(u_t|x_t)^{\top}\Big]
\end{eqnarray}
under a {\it linear} policy
with {\it additive} Gaussian noise \cite{rajeswaran2017towards}, in that 
\begin{eqnarray}\label{eq:Gaussian_policy}
\pi^i_{K_t^i}(u^i_t = u|x_t) = {\rm det} (2\pi\sigma^2 I)^{-1/2} \exp\left(-\frac{1}{2\sigma^2}\left(u-K_t^ix_t\right)^{\top}\left(u-K_t^ix_t\right)\right). %\mathcal{N}(K_t^ix_t,\sigma^2 I).
\end{eqnarray}

By a similar analysis to that in \cite{FGKM2018}, we can show that the Fisher information matrix of size $k_id \times k_id$, which is indexed as $[G_K]_{(j,q),(j',q')}$ where $j,j'\in\{1,2,\cdots,k_i\}$ and $q,q';\in\{1,2,\cdots,d\}$, has a block diagonal form where the only non-zeros blocks are
$[G_{K_t^i}]_{(j,\cdot),(j,\cdot)} = \Sigma_{t}^{\pmb{K}}=\mathbb{E}[x_t^{\pmb{K}}(x_t^{\pmb{K}})^{\top}]$ (this is the block corresponding to the i-th coordinate of the action, as j ranges from 1
to $k_i$).  Hence \eqref{eq:natural_fisher_version} is equivalent to the following updating rule
$K_t^{i\prime} \longleftarrow K_t^{i} - \eta \nabla _{K_t^{i}}C^i(\pmb{K}) \Big(\Sigma_{t}^{\pmb{K}}\Big)^{-1}$,
which is equivalent to \ref{eqn:nonzero_grad_update_rule}.

\begin{algorithm}[H]
\caption{\textbf{Natural Policy Gradient Method with Known Parameters}}
\label{alg:NPG_known}
\begin{algorithmic}[1]
    \STATE \textbf{Input}: Number of iterations $M$, time horizon $T$, initial policies $\pmb{K}^{(0)}=(\pmb{K}^{1,(0)},\cdots,\pmb{K}^{N,(0)})$, step size $\eta$, model parameters $\{A_t\}_{t=0}^{T-1}$, $\{B_t^i\}_{t=0}^{T-1}$, $\{Q_t^i\}_{t=0}^{T}$, and $\{R_t^i\}_{t=0}^{T-1}$ ($i=1,\cdots,N$).
        \FOR {$m\in\{1, \ldots, M\}$}
         \FOR {$t\in\{T-1, \ldots, 0\}$}
         \FOR {$i\in\{1, \ldots, N\}$}
          \STATE Calculate the matrix $P_{t,i}^{\pmb{K}^{(m-1)}}$ by
          \begin{eqnarray}\label{alg:cal_riccati}
              P_{t,i}^{\pmb{K}^{(m-1)}} &=& Q_t^i + (K_t^{i,(m-1)})^\top R_t^{i} K_t^{i,(m-1)}\nonumber \\
              &&+ \left(A_t -\sum_{i=1}^N B_t^iK_t^{i,(m-1)}\right)^\top P_{t+1,i}^{\pmb{K}^{(m-1)}}\left(A_t -\sum_{i=1}^N B_t^iK_t^{i,(m-1)} \right)
          \end{eqnarray}
          with $P_{T,i}^{\pmb{K}^{(m-1)}}=Q_T^i$.
          \STATE Calculate the matrix $E_t^i$ by 
          \begin{equation}\label{alg:cal_grad}
              E_{t,i}^{\pmb{K}^{(m-1)}} = R_t^{i}K_t^{i,(m-1)}-(B_t^i)^\top P_{t+1,i}^{\pmb{K}^{(m-1)}}\left(A_t-\sum_{i=1}^NB_t^iK_t^{i,(m-1)}\right).
          \end{equation}   
             \STATE Update the policies using the natural policy gradient updating rule: 
            \begin{equation}\label{eqn:algo_known_update}
                 K_t^{i,(m)} = K_t^{i,(m-1)} - 2\eta E_{t,i}^{\pmb{K}^{(m-1)}}.
             \end{equation}
              \ENDFOR
            \ENDFOR
        \ENDFOR
 \STATE Return the iterates $\pmb{K}^{(M)}=(\pmb{K}^{1,(M)},\cdots,\pmb{K}^{N,(M)})$.
\end{algorithmic}
\end{algorithm}
\begin{Remark}
{\rm In Algorithm \ref{alg:NPG_known}, during iteration $m$, each player first calculates the solution to the backward Riccati equation in \eqref{alg:cal_riccati} based on the observations from the previous iteration $m-1$, and obtains the policy gradients in \eqref{alg:cal_grad}. The natural gradient method is then applied in \eqref{eqn:algo_known_update} to update the policy for that iteration. Note that \eqref{eqn:algo_known_update} does not involve $(\Sigma_t^{\pmb{K}})^{-1}$ since, as we show later in Lemma \ref{lemma:nonzero_policygrad}, we have $2E_{t,i}^{\pmb{K}}=\nabla_{K_t^i}C^i(\pmb{K})(\Sigma_t^{\pmb{K}})^{-1}$.

A straightforward analysis shows that the computational complexity of Algorithm \ref{alg:NPG_known} is\\ $\mathcal{O}(MNT \max(k_i^2 d, d^3))$ and that it requires $(Td \sum_{i=1}^N k_i)$  storage units for $\pmb{K}$.}
\end{Remark}

We now introduce the main assumption and the main result for the natural policy gradient method applied to the class of general-sum LQ games. To start we define $\rho^*$ as
\begin{eqnarray}\label{eqn:defn_rho_star}
\rho^* := \max\left\{\max_{0\leq t \leq T-1}\left\|A_t-\sum_{i=1}^N B_t^iK_t^{i*}\right\|, 1+ \delta\right\},
\end{eqnarray}
for some small constant $\delta>0$, and define $\psi := \max_i\{C^i(\pmb{K}^{i,(0)},\pmb{K}^{-i*})-C^{i}(\pmb{K}^*)\}$. 
Now set
\begin{equation}\label{eqn:defn_rho_bar}
\bar{\rho}:= \rho^*+N\gamma_B\sqrt{\frac{T\psi}{\sx\srmin}} + \frac{1}{20 T^2}.
\end{equation}

\begin{Assumption}[System Noise]\label{ass:noise_cond} The system parameters satisfy the following inequality
\begin{equation}\label{eq:system_noise}
\frac{(\sx)^5}{\|\Sigma_{\pmb{K}^{*}}\| } >20(N-1)^2\,T^2\,d\, \frac{(\gamma_B)^4(\max_i\{C^i(\pmb{K}^*)\}+\psi)^4}{\sqmin^2\srmin^2}\left(\frac{\bar{\rho}^{2T}-1}{\bar{\rho}^2-1}\right)^2.
\end{equation}
%where $\sx$ is defined in \eqref{defn_sx_nonzero}, $\srmin$ and $\sqmin$ are defined in \eqref{eqn:defn_srmin_sqmin}, $\gamma_B$ is defined in \eqref{eqn:defn_gammaAgammaBgammaR}, and $\bar{\rho}$ is defined in \eqref{eqn:defn_rho_bar}.
\end{Assumption}
Note that $\sx$ appears on both sides of \eqref{eq:system_noise} (indirectly through $\bar{\rho}$ on the R.H.S.) and when $\sx$ increases, the LHS of \eqref{eq:system_noise} increases while the RHS of \eqref{eq:system_noise} decreases. Therefore, Assumption \ref{ass:noise_cond} requires $\sx$ to be large enough such that inequality \eqref{eq:system_noise} holds. This can be interpreted as ensuring that the system needs a certain level of noise for the learning agents to find the correct direction towards the Nash equilibrium.  Assumption \ref{ass:noise_cond} also imposes conditions on the initial policy $\pmb{K}^{(0)}$ as $\psi$ depends on the difference between $C^i(\pmb{K}^{i,(0)},\pmb{K}^{-i*})$ and $C^{i}(\pmb{K}^*)$.

Below we provide two examples such that Assumption~\ref{ass:noise_cond} is satisfied. In Section \ref{sec:numerical_experiments} we will show that, at least in some circumstances, the natural policy gradient method leads to the Nash equilibrium even when Assumption~\ref{ass:noise_cond} is violated.  This further demonstrates the power of the natural policy gradient method in practice.
\begin{Remark}[Examples and Discussion of Assumption \ref{ass:noise_cond}.]
\label{app:sys_noise_ass}{\rm
\begin{enumerate}
    \item A two-player example where Assumption \ref{ass:noise_cond} is satisfied:
    \begin{itemize}
%     \item Parameters: 
%     \[
% A_t = 
% \begin{bmatrix}
% 0.2 & -0.1\\
% 0.05 & 0.2
% \end{bmatrix},\quad
% B_t^1 = 
% \begin{bmatrix}
% 0.1 \\
% 0.02
% \end{bmatrix},\quad
% B_t^2 = 
% \begin{bmatrix}
% 0.04 \\
% -0.1
% \end{bmatrix},\quad
% W = 
% \begin{bmatrix}
% 0.3 & 0.1 \\
% 0.1 & 0.2
% \end{bmatrix},
% \]
%     \[
% Q_T^1=Q_T^2 =Q_t^1=Q_t^2 = 
% \begin{bmatrix}
% 0.1 & 0\\
% 0 & 0.1
% \end{bmatrix},\quad
% R_t^1(t)=R_t^2(t)=0.3,
% \]
% and $T=1$.
 \item Parameters: 
    \[
A_t = 
\begin{bmatrix}
0.1 & -0.05\\
-0.05 & 0.1
\end{bmatrix},\quad
B_t^1 = 
\begin{bmatrix}
0.04 \\
0.03
\end{bmatrix},\quad
B_t^2 = 
\begin{bmatrix}
0.01 \\
-0.05
\end{bmatrix},\quad
W = 
\begin{bmatrix}
0.2 & 0.05 \\
0.05 & 0.1
\end{bmatrix},
\]
    \[
Q_T^1=Q_T^2 =Q_t^1=Q_t^2 = 
\begin{bmatrix}
0.1 & -0.01\\
-0.01 & 0.1
\end{bmatrix},\quad
R_t^1(t)=R_t^2(t)=0.35,
\]
and $T=2$.
    \item Initialization: Take $x_0=(x_0^1,x_0^2)$ where $x_0^1, x_0^2$ are independent and sampled from  $\mathcal{N}(0.25,0.2)$ and $\mathcal{N}(0.4,0.3)$ respectively. The initial policies are $\pmb{K}^{1,(0)}=\pmb{K}^{2,(0)}=(0.2,0.01)$.
\end{itemize}
\item A three-player LQ game example where Assumption \ref{ass:noise_cond} is satisfied:
    \begin{itemize}
    \item Parameters: 
    \[
A_t = 
\begin{bmatrix}
0.05 & -0.1 & 0.1\\
0.1 & 0.2 & -0.06\\
-0.02 & 0.03 &0.1
\end{bmatrix},\quad
B_t^1 = 
\begin{bmatrix}
0.05 \\
0.01 \\
-0.01
\end{bmatrix},\quad
B_t^2 = 
\begin{bmatrix}
0.01 \\
-0.05\\
-0.02
\end{bmatrix},\quad
B_t^3 = 
\begin{bmatrix}
-0.02 \\
0.01\\
0.05
\end{bmatrix},
\]
\[
W = 
\begin{bmatrix}
0.1 & 0.01 & 0.02 \\
0.01 & 0.2 & 0.01\\
0.02 & 0.01 &  0.1
\end{bmatrix},\quad
Q_T^1=Q_T^2 =Q_T^3=Q_t^1=Q_t^2 = Q_t^3= 
\begin{bmatrix}
0.2 & 0 & 0\\
0 & 0.2 & 0 \\
0 & 0 & 0.2 \\
\end{bmatrix},
\]
$R_t^1(t)=R_t^2(t)=0.5$, $R_t^3(t)=0.6$, and $T=1$.
    \item Initialization: Take $x_0=(x_0^1,x_0^2,x_0^3)$ where $x_0^1, x_0^2, x_0^3$ are independent and sampled from  $\mathcal{N}(0.3,0.2)$ and $\mathcal{N}(0.2,0.3)$, and $\mathcal{N}(0.3,0.2)$ respectively. The initial policies are $\pmb{K}^{1,(0)}=(0,-0.01,0)$, $\pmb{K}^{2,(0)}=(0,-0.01,0.01)$, and $\pmb{K}^{3,(0)}=(0,-0.001,0)$.
\end{itemize}

\item In most large population games, the individual contribution to the joint dynamics scales like $\frac{1}{N}$ \cite{huang2006large,lasry2007mean}. In this setting, we denote by $x_t^{(N)}$, the state dynamics with $N$ agents, which follows 
\begin{eqnarray}
x_{t+1}^{(N)} = A^{(N)}_t x^{(N)}_t + \sum_{i=1}^N B^{(N),i}_t u_t^{(N),i} +w^{(N)}_t.
\end{eqnarray}
If $B^{(N),i}_t = \mathcal{O}(\frac{1}{N})$, in the sense that there exists a constant $C_B>0$ (which does not depend on $N$) and $N_0 \in \mathbb{N}_+$ such that
\begin{eqnarray}
\left\|B^{(N),i}_t\right\| \leq \frac{C_B}{N}
\end{eqnarray}
for all $0\leq t \leq T$ and $N \ge N_0$, then we have
\begin{eqnarray}
\gamma_B^{(N)}: = \max_i \Big\{\max_{t=0,1,\cdots,T-1}\|B^{(N),i}_t\|\Big\} \leq \frac{C_B}{N}
\end{eqnarray}
for $N \ge N_0$. In this case, Assumption \ref{ass:noise_cond} is reduced to the following condition:
\begin{eqnarray}
\frac{(\underline{\sigma}_{\pmb{X}})^5}{\|\Sigma_{\pmb{K}^*}\|} > 20\frac{(N-1)^2}{N^4} T^2 d \frac{(C_B)^4 (\max_i \{C^i(\pmb{K}^*)\}+\psi)^4}{\underline{\sigma}_{\pmb{Q}}^2\underline{\sigma}_{\pmb{R}}^2} \left( \frac{\bar{\rho}^{2T}-1}{\bar{\rho}^{2}-1}\right)^2 = \mathcal{O} \left(\frac{1}{N^2}\right).\label{eq:noise}
\end{eqnarray}
The RHS of condition \eqref{eq:noise} decays quadratically in $N$ since one can check that $C^i(\pmb{K}^*)$, $\psi$ and $\bar{\rho}$ do not grow with respect to $N$ in this case. Hence  condition \eqref{eq:noise} is automatically satisfied in the large population regime.
\item In addition, if we focus on the set of stabilizing policies
\begin{eqnarray}
\Omega:=\{\pmb{K}: \|A_t-B_tK_t\|\leq \rho<1,\quad \forall t=0,1,\cdots,T-1\},
\end{eqnarray}
then the term $\left( \frac{\bar{\rho}^{2T}-1}{\bar{\rho}^{2}-1}\right)^2$ on the RHS of \eqref{eq:system_noise} can be replaced by $\left( \frac{1}{1-\rho}\right)^2$, which is independent of the horizon $T$.

\item Furthermore, $T^2$  on the RHS of \eqref{eq:system_noise} can be removed when the dynamics and the cost parameters are time independent, and the agent is searching for time-independent linear policies, namely 
\begin{eqnarray*}
R^i_t \equiv R^i, Q^i_t \equiv Q^i, A^i_t \equiv A^i, B_t^i\equiv B_,  K^i_t \equiv K^i, \quad t=0,1,\cdots,T-1, i=1,2,\cdots,N.
\end{eqnarray*}
\end{enumerate}
 }
\end{Remark}

% \begin{Remark}%[Assumption \ref{ass:noise_cond} on System Noise]
% \label{app:sys_noise_ass}{\rm
% An example where Assumption \ref{ass:noise_cond} is satisfied:
%     \begin{itemize}
%     \item Parameters: 
%     \[
% A_t = 
% \begin{bmatrix}
% 0.2 & 0\\
% 0 & 0.2
% \end{bmatrix},\quad
% B_t^1 = 
% \begin{bmatrix}
% 0.1 \\
% 0
% \end{bmatrix},\quad
% B_t^2 = 
% \begin{bmatrix}
% 0 \\
% -0.1
% \end{bmatrix},\quad
% W = 
% \begin{bmatrix}
% 0.2 & 0 \\
% 0 & 0.2
% \end{bmatrix},
% \]
%     \[
% Q_T^1=Q_T^2 =Q_t^1=Q_t^2 = 
% \begin{bmatrix}
% 0.1 & 0\\
% 0 & 0.1
% \end{bmatrix},\quad
% R_t^1(t)=R_t^2(t)=0.3,
% \]
% and $T=1$.
%     \item Initialization: Take $x_0=(x_0^1,x_0^2)$ where $x_0^1, x_0^2$ are independent and sampled from  $\mathcal{N}(0.3,0.2)$ and $\mathcal{N}(0.4,0.3)$ respectively. The initial policies are $\pmb{K}^{1,(0)}=\pmb{K}^{2,(0)}=(0.2,0.01)$.
% \end{itemize} }
% \end{Remark}

\begin{Theorem}[Global Convergence of the Natural Policy Gradient Method]\label{thm:local_conv_NPG}
Assume Assumptions~\ref{ass:nonzero_cost}, \ref{ass:nonzero_initial_noise}, \ref{ass:nonzero_exist_sol} hold. We also assume Assumption \ref{ass:noise_cond} so that 
\begin{equation}\label{eqn:defn_hat_alpha}
\widehat{\alpha} := \frac{\sx\srmin}{\|\Sigma_{\pmb{K}^{*}}\|}-20\,(N-1)^2\,T^2\,d\, \frac{(\gamma_B)^4(\max_i\{C^i(\pmb{K}^*)\}+\psi)^4}{\sx^4\sqmin^2\srmin}\left(\frac{\bar{\rho}^{2T}-1}{\bar{\rho}^2-1}\right)^2 >0.
\end{equation}
Then there exists an $0<\eta_0\in\HH(\frac{1}{\sum_{i=1}^NC^i(\pmb{K}^{i,(0)},\pmb{K}^{-i*})+1})$ such that for $\epsilon>0, 0<\eta<\eta_0$,
%{\color{red}If we choose $\eta$ such that $\eta\in (0,\frac{1}{\widehat{\alpha}})$},
the cost function of the natural gradient method \eqref{eqn:nonzero_grad_update_rule} satisfies
\[
\sum_{i=1}^N\Big( C^i(\pmb{K}^{i,(M)},\pmb{K}^{-i*})-C^i(\pmb{K}^{*}) \Big)\leq\epsilon,
\]
whenever $M\geq \frac{1}{\widehat{\alpha}\eta}\log(\frac{\sum_{i=1}^N(C^i(\pmb{K}^{i,(0)},\pmb{K}^{-i*})-C^{i}(\pmb{K}^*))}{\epsilon}).$
\end{Theorem}
The convergence rate developed in Theorem \ref{thm:local_conv_NPG} is usually referred to as a linear convergence  rate \cite{nocedal2006numerical} as we show that $\frac{\sum_{i=1}^N (C^i(\pmb{K}^{i,(m+1)},\pmb{K}^{-i*})-C^i(\pmb{K}^{*})) }{\sum_{i=1}^N ( C^i(\pmb{K}^{i,(m)},\pmb{K}^{-i*})-C^i(\pmb{K}^{*}))}\leq r $ for some $r\in (0,1)$. This leads to the result, seen in point 2. of the Theorem, that $\sum_{i=1}^N \left(C^i(\pmb{K}^{i,(M)},\pmb{K}^{-i*})-C^i(\pmb{K}^{*})\right) = O(e^{-M})$ and therefore this is also called an exponential 
convergence rate in the literature \cite{nocedal2006numerical}.

Compared to directly evaluating $C^i(\pmb{K}^{(m)})$ along the training process $m=1,2,\cdots,M$, we observe that it is more natural to analyze the cost $C^i(\pmb{K}^{i,(m)},\pmb{K}^{-i*})$, with player $i$ taking policy $\pmb{K}^{i,(m)}$ and (as if) the other players were using optimal policies, as  conditional optimality is the essence of a Nash equilibrium. This view point is the key to establishing the global convergence result and this criterion is also used in all the key lemmas as well as the proof of the main theorem.

Note that our result can be easily generalized to the setting with agent-dependent learning rates as long as $\eta_i<\eta_0$ holds for each agent's step size $\eta_i$.

The proof of Theorem \ref{thm:local_conv_NPG} relies on the regularity of the LQ game problem, some properties of the gradient descent dynamics, and the perturbation analysis of the covariance matrix of the controlled dynamics.

\paragraph{Vanilla Policy Gradient Method.} By modifying parts of the proof of the convergence result for the natural policy gradient method, we can extend the convergence result to the case of the vanilla policy gradient method, where the updating rule for each player $i$ ($i=1,\cdots,N$), in a sequence of games for $m=1,\dots, M$, is given by
\begin{equation}\label{eqn:nonzero_grad_update_rule_vani}
    K_{t}^{i,(m)} = K_t^{i,(m-1)} - \eta\nabla_{K_t^i}C^i(\pmb{K}^{(m-1)}),\quad\forall\, 0\leq t\leq T-1.
\end{equation}

We now give the system noise assumption and the global convergence result for the vanilla policy gradient method. %in Assumption \ref{ass:noise_cond_vani} and Theorem \ref{thm:conv_PG_vani}.

\begin{Assumption}[System Noise for the Vanilla Policy Gradient]\label{ass:noise_cond_vani} The system parameters satisfy the following two inequalities
\begin{equation}\label{eq:system_noise_vani_1}
\frac{(\sx)^7}{\|\Sigma_{\pmb{K}^{*}}\| } >20(N-1)^2\,T^2\,d\, \frac{(\gamma_B)^4(\max_i\{C^i(\pmb{K}^*)\}+\psi)^4}{\sqmin^2\srmin^2}\left(\frac{\bar{\rho}^{2T}-1}{\bar{\rho}^2-1}\right)^2\left(\bar{\rho}^{2T}\|\Sigma_0\|+(\bar{\rho}^{2T}+1)\|W\|\right)^2,
\end{equation}
and
\begin{equation}\label{eq:system_noise_vani_2}
(\sx)^2 > \frac{5(\max_i\{C^i(\pmb{K}^*)\}+\psi)}{\sqmin}\left(\frac{\max_i\{C^i(\pmb{K}^*)\}+\psi}{\sqmin}+\bar{\rho}^{2T}\|\Sigma_0\|+(\bar{\rho}^{2T}+1)\|W\|\right).
\end{equation}
\end{Assumption}

\begin{Theorem}[Global Convergence of the Vanilla Policy Gradient Method]\label{thm:conv_PG_vani}
Assume Assumptions~\ref{ass:nonzero_cost}, \ref{ass:nonzero_initial_noise}, \ref{ass:nonzero_exist_sol} hold. We also assume Assumption \ref{ass:noise_cond_vani} so that 
\begin{eqnarray}\label{eqn:defn_hat_alpha_vani}
\widetilde{\alpha}\hspace{-0.2cm} &:=&\hspace{-0.2cm} \frac{\sx^2\srmin}{\|\Sigma_{\pmb{K}^{*}}\|}-20\,(N-1)^2\,T^2\,d\, \frac{(\gamma_B)^4(\max_i\{C^i(\pmb{K}^*)\}+\psi)^4}{\sx^5\sqmin^2\srmin}\left(\frac{\bar{\rho}^{2T}-1}{\bar{\rho}^2-1}\right)^2\hspace{-0.2cm}\left(\bar{\rho}^{2T}\|\Sigma_0\|+(\bar{\rho}^{2T}+1)\|W\|\right)^2\nonumber\\
&>&\hspace{-0.2cm} 0.
\end{eqnarray}
Then there exists an $0<\widetilde{\eta}_0\in\HH(\frac{1}{\sum_{i=1}^NC^i(\pmb{K}^{i,(0)},\pmb{K}^{-i*})+1})$ such that for $\epsilon>0, 0<\eta<\widetilde{\eta}_0$,
%{\color{red}If we choose $\eta$ such that $\eta\in (0,\frac{1}{\widehat{\alpha}})$},
the cost function of the vanilla policy gradient method \eqref{eqn:nonzero_grad_update_rule_vani} satisfies
\[
\sum_{i=1}^N\Big( C^i(\pmb{K}^{i,(M)},\pmb{K}^{-i*})-C^i(\pmb{K}^{*}) \Big)\leq\epsilon,
\]
whenever $M\geq \frac{1}{\widetilde{\alpha}\eta}\log(\frac{\sum_{i=1}^N(C^i(\pmb{K}^{i,(0)},\pmb{K}^{-i*})-C^{i}(\pmb{K}^*))}{\epsilon}).$
\end{Theorem}

Since the key step to prove the global convergence of the natural policy gradient method (Theorem \ref{thm:local_conv_NPG}) is the one-step contraction lemma (Lemma \ref{lemma:one_step_loc_conv}), we also establish the one-step contraction result for the vanilla version by modifying some parts of Lemma \ref{lemma:one_step_loc_conv}, see Lemma \ref{lemma:one_step_cont_vani} and its proof in Appendix \ref{app:vani}. To prove Theorem \ref{thm:conv_PG_vani}, it suffices to show that Lemma \ref{lemma:one_step_cont_vani} holds, and the rest of the proof follows the same arguments as in the proof of Theorem \ref{thm:local_conv_NPG} for the natural version.

We explain the differences between the system noise assumptions required for the vanilla and the natural policy gradient methods in Remark \ref{remark:comp_noi_cond} in Appendix \ref{app:vani}.

\subsection{Regularity of the LQ game and Properties of the Gradient Descent Dynamics}

We begin with the analysis of some properties of the N-player general-sum LQ game \eqref{eqn:nonzero_State}-\eqref{eqn:nonzero_cost_at_t}. Our aim is to establish two key results Lemma \ref{lemma:grad_domi} and Lemma \ref{lemma:nonzero_almost_smoothness} which provide the gradient dominance condition and a smoothness condition on the cost function $C^i(\pmb{K})$ of player $i$ with respect to the joint policy $\pmb{K}$ from all players, respectively.

%To start we show in Lemma \ref{lemma:PK1K2_posi_defi} that the Ricatti system  $\{P^{\pmb{K}}_{t,i}\}_{t=0}^T$ for each player $i$, which is induced by a joint policy $\pmb{K}$ from all $N$ players is well defined;  Lemma \ref{lemma:nonzero_policygrad} gives a representation of the gradient term and Lemma \ref{lemma:nonzero_cost_diff} gives a representation of the cost difference using advantage functions; Lemma \ref{lemma:grad_domi} and Lemma \ref{lemma:nonzero_almost_smoothness} provide the gradient dominance condition and a smoothness condition on the cost function $C^i(\pmb{K})$ of player $i$ with respect to the joint policy $\pmb{K}$ from all players, respectively; and finally, Lemma \ref{lemma:nonzero_bds_P_Sigma} gives two useful upper bounds on the Ricatti system and state covariance matrices.

%Let us start with the analysis of some properties of the $N$-player general-sum LQG \eqref{eqn:nonzero_State}-\eqref{eqn:nonzero_cost_at_t}. 

In the finite horizon setting, define $P_{t,i}^{\pmb{K}}$ as the solution to 
\begin{equation}\label{eqn:lyapunov}
    P_{t,i}^{\pmb{K}} = Q_t^i + (K_t^{i})^\top R_t^{i} K_t^{i} + \left(A_t - \sum_{j=1}^NB_t^jK_t^{j} \right)^\top P_{t+1,i}^{\pmb{K}}\left(A_t - \sum_{j=1}^NB_t^jK_t^{j} \right),
\end{equation}
with terminal condition
\[
P_{T,i}^{\pmb{K}} = Q_{T}^i.
\]

\begin{Lemma}\label{lemma:PK1K2_posi_defi}
Assume Assumptions \ref{ass:nonzero_cost} holds. Then for $i=1,\cdots,N$, $t=0,\dots,T$ the matrices $P_{t,i}^{\pmb{K}}$ defined in \eqref{eqn:lyapunov} are positive definite.
\end{Lemma}

\begin{proof}
We prove that the terms in the sequence $\{P_{t,i}^{\pmb{K}}\}_{t=0}^T$ are positive definite for $i=1,\cdots,N$ by backward induction. For $t=T$,   $P_{T,i}^{\pmb{K}} = Q_T^i$ is positive definite since $Q_T^i$ is positive definite. Assume $P_{t+1,i}^{\pmb{K}}$ is positive definite for some $t+1$, then take any $z \in \mathbb{R}^d$ such that $z \neq 0$,
\begin{equation*}
    \begin{split}
        z^{\top} P_{t,i}^{\pmb{K}} z &= z^\top Q_t^i z + z^\top(K_t^{i})^\top R_t^{i} K_t^{i}z  
        + z^\top\left(A_t - \sum_{j=1}^N B_t^jK_t^{j}\right)^\top P_{t+1,i}^{\pmb{K}}\left(A_t - \sum_{j=1}^N B_t^jK_t^{j}\right)z > 0.
    \end{split}
\end{equation*}
The last inequality holds since $z^\top Q_t^i z>0$ and the other terms are non-negative under Assumption \ref{ass:nonzero_cost}. By backward induction, we have $P_{t,i}^{\pmb{K}}$ positive definite, $\forall\,t=0,1,\cdots,T$.
\end{proof}

Player $i$'s cost under the set of policies $\pmb{K}$ can be rewritten as
\[
C^i(\pmb{K}) = \mathbb{E}\left[x_0^\top P_{0,i}^{\pmb{K}} x_0+N_{0,i}^{\pmb{K}}\right],
\]
where the expectation is taking with respect to the initial state $x_0$ and the $ N_{t,i}^{\pmb{K}}$ are defined backwards:
\begin{equation}\label{eqn:defn_Nt_nonzero}
    N_{t,i}^{\pmb{K}} = N_{t+1,i}^{\pmb{K}} +\Tr(WP_{t+1,i}^{\pmb{K}}),\quad N_{T,i}^{\pmb{K}}=0.
\end{equation}
To see this,
\begin{eqnarray*}
     \mathbb{E}\left[x_0^\top P_{0,i}^{\pmb{K}}x_0+N_{0,i}^{\pmb{K}}\right] &=& \mathbb{E}\left[x_0^\top Q_0^i x_0+ x_0^\top(K_0^{i})^\top R_0^{i} K_0^{i} x_0  + \sum_{t=0}^{T-1}w_t^\top P_{t+1,i}^{\pmb{K}} w_t\right.\\
     && \quad + \left.x_0^\top(A_0 - \sum_{j=1}^NB_0^jK_0^{j})^\top P_{1,i}^{\pmb{K}}(A_0 - \sum_{j=1}^NB_0^jK_0^{j})x_0 \right]\\
     &=& \mathbb{E}\left[x_0^\top Q_0^i x_0 + (u_0^i)^\top R_0^{i}u_0^i+  x_1^\top P_{1,i}^{\pmb{K}}x_1 + \sum_{t=1}^{T-1}w_t^\top P_{t+1,i}^{\pmb{K}} w_t\right]\\
     &=& \mathbb{E}\left[c_0^i(x_0,u_0^i) + x_1^\top P_{1,i}^{\pmb{K}}x_1 + \sum_{t=1}^{T-1}w_t^\top P_{t+1,i}^{\pmb{K}} w_t\right]\\
     &=& \mathbb{E}\left[c_0^i(x_0,u_0^i) + c_1^i(x_1,u_1^i) +x_2^\top P_{2,i}^{\pmb{K}}x_2 + \sum_{t=2}^{T-1}w_t^\top P_{t+1,i}^{\pmb{K}} w_t\right]=\mathbb{E}\left[\sum_{t=0}^{T}c_t^i(x_t,u_t^i)\right].
\end{eqnarray*}
In addition, for $t=0,1,\cdots,T-1$ and $i=1,\cdots,N$, define 
\begin{equation}\label{eqn:defn_Ei}
    E_{t,i}^{\pmb{K}} = R_t^{i}K_t^{i}-(B_t^i)^\top P_{t+1,i}^{\pmb{K}}\left(A_t-\sum_{j=1}^NB_t^jK_t^{j}\right).
\end{equation}
Then we have the following representation of the gradient terms.

\begin{Lemma}\label{lemma:nonzero_policygrad}
The policy gradients have the following representation: for $t=0,1,\cdots,T-1$,
\begin{equation}\label{eqn:poli_grad_exp}
    \nabla_{K_t^i}C(\pmb{K}) = 2\left(R_t^{i}K_t^{i}-(B_t^i)^\top P_{t+1,i}^{\pmb{K}}\left(A_t-\sum_{j=1}^NB_t^jK_t^{j}\right)\right)\Sigma_t^{\pmb{K}} = 2E_{t,i}^{\pmb{K}}\,\Sigma_t^{\pmb{K}},
\end{equation}
where $E_{t,i}^{\pmb{K}}$ is defined in \eqref{eqn:defn_Ei}.
\end{Lemma}

\begin{proof}
Expressing the cost function in terms of $K^i_t$ and suppressing the arguments of the cost $c^i_t$,
\begin{eqnarray*}
C^i(\pmb{K}) &=& \mathbb{E}\left[\sum_{s=0}^{t-1}c_s^i + c_t^i+x_{t+1}^\top P_{t+1,i}^{\pmb{K}}x_{t+1}+\sum_{s=t+1}^{T-1}w_s^\top P_{s+1,i}^{\pmb{K}} w_s\right]\\
&=& \mathbb{E}\left[\sum_{s=0}^{t-1}c_s^i + x_t^\top(Q_t^i + (K_t^{i})^\top R_t^{i} K_t^{i})x_t+\sum_{s=t+1}^{T-1}w_s^\top P_{s+1,i}^{\pmb{K}} w_s\right.\\
&& +\left.\left(x_t^\top\left(A_t-\sum_{j=1}^NB_t^jK_t^{j}\right)^\top +w_t^\top\right) P_{t+1,i}^{\pmb{K}}\left(\left(A_t-\sum_{j=1}^NB_t^jK_t^{j}\right)x_t+w_t\right)\right].
\end{eqnarray*}
Therefore, for $i=1,\dots,N$, the gradients are given by
\begin{equation*}
    \nabla_{K_t^{i}}C^i(\pmb{K}) = 2\left(R_t^{i}K_t^{i}-(B_t^i)^\top P_{t+1,i}^{\pmb{K}}\left(A_t-\sum_{j=1}^NB_t^jK_t^{j}\right)\right)\Sigma_t^{\pmb{K}}.
\end{equation*}
\end{proof}

We can write the value function for player $i$, $V_{\pmb{K}}^i(x,t)$ for $t=0,1,\cdots,T-1$, as 
\begin{equation*}
    V_{\pmb{K}}^i(x,t) = \mathbb{E}_{\pmb{w}}\left.\left[\sum_{s=t}^{T-1}(x_s^\top Q_s^ix_s+(u_s^i)^\top R_s^{i} u_s^i)+x_T^{\top}Q_T^ix_T\right|x_t = x\right]=x^{\top}P_{t,i}^{\pmb{K}} x +N_{t,i}^{\pmb{K}},
\end{equation*}
 with terminal condition 
\begin{equation}\label{eqn:terminal_V}
    V_{\pmb{K}}^i(x,T)=x^{\top}Q_T^ix, 
\end{equation}
where $N_{t,i}^{\pmb{K}}$ is defined in \eqref{eqn:defn_Nt_nonzero}, and $\mathbb{E}_{\pmb{w}}$ denotes the expectation over the noise $\pmb{w}$.
We then define the $Q$ function, $Q_{\pmb{K}}^i(x,u,t)$ for the Markovian control $u=(u^1,\cdots,u^N)$ to be
\begin{equation*}
    Q_{\pmb{K}}^i(x,u,t)=x^\top Q_t^ix+(u^i)^\top R_t^{i} u^i +  \mathbb{E}_{w_{t}}\left[ V_{\pmb{K}}^i\Big(A_tx+\sum_{j=1}^NB_t^ju^j+w_{t},t +1\Big)\right],
\end{equation*}
and the advantage function to be
\begin{equation*}
    A_{\pmb{K}}^i(x,u,t) = Q_{\pmb{K}}^i(x,u,t)-V_{\pmb{K}}^i(x,t).
\end{equation*}
Here $\mathbb{E}_{w_t}$ denotes the expectation taken with respect to $w_t$.
Note that $C^i(\pmb{K}) = \mathbb{E}[V_{\pmb{K}}^i(x_0,0)]$. We write
\begin{equation*}
    (\pmb{K}^{i\prime},\pmb{K}^{-i}):=(\pmb{K}^1,\cdots,\pmb{K}^{i-1},\pmb{K}^{i\prime},\pmb{K}^{i+1},\cdots,\pmb{K}^N),
\end{equation*}
and the control sequences under the policy $(\pmb{K}^{i\prime},\pmb{K}^{-i})$ as $\{u_t^{i\prime,-i}\}_{t=0}^{T-1}$ with
\begin{equation}\label{eqn:defn_pert_u}
   u_t^{i\prime}=-K_t^{i\prime}x_t^{\pmb{K}^{i\prime},\pmb{K}^{-i}} ,\quad\text{and}\quad u_t^j=-K_t^jx_t^{\pmb{K}^{i\prime},\pmb{K}^{-i}},\quad j\neq i,
\end{equation}
where $\{x_t^{\pmb{K}^{i\prime},\pmb{K}^{-i}}\}_{t=0}^T$ is the state trajectory under $(\pmb{K}^{i\prime},\pmb{K}^{-i})$. Then we can write the difference between the cost functions of $\pmb{K}=(\pmb{K}^1,\cdots,\pmb{K}^N)$ and $(\pmb{K}^{i\prime},\pmb{K}^{-i})$ in terms of advantage functions.

\begin{Lemma}[Cost Difference]\label{lemma:nonzero_cost_diff}
Assume $\pmb{K}$ and $(\pmb{K}^{i\prime},\pmb{K}^{-i})$ have finite costs. Then
\begin{equation}
    V_{\pmb{K}^{i\prime},\pmb{K}^{-i}}^i(x,0)-V_{\pmb{K}}^i(x,0) = \mathbb{E}_{\pmb{w}}\left[\sum_{t=0}^{T-1} A_{\pmb{K}}^i(x_t^{\pmb{K}^{i\prime},\pmb{K}^{-i}},u_t^{i\prime,-i},t)\right],
\end{equation}
where $\{u_t^{i\prime,-i}\}_{t=0}^{T-1}$ is defined in \eqref{eqn:defn_pert_u} with $x_0^{\pmb{K}^{i\prime},\pmb{K}^{-i}}=x$.
For player $i$,
\begin{eqnarray*}
     A_{\pmb{K}}^i(x_t^{\pmb{K}^{i\prime},\pmb{K}^{-i}},u_t^{i\prime,-i},t)
    &=& (x_t^{\pmb{K}^{i\prime},\pmb{K}^{-i}})^{\top}(K_t^{i\prime}-K_t^i)^{\top}(R_t^{i}+(B_t^i)^{\top}P_{t +1,i}^{\pmb{K}}B_t^i)(K_t^{i\prime}-K_t^i)\,x_t^{\pmb{K}^{i\prime},\pmb{K}^{-i}}\\ &&+\,2(x_t^{\pmb{K}^{i\prime},\pmb{K}^{-i}})^{\top}(K_t^{i\prime}-K_t^i)^{\top}E_{t,i}^{\pmb{K}}\,x_t^{\pmb{K}^{i\prime},\pmb{K}^{-i}}.
\end{eqnarray*}
\end{Lemma}

\begin{proof} Write $x_t^\prime=x_t^{\pmb{K}^{i\prime},\pmb{K}^{-i}}$ for simplicity and denote by $c_t^{i\prime}(x)$ the (instantaneous) cost of player $i$ generated by $(\pmb{K}^{i\prime},\pmb{K}^{-i})$ with a single trajectory starting from $x^\prime_0=x$. That is,
\[
    c_t^{i\prime}(x_t^\prime) =(x_t^\prime)^\top Q_t^ix_t^\prime+(K_t^{i\prime}x_t^{i\prime})^\top R_t^{i}K_t^{i\prime}x_t^{i\prime}, \quad t = 0,1,\cdots,T-1,
\]
and
\begin{equation*}
    c_T^{i\prime}(x_T^\prime) = (x_T^{\prime})^\top Q_T^ix_T^{\prime},
\end{equation*}
with 
\[
u_{t}^{i\prime} = - K_t^{i\prime} x_t^\prime\quad, x_{t+1}^\prime = A_tx_t^\prime +B_t^iu_t^{i\prime} +\sum_{j=1,j\neq i}^NB_t^ju_t^{j} + w_t,\,\,\,\, x_0^\prime=x.
\]
Therefore
\begin{eqnarray*}
  V_{\pmb{K}^{i\prime},\pmb{K}^{-i}}^i(x,0)-V_{\pmb{K}}^i(x,0)  &=& \mathbb{E}_{\pmb{w}}\left[\sum_{t=0}^{T} c_t^{i\prime}(x_t^\prime)\right] -V_{\pmb{K}}^i(x,0)\\
     & =& \mathbb{E}_{\pmb{w}}\left[ \sum_{t=0}^{T} \left(c_t^{i\prime}(x_t^\prime)-V_{\pmb{K}}^i(x_t',t)\right) +\sum_{t=1}^T V_{\pmb{K}}^i(x_t',t)\right]\\
     & = & \mathbb{E}_{\pmb{w}}\left[\sum_{t=0}^{T-1}\left( c_t^{i\prime}(x_t^\prime)+V_{\pmb{K}}^i(x_{t+1}',t+1)-V_{\pmb{K}}^i(x_t',t)\right)\right]\\
     & = &\mathbb{E}_{\pmb{w}}\left.\left[ \sum_{t=0}^{T-1}\left( Q_{\pmb{K}}^i(x_t^\prime,u_t^{i\prime,-i},t)-V_{\pmb{K}}^i(x_t^\prime,t)\right)\right|x_0=x\right] \\ &=& \mathbb{E}_{\pmb{w}}\left.\left[\sum_{t=0}^{T-1} A_{\pmb{K}}^i(x_t^\prime,u_t^{i\prime,-i},t)\right|x_0=x\right],
\end{eqnarray*}
where the third equality holds since $c_T^{i\prime}(x_T^\prime)=V_{\pmb{K}}^i(x_T^{\prime},T)$ by \eqref{eqn:terminal_V} with the same single trajectory.

For player $i$, the advantage function is given by
\begin{eqnarray*}
    && A_{\pmb{K}}^i(x_t^\prime,u_t^{i\prime,-i},t) 
    = Q_{\pmb{K}}^i(x_t^\prime,u_t^{i\prime,-i},t)-V_{\pmb{K}}^i(x_t^\prime,t) \\
    & =& (x_t^\prime)^{\top}(Q_t^i+(K_t^{i\prime})^{\top}R_t^{i} K_t^{i\prime})x_t^\prime + \mathbb{E}_{w_t} \left[V_{\pmb{K}}^i((A_t-B_t^iK_t^{i\prime}-\sum_{j=1,j\neq i}^NB_t^jK_t^{j})x_t^\prime+w_{t},t +1)\right]-V_{\pmb{K}}^i(x_t^\prime,t) \\
    &=& (x_t^\prime)^{\top}(Q_t^i+(K_t^{i\prime})^{\top}R_t^{i} K_t^{i\prime})x_t^\prime + \left((x_t^\prime)^{\top}\left(A_t-B_t^iK_t^{i\prime}-\sum_{j=1,j\neq i}^NB_t^jK_t^{j}\right)^{\top} \right.\\
    && \left.P_{t+1,i}^{\pmb{K}}\left(A_t-B_t^iK_t^{i\prime}-\sum_{j=1,j\neq i}^NB_t^jK_t^{j}\right)x_t^\prime+\Tr(WP_{t+1,i}^{\pmb{K}})+N_{t+1,i}^{\pmb{K}}\right) - \left( (x_t^\prime)^{\top}P_{t,i}^{\pmb{K}}x_t^\prime + N_{t,i}^{\pmb{K}}\right)\\
    &=& (x_t^\prime)^{\top}(Q_t^i+(K_t^{i\prime}-K_t^i+K_t^i)^{\top}R_t^{i}(K_t^{i\prime}-K_t^i+K_t^i))x_t^{\prime} \\
    && + (x_t^\prime)^{\top}\left(A_t- B_t^i(K_t^{i\prime}-K_t^i) - \sum_{j=1}^NB_t^jK_t^{j}\right)^{\top}P_{t +1,i}^{\pmb{K}}\left(A_t- B_t^i(K_t^{i\prime}-K_t^i) - \sum_{j=1}^NB_t^jK_t^{j}\right)x_t^\prime \\
    && - (x_t^\prime)^{\top}\left(Q_t^i+(K_t^{i})^{\top}R_t^{i} K_t^{i}+\left(A_t-\sum_{j=1}^NB_t^jK_t^{j}\right)^{\top}P_{t+1,i}^{\pmb{K}}\left(A_t-\sum_{j=1}^NB_t^jK_t^{j}\right)\right)x_t^{\prime} \\
    &=& (x_t^\prime)^{\top}(K_t^{i\prime}-K_t^i)^{\top}(R_t^{i}+(B_t^i)^{\top}P_{t +1,i}^{\pmb{K}}B_t^i)(K_t^{i\prime}-K_t^i)x_t^\prime +\,2(x_t^{\prime})^{\top}(K_t^{i\prime}-K_t^i)^{\top}E_{t,i}^{\pmb{K}}\,x_t^\prime.
\end{eqnarray*}

\end{proof}

Note that the derivations of Lemmas \ref{lemma:PK1K2_posi_defi},
\ref{lemma:nonzero_policygrad} and \ref{lemma:nonzero_cost_diff} are largely inspired by \cite{FGKM2018,hambly2020policy}. However, the final expressions are different since our setting is different from both \cite{FGKM2018} and \cite{hambly2020policy}. For completeness, we provide the proofs and derivations in our context.

For the policy gradient method \cite{FGKM2018,hambly2020policy} in the single-agent setting, gradient domination and smoothness of the objective function are two key conditions to guarantee the global convergence of the gradient descent methods. This is also the case for the N-player game setting. The gradient dominance condition for each player $i$ is proved in Lemma \ref{lemma:grad_domi}, which indicates that for a policy $\pmb{K}$, the
distance between $C^i(\pmb{K})$ and the optimal cost $C^i(\pmb{K}^*)$ is bounded by the sum of the magnitudes of the gradients $\nabla_tC^i(\pmb{K})$ for $t = 0, 1, \cdots , T-1$. The smoothness condition for each player $i$ is proved in Lemma \ref{lemma:nonzero_almost_smoothness} where the difference between $C^i(\pmb{K}^{i\prime},\pmb{K}^{-i})$ and $ C^i(\pmb{K})$ can be rewritten as a function of $\pmb{K}^{i\prime}-\pmb{K}^{i}$.

\begin{Lemma}[Gradient Dominance]\label{lemma:grad_domi}
Assume Assumptions \ref{ass:nonzero_cost}, \ref{ass:nonzero_initial_noise}, and \ref{ass:nonzero_exist_sol} hold. Then, for player $i$, we have
\begin{eqnarray*}
     C^i(\pmb{K}^i,\pmb{K}^{-i*})-C^i(\pmb{K}^{*})
     &\leq& \frac{\|\Sigma_{\pmb{K}^{*}}\|}{\srmin}\sum_{t=0}^{T-1}\Tr\Big((E_{t,i}^{\pmb{K}^i,\pmb{K}^{-i*}})^\top E_{t,i}^{\pmb{K}^i,\pmb{K}^{-i*}}\Big)\\
     &\leq& \frac{\|\Sigma_{\pmb{K}^{*}}\|}{4\srmin\,(\sx)^2}\sum_{t=0}^{T-1}\Tr\Big(\nabla_{K_t^i}C^i(\pmb{K}^i,\pmb{K}^{-i*})^\top\nabla_{K_t^i}C^i(\pmb{K}^i,\pmb{K}^{-i*})\Big),\\
\end{eqnarray*}
and
\begin{eqnarray*}
     C^i(\pmb{K}^i,\pmb{K}^{-i*})-C^i(\pmb{K}^{*}) 
     &\geq& \sx\sum_{t=0}^{T-1}\frac{1}{\|R_t^{i}+(B_t^i)^{\top}P_{t+1,i}^{\pmb{K}^i,\pmb{K}^{-i*}}B_t^i\|}\Tr\big((E_{t,i}^{\pmb{K}^i,\pmb{K}^{-i*}})^\top E_{t,i}^{\pmb{K}^i,\pmb{K}^{-i*}}\big)\\
     &\geq& \frac{\sx}{4\|\Sigma_{\pmb{K}^i,\pmb{K}^{-i*}}\|^2}\sum_{t=0}^{T-1}\frac{\Tr\big(\nabla_{K_t^i}C^i(\pmb{K}^i,\pmb{K}^{-i*})^\top\nabla_{K_t^i}C^i(\pmb{K}^i,\pmb{K}^{-i*})\big)}{\|R_t^{i}+(B_t^i)^{\top}P_{t+1,i}^{\pmb{K}^i,\pmb{K}^{-i*}}B_t^i\|}.
\end{eqnarray*}
\end{Lemma}
\begin{proof}
By Lemma 
\ref{lemma:nonzero_cost_diff}, we have
\begin{eqnarray*}
     &&A_{\pmb{K}}^i(x_t^{\pmb{K}^{i\prime},\pmb{K}^{-i}},u_t^{i\prime,-i},t)\\
    &=& (x_t^{\pmb{K}^{i\prime},\pmb{K}^{-i}})^{\top}(K_t^{i\prime}-K_t^i)^{\top}(R_t^{i}+(B_t^i)^{\top}P_{t +1,i}^{\pmb{K}}B_t^i)(K_t^{i\prime}-K_t^i)\,x_t^{\pmb{K}^{i\prime},\pmb{K}^{-i}}\\ &&+\,2(x_t^{\pmb{K}^{i\prime},\pmb{K}^{-i}})^{\top}(K_t^{i\prime}-K_t^i)^{\top}E_{t,i}^{\pmb{K}}\,x_t^{\pmb{K}^{i\prime},\pmb{K}^{-i}}\\
    &=& \Tr\big(x_t^{\pmb{K}^{i\prime},\pmb{K}^{-i}}(x_t^{\pmb{K}^{i\prime},\pmb{K}^{-i}})^{\top}(K_t^{i\prime}-K_t^i)^{\top}(R_t^{i}+(B_t^i)^{\top}P_{t +1,i}^{\pmb{K}}B_t^i)(K_t^{i\prime}-K_t^i)\big)\\ &&+\,2\Tr\big(x_t^{\pmb{K}^{i\prime},\pmb{K}^{-i}}(x_t^{\pmb{K}^{i\prime},\pmb{K}^{-i}})^{\top}(K_t^{i\prime}-K_t^i)^{\top}E_{t,i}^{\pmb{K}}\big)\\
    &=& \Tr\big(x_t^{\pmb{K}^{i\prime},\pmb{K}^{-i}}(x_t^{\pmb{K}^{i\prime},\pmb{K}^{-i}})^{\top}(K_t^{i\prime}-K_t^i+(R_t^{i}+(B_t^i)^{\top}P_{t +1,i}^{\pmb{K}}B_t^i)^{-1}E_{t,i}^{\pmb{K}})^{\top}(R_t^{i}+(B_t^i)^{\top}P_{t +1,i}^{\pmb{K}}B_t^i)\\
    &&(K_t^{i\prime}-K_t^i+(R_t^{i}+(B_t^i)^{\top}P_{t +1,i}^{\pmb{K}}B_t^i)^{-1}E_{t,i}^{\pmb{K}})\big)\\
    &&-\Tr\big(x_t^{\pmb{K}^{i\prime},\pmb{K}^{-i}}(x_t^{\pmb{K}^{i\prime},\pmb{K}^{-i}})^{\top}(E_{t,i}^{\pmb{K}})^\top(R_t^{i}+(B_t^i)^{\top}P_{t +1,i}^{\pmb{K}}B_t^i)^{-1}E_{t,i}^{\pmb{K}}\big)\\
    &\geq&-\Tr\big(x_t^{\pmb{K}^{i\prime},\pmb{K}^{-i}}(x_t^{\pmb{K}^{i\prime},\pmb{K}^{-i}})^{\top}(E_{t,i}^{\pmb{K}})^\top(R_t^{i}+(B_t^i)^{\top}P_{t +1,i}^{\pmb{K}}B_t^i)^{-1}E_{t,i}^{\pmb{K}}\big),
\end{eqnarray*}
with equality in the last line when $K_t^{i\prime}=K_t^i-(R_t^{i}+(B_t^i)^{\top}P_{t +1,i}^{\pmb{K}}B_t^i)^{-1}E_{t,i}^{\pmb{K}}$. Then, by Lemma \ref{lemma:nonzero_cost_diff}, letting $\{x_t^*\}_{t=0}^T$ and $u_t^{*}=(u_t^{1*},\cdots,u_t^{N*})$ with $u_t^{i*}=-K_t^{i*}x_t^*$ denote the state and control sequences induced by the set of optimal policies $\pmb{K}^{*}$, we have
\begin{eqnarray}
     && C^i(\pmb{K}^i,\pmb{K}^{-i*})-C^i(\pmb{K}^{*})= - \mathbb{E}\left[\sum_{t=0}^{T-1}A_{\pmb{K}^i,\pmb{K}^{-i*}}^i(x_t^{*},u_t^{*},t)\right]\nonumber\\
     &\leq& \mathbb{E}\left[\sum_{t=0}^{T-1}\Tr\Big(x_t^*(x_t^*)^\top (E_{t,i}^{\pmb{K}^i,\pmb{K}^{-i*}})^\top(R_t^{i}+(B_t^i)^{\top}P_{t+1,i}^{\pmb{K}^i,\pmb{K}^{-i*}}B_t^i)^{-1}E_{t,i}^{\pmb{K}^i,\pmb{K}^{-i*}}\Big)\right]\nonumber\\
     &\leq& \|\Sigma_{\pmb{K}^{*}}\|\sum_{t=0}^{T-1}\Tr\Big((E_{t,i}^{\pmb{K}^i,\pmb{K}^{-i*}})^\top(R_t^{i}+(B_t^i)^{\top}P_{t+1,i}^{\pmb{K}^i,\pmb{K}^{-i*}}B_t^i)^{-1}E_{t,i}^{\pmb{K}^i,\pmb{K}^{-i*}}\Big)\label{eqn:inte_trace_ineq}\\
     &\leq& \frac{\|\Sigma_{\pmb{K}^{*}}\|}{\srmin}\sum_{t=0}^{T-1}\Tr\Big((E_{t,i}^{\pmb{K}^i,\pmb{K}^{-i*}})^\top E_{t,i}^{\pmb{K}^i,\pmb{K}^{-i*}}\Big)\nonumber\\
     &=& \frac{\|\Sigma_{\pmb{K}^{*}}\|}{4\srmin}\sum_{t=0}^{T-1}\Tr\Big((\Sigma_{t}^{\pmb{K}^i,\pmb{K}^{-i*}})^{-1}\nabla_{K_t^i}C^i(\pmb{K}^i,\pmb{K}^{-i*})^\top\nabla_{K_t^i}C^i(\pmb{K}^i,\pmb{K}^{-i*})(\Sigma_{t}^{\pmb{K}^i,\pmb{K}^{-i*}})^{-1}\Big)\nonumber\\
     &\leq& \frac{\|\Sigma_{\pmb{K}^{*}}\|}{4\srmin\,(\sx)^2}\sum_{t=0}^{T-1}\Tr\Big(\nabla_{K_t^i}C^i(\pmb{K}^i,\pmb{K}^{-i*})^\top\nabla_{K_t^i}C^i(\pmb{K}^i,\pmb{K}^{-i*})\Big)\nonumber.
\end{eqnarray}
Note that \eqref{eqn:inte_trace_ineq} holds as $\Tr(AB) \leq \sigma_{\max}(A) \Tr(B)$ for any matrix $A$ and real symmetric positive semi-definite matrices $B$ of the same size, \cite{Saniuk1987}. From \cite[Lemma 1]{wang1986trace}, for any symmetric matrix $A$ and any symmetric positive semi-definite matrix $B$, it holds that
\begin{equation}\label{eqn:trace_inequality}
    \sigma_{\min}(A) \Tr(B) \leq \Tr(AB) \leq \sigma_{\max}(A)\Tr(B).
\end{equation}
These bounds will be used in several places.
For the lower bound, consider $K_t^{i\prime}=K_t^i-(R_t^{i}+(B_t^i)^{\top}P_{t +1,i}^{\pmb{K}}B_t^i)^{-1}E_{t,i}^{\pmb{K}}\,$. Using $C^i(\pmb{K}^{i\prime},\pmb{K}^{-i*})\geq C^i(\pmb{K}^{*})$ and letting $\{u_t^{\pmb{K}^{i\prime},\pmb{K}^{-i*}}\}_{t=0}^{T-1}$ denote the control 
sequence induced by $(\pmb{K}^{i\prime},\pmb{K}^{-i*})$, by Lemma \ref{lemma:nonzero_cost_diff} we have
\begin{eqnarray*}
     &&C^i(\pmb{K}^i,\pmb{K}^{-i*})-C^i(\pmb{K}^{*})\\
     &\geq& C^i(\pmb{K}^i,\pmb{K}^{-i*})-C^i(\pmb{K}^{i\prime},\pmb{K}^{-i*})\\
     &=& -\mathbb{E}\left[\sum_{t=0}^{T-1}A_{\pmb{K}^i,\pmb{K}^{-i*}}^i(x_t^{\pmb{K}^{i\prime},\pmb{K}^{-i*}},u_t^{\pmb{K}^{i\prime},\pmb{K}^{-i*}}, t)\right]\\
     &=& \mathbb{E}\left[\sum_{t=0}^{T-1}\Tr\Big(x_t^{\pmb{K}^{i\prime},\pmb{K}^{-i*}}(x_t^{\pmb{K}^{i\prime},\pmb{K}^{-i*}})^\top (E_{t,i}^{\pmb{K}^i,\pmb{K}^{-i*}})^\top(R_t^{i}+(B_t^i)^{\top}P_{t+1,i}^{\pmb{K}^i,\pmb{K}^{-i*}}B_t^i)^{-1}E_{t,i}^{\pmb{K}^i,\pmb{K}^{-i*}}\Big)\right]\\
     &\geq& \sx\sum_{t=0}^{T-1}\frac{1}{\|R_t^{i}+(B_t^i)^{\top}P_{t+1,i}^{\pmb{K}^i,\pmb{K}^{-i*}}B_t^i\|}\Tr\big((E_{t,i}^{\pmb{K}^i,\pmb{K}^{-i*}})^\top E_{t,i}^{\pmb{K}^i,\pmb{K}^{-i*}}\big)\\
     &\geq& \frac{\sx}{4\|\Sigma_{\pmb{K}^i,\pmb{K}^{-i*}}\|^2}\sum_{t=0}^{T-1}\frac{\Tr\big(\nabla_{K_t^i}C^i(\pmb{K}^i,\pmb{K}^{-i*})^\top\nabla_{K_t^i}C^i(\pmb{K}^i,\pmb{K}^{-i*})\big)}{\|R_t^{i}+(B_t^i)^{\top}P_{t+1,i}^{\pmb{K}^i,\pmb{K}^{-i*}}B_t^i\|}.
\end{eqnarray*}
\end{proof}

As in the single-agent case \cite{FGKM2018,hambly2020policy}, we now provide an expression for $ C^i(\pmb{K}^{i\prime},\pmb{K}^{-i}) - C^i(\pmb{K})$, which is easier to analyze.
\begin{Lemma}[Almost Smoothness]\label{lemma:nonzero_almost_smoothness}
For two sets of policies $\pmb{K}$ and $(\pmb{K}^{i\prime},\pmb{K}^{-i})$, we have that
\begin{eqnarray*}
 C^i(\pmb{K}^{i\prime},\pmb{K}^{-i}) - C^i(\pmb{K}) 
 &=&\sum_{t=0}^{T-1}\Big[\Tr\big(\Sigma_t^{\pmb{K}^{i\prime},\pmb{K}^{-i}} (K_t^{i\prime}-K_t^i)^{\top}(R_{t}^{i}+(B_t^i)^{\top}P_{t+1,i} ^{\pmb{K}}B_t^i)(K_t^{i\prime}-K_t^i)\big)\\
 &&\quad+2\Tr\big(\Sigma_t^{\pmb{K}^{i\prime},\pmb{K}^{-i}} (K_t^{i\prime}-K_t^i)^{\top}E_{t,i}^{\pmb{K}}\big)\Big].
\end{eqnarray*}
\end{Lemma}
\begin{proof}
By Lemma \ref{lemma:nonzero_cost_diff},
\begin{eqnarray*}
 C^i(\pmb{K}^{i\prime},\pmb{K}^{-i}) - C^i(\pmb{K}) &=& \mathbb{E}\left[\sum_{t=0}^{T-1}A_{\pmb{K}}^i(x_t^{\pmb{K}^{i\prime},\pmb{K}^{-i}},u_t^{i\prime,-i},t)\right]\\
 &=&\sum_{t=0}^{T-1}\left[\Tr\big(\Sigma_t^{\pmb{K}^{i\prime},\pmb{K}^{-i}} (K_t^{i\prime}-K_t^i)^{\top}(R_{t}^{i}+(B_t^i)^{\top}P_{t+1,i} ^{\pmb{K}}B_t^i)(K_t^{i\prime}-K_t^i)\big)\right.\\
 && \left.+2\Tr\big(\Sigma_t^{\pmb{K}^{i\prime},\pmb{K}^{-i}} (K_t^{i\prime}-K_t^i)^{\top}E_{t,i}^{\pmb{K}}\big)\right].
\end{eqnarray*}
\end{proof}
% \begin{Lemma}[Almost Smoothness]\label{lemma:nonzero_almost_smoothness}
% It holds that
% \begin{eqnarray*}
%  C^i(\pmb{K}^{i\prime},\pmb{K}^{-i*}) - C^i(\pmb{K}^{i},\pmb{K}^{-i*}) 
%  &=&\sum_{t=0}^{T-1}\Big[\Tr\big(\Sigma_t^{\pmb{K}^{i\prime},\pmb{K}^{-i*}} (K_t^{i\prime}-K_t^i)^{\top}(R_{t}^{i}+(B_t^i)^{\top}P_{t+1,i} ^{\pmb{K}^i,\pmb{K}^{-i*}}B_t^i)(K_t^{i\prime}-K_t^i)\big)\\
%  && +2\Tr\big(\Sigma_t^{\pmb{K}^{i\prime},\pmb{K}^{-i*}} (K_t^{i\prime}-K_t^i)^{\top}E_{t,i}^{\pmb{K}^i,\pmb{K}^{-i*}}\big)\Big].
% \end{eqnarray*}
% \end{Lemma}
% \begin{proof}
% By Lemma \ref{lemma:nonzero_cost_diff},
% \begin{eqnarray*}
%  C^i(\pmb{K}^{i\prime},\pmb{K}^{-i*}) - C^i(\pmb{K}^{i},\pmb{K}^{-i*}) &=& \mathbb{E}\big[\sum_{t=0}^{T-1}A_{\pmb{K}^i,\pmb{K}^{-i*}}^i(x_t^{\pmb{K}^{i\prime},\pmb{K}^{-i*}},u_t^{i\prime,-i*},t)\big]\\
%  &=&\sum_{t=0}^{T-1}\Big[\Tr\big(\Sigma_t^{\pmb{K}^{i\prime},\pmb{K}^{-i*}} (K_t^{i\prime}-K_t^i)^{\top}(R_{t}^{i}+(B_t^i)^{\top}P_{t+1,i} ^{\pmb{K}^i,\pmb{K}^{-i*}}B_t^i)(K_t^{i\prime}-K_t^i)\big)\\
%  && +2\Tr\big(\Sigma_t^{\pmb{K}^{i\prime},\pmb{K}^{-i*}} (K_t^{i\prime}-K_t^i)^{\top}E_{t,i}^{\pmb{K}^i,\pmb{K}^{-i*}}\big)\Big].
% \end{eqnarray*}
% \end{proof}
\begin{Lemma}\label{lemma:nonzero_bds_P_Sigma}
Assume Assumptions \ref{ass:nonzero_cost} and \ref{ass:nonzero_initial_noise} hold. Then for $t=0,1,\dots,T$ and $i=1,\dots,N$
\[\left\|P_{t,i}^{\pmb{K}}\right\|\leq \frac{C^i(\pmb{K})}{\sx} ,\,\,\|\Sigma_{\pmb{K}}\|\leq \frac{C^i(\pmb{K})}{\sqmin}.
\]
%where $\sx$ and $\sqmin$ are defined in \eqref{defn_sx_nonzero} and \eqref{eqn:defn_srmin_sqmin}.
\end{Lemma}

\begin{proof}
By the trace inequality \eqref{eqn:trace_inequality}, it is straightforward to see that
\begin{equation*}
    C^i(\pmb{K})\geq\mathbb{E}[x_t^{\top}P_{t,i}^{\pmb{K}}x_t] \geq \left\|P_{t,i}^{\pmb{K}}\right\|\sigma_{\min}(\mathbb{E}[x_tx_t^{\top}])\geq \sx\left\|P_{t,i}^{\pmb{K}}\right\|
\end{equation*}
and
\begin{eqnarray*}
    C^i(\pmb{K})&=&\sum_{t=0}^{T-1}\Tr\big(\mathbb{E}[x_tx_t^{\top}](Q_t^i+(K_t^i)^{\top}R_t^{i}K_t^i)\big)+\Tr\big(\mathbb{E}[x_Tx_T^{\top}]Q_T^i\big)\\
    &\geq& \min_{t\in[0,T]} \sigma_{\min}(Q_t^i)\cdot\Tr(\Sigma_{\pmb{K}})\\
    &\geq& \sqmin\cdot\|\Sigma_{\pmb{K}}\|.
\end{eqnarray*}
Then the statements in Lemma \ref{lemma:nonzero_bds_P_Sigma} follow since under Assumptions \ref{ass:nonzero_cost} and \ref{ass:nonzero_initial_noise}, we have $\sx>0$ and $\sqmin>0$.
\end{proof}

\subsection{Perturbation Analysis of the State Covariance Matrix}

Our aim in this section is to provide an explicit control of the change in the state covariance matrix after a change in policy $\pmb{K}^i$. We begin by defining two linear operators on symmetric matrices. For {$X\in\mathbb{R}^{d\times d}$} we set
\begin{equation}\label{eqn:defn_F_t}
    \mathcal{F}_{K_t}(X):=\left(A_t-\sum_{j=1}^{N}B_t^jK_t^j      \right)X\left(A_{t-i}-\sum_{j=1}^{N}B_{t-i}^jK_{t-i}^j\right)^{\top}, 
\end{equation}
and
\begin{equation*}
    \mathcal{T}_{\pmb{K}}(X):=X+\sum_{t=0}^{T-1} \Pi_{i=0}^t \left(A_t-\sum_{j=1}^{N}B_t^jK_t^j\right)\,X\,\Pi_{i=0}^t \left(A_{t-i}-\sum_{j=1}^{N}B_{t-i}^jK_{t-i}^j\right)^{\top},
\end{equation*}
where $\Pi$ denotes the multiplication from the left.
If we write $\mathcal{G}_{t}^{\pmb{K}} =\mathcal{F}_{K_t}\circ\mathcal{F}_{K_{t-1}}\circ\cdots\circ \mathcal{F}_{K_0}$, then the following relationships hold: %then {\color{red}($\mathcal{G}_{K_t}$ and $\mathcal{G}_{t}^{i'}$ are a bit confusing. better notation?)}
\begin{equation}\label{eqn:defn_Gt}
\mathcal{G}_{t}^{\pmb{K}}(X) = \mathcal{F}_{K_t}\circ \mathcal{G}_{t-1}^{\pmb{K}}(X) = \Pi_{i=0}^t \left(A_t-\sum_{j=1}^{N}B_t^jK_t^j\right)\,X\,\Pi_{i=0}^t \left(A_{t-i}-\sum_{j=1}^{N}B_{t-i}^jK_{t-i}^j\right)^{\top},
\end{equation}
and 
\begin{equation}\label{eqn:defn_T_K1K1}
    \mathcal{T}_{\pmb{K}}(X) = X\,+\,\sum_{t=0}^{T-1} \mathcal{G}_{t}^{\pmb{K}}(X).
\end{equation}
When the policy $\pmb{K}$ is clear we will write $\mathcal{G}_{t}=\mathcal{G}_{t}^{\pmb{K}}$ and $\mathcal{F}_{t}=\mathcal{F}_{K_t}$.

We also define the induced norm for these operators as
\begin{equation}\label{defn_operator_norm}
    \|T\|:=\sup_X \frac{\|T(X)\|}{\|X\|},
\end{equation}
where $T=\mathcal{F}_{K_t},\mathcal{G}_t^{\pmb{K}}, \mathcal{T}_{\pmb{K}}$ and the supremum is over all symmetric matrix $X$ with non-zero spectral norm.

We first show the relationship between the operator  $\mathcal{T}_{\pmb{K}}$ and the quantity $\Sigma_{\pmb{K}}$.

\begin{Proposition}\label{prop:Sigma_Gamma_relation}
For $T \ge 2$, we have that
\begin{equation}
  \Sigma_{\pmb{K}}  =\mathcal{T}_{\pmb{K}}(\Sigma_0)+\Delta_{\pmb{K}}(W),
\end{equation}
where
\begin{equation*}
    \Delta_{\pmb{K}}(W) = \sum_{t=1}^{T-1} \sum_{s=1}^{t}\,D_{t,s} W D_{t,s}^{\top}+T\,W,
\end{equation*}
with $D_{t,s} = \Pi_{u=s}^{t} (A_u-\sum_{j=1}^NB_u^jK_{u}^j)$ (for $s=1,2,\cdots,t$),  and $ \Sigma_0=\mathbb{E}\left[x_0x_0^{\top}\right]$.
\end{Proposition}

\begin{proof}
Recall that $\Sigma_{t}^{\pmb{K}} =  \mathbb{E}\left[x_t x_t^{\top}\right]$ and note that 
\begin{eqnarray*}
\Sigma_1^{\pmb{K}} &=& \mathbb{E}[x_{1}x_{1}^\top] = \mathbb{E}\left[\left(\left(A_0-\sum_{j=1}^NB_0^jK_0^j\right)x_0+w_0\right)\left(\left(A_0-\sum_{j=1}^NB_0^jK_0^j\right)x_0+w_0\right)^\top\right]\\ \nonumber
&=& \left(A_0-\sum_{j=1}^NB_0^jK_0^j\right)\,\Sigma_0 \left(A_0-\sum_{j=1}^NB_0^jK_0^j\right)^\top +W = \mathcal{G}_0(\Sigma_0) +W.
\end{eqnarray*}
We first prove that for $t=2,3,\cdots,T$
\begin{eqnarray}\label{eqn:intermedidate1}
\Sigma_t^{\pmb{K}} = \mathcal{G}_{t-1}(\Sigma_0) + \sum_{s=1}^{t-1} D_{t-1,s} W D_{t-1,s}^{\top} +W.
\end{eqnarray}
%\noindent To see this, we have when $t=2$,
%\begin{eqnarray*}
%\Sigma_2^{\pmb{K}} &=& \mathbb{E}[x_{2}x_{2}^\top] = \mathbb{E}\left[\left(\left(A_1-\sum_{j=1}^NB_1^jK_1^j\right) x_1+w_1\right)\left(\left(A_1-\sum_{j=1}^NB_1^jK_1^j\right)x_1+w_1\right)^\top\right]\\
%&=& \left(A_1-\sum_{j=1}^NB_1^jK_1^j\right)\Sigma_1^{\pmb{K}} \left(A_1-\sum_{j=1}^NB_1^jK_1^j\right)^\top +W \\
%&=& \mathcal{G}_1(\Sigma_0) + \left(A_1-\sum_{j=1}^NB_1^jK_1^j\right)W\left(A_1-\sum_{j=1}^NB_1^jK_1^j\right)^\top + W,
%\end{eqnarray*}
%which satisfies \eqref{eqn:intermedidate1}. 
We have the result for $t=1$, so assume \eqref{eqn:intermedidate1} holds for $t\leq k$. Then for $t =  k+1$,
\begin{eqnarray*}
\mathbb{E}[x_{t+1}x_{t+1}^\top] &= &\mathbb{E}\left[\left(\left(A_t-\sum_{j=1}^NB_t^jK_t^j\right)x_t+w_t\right)\left(\left(A_t-\sum_{j=1}^NB_t^jK_t^j\right)x_t+w_t\right)^\top\right]\\
&=& \left(A_t-\sum_{j=1}^NB_t^jK_t^j\right)\Sigma_t^{\pmb{K}} \left(A_t-\sum_{j=1}^NB_t^jK_t^j\right)^\top +W
%&= & (A-B\,K_t)\left(  \mathcal{G}_{t-1}(\Sigma_0) + \sum_{s=1}^{t-1} D_{t-1,s} W D_{t-1,s}^{\top} +W\right) \left(A-B\,K_t\right)^\top +W.\\
%&= & \mathcal{G}_{t}(\Sigma_0) + \sum_{s=1}^{t-1} D_{t,s} W D_{t,s}^{\top} +D_{t,t} W D_{t,t}^{\top}+W\\
=  \mathcal{G}_{t}(\Sigma_0) + \sum_{s=1}^{t} D_{t,s} W D_{t,s}^{\top} +W.
\end{eqnarray*}
Therefore \eqref{eqn:intermedidate1} holds,  $\forall\,t=1,2,\cdots,T$. Finally,
\[
\Sigma_{\pmb{K}} = \sum_{t=0}^T \Sigma_t^{\pmb{K}} = \Sigma_0
+\sum_{t=0}^{T-1}\mathcal{G}_t(\Sigma_0) +\sum_{t=1}^{T-1}\sum_{s=1}^t D_{t,s}WD_{t,s}^{\top} + TW = \mathcal{T}_{\pmb{K}}(\Sigma_0)+\Delta_{\pmb{K}}(W).
\]
\end{proof}

Given two policies $\pmb{K}$ and $\pmb{K}^{\prime}=(\pmb{K}^{1\prime},\cdots,\pmb{K}^{N\prime})$, let us define
\begin{equation}\label{eqn:defn_rho}
\begin{split}
    \rho_{\pmb{K},\pmb{K}^\prime} 
    &:= \max\Bigg\{\max_i\bigg\{\max_{0\leq t \leq T-1}\bigg\|A_t-B_t^iK_t^{i}-\sum_{j=1,j\neq i}^N B_t^jK_t^{j*}\bigg\|\bigg\},\max_{0\leq t \leq T-1}\bigg\|A_t-\sum_{i=1}^N B_t^iK_t^i\bigg\|,\\
    &\qquad\qquad\max_i\bigg\{\max_{0\leq t \leq T-1}\bigg\|A_t-B_t^iK_t^{i\prime}-\sum_{j=1,j\neq i}^N B_t^jK_t^{j*}\bigg\|\bigg\},\\
&\qquad\qquad \max_i\bigg\{\max_{0\leq t \leq T-1}\bigg\|A_t-B_t^iK_t^{i\prime}-\sum_{j=1,j\neq i}^N B_t^jK_t^{j}\bigg\|\bigg\},\,1+\xi\Bigg\},
\end{split}
\end{equation}
for some small constant $\xi>0$.

For any given policy $\pmb{K}$,  $\max_{0\leq t \leq T-1}\left\|A_t-\sum_{j=1}^N B_t^jK_t^{j}\right\|$ measures the radius of the state dynamics under policy $\pmb{K}$. Thus $ \rho_{\pmb{K},\pmb{K}^\prime}$ defined in  \eqref{eqn:defn_rho} is the maximum radius of the policies $\pmb{K}$, $(\pmb{K}^i,\pmb{K}^{-i*})$,  $(\pmb{K}^{i\prime},\pmb{K}^{-i*})$, and $(\pmb{K}^{i\prime},\pmb{K}^{-i})$ $(i=1,2,\cdots,N)$. We will show later that the value (or the upper bound) of $ \rho_{\pmb{K},\pmb{K}^\prime}$ plays an essential role in the convergence analysis.

\begin{Remark}{\rm
By the definition of $\rho_{\pmb{K},\pmb{K}^\prime}$ in \eqref{eqn:defn_rho}, we have $\rho_{\pmb{K},\pmb{K}^\prime} \ge 1+\xi>1$. This regularization term $1+\xi$ is introduced to simplify the presentation. Alternatively, if we remove this term from the definition of $\rho_{\pmb{K},\pmb{K}^\prime}$, a similar analysis can still be carried out by considering the different cases: $\rho_{\pmb{K},\pmb{K}^\prime}<1$, $\rho_{\pmb{K},\pmb{K}^\prime}=1$ and $\rho_{\pmb{K},\pmb{K}^\prime}>1$.}
%$\rho_{\pmb{K},\pmb{K}^\prime} := \max\Big\{\max_{0\leq t \leq T-1}\|A_t-\sum_{i=1}^N B_t^iK_t^i\|,\max_i\big\{\max_{0\leq t \leq T-1}\|A_t-\sum_{j=1,j\neq i}^N B_t^jK_t^{j}-B_t^iK_t^{i\prime}\|\big\},\max_i\big\{\max_{0\leq t \leq T-1}\|A_t-\sum_{j=1,j\neq i}^N B_t^jK_t^{j*}-B_t^iK_t^{i}\|\big\},\max_i\big\{\max_{0\leq t \leq T-1}\|A_t-\sum_{j=1,j\neq i}^N B_t^jK_t^{j*}-B_t^iK_t^{i\prime}\|\big\}\Big\}$, a similar analysis can still be carried out by considering the different cases: $\rho_{\pmb{K},\pmb{K}^\prime}<1$, $\rho_{\pmb{K},\pmb{K}^\prime}=1$ and $\rho_{\pmb{K},\pmb{K}^\prime}>1$.}
\end{Remark}
We now provide an upper bound for $\rho_{\pmb{K},\pmb{K}^{\prime}}$.
\begin{Lemma}\label{lemma:rho_upper_bd}
Assume Assumption \ref{ass:nonzero_exist_sol} holds. Then, 
\begin{eqnarray}
\rho_{\pmb{K},\pmb{K}^\prime}
% &\leq& \rho^* +N\gamma_B\sqrt{\frac{T}{\sx\srmin} \max_i\{C^i(\pmb{K}^{i},\pmb{K}^{-i*})- C^i(\pmb{K}^{*}),C^i(\pmb{K}^{i\prime},\pmb{K}^{-i*})- C^i(\pmb{K}^{*})\}}\\
&\leq& \rho^* +N\gamma_B\sqrt{\frac{T}{\sx\srmin} \max_i\Big\{C^i(\pmb{K}^{i},\pmb{K}^{-i*})- C^i(\pmb{K}^{*})\Big\}}+\gamma_B\max_i\max_t\Big\{\|K_t^{i\prime}-K_t^i\|\Big\}\label{eqn:rho_upper_bd}
\end{eqnarray}
where %$\rho_{\pmb{K},\pmb{K}^\prime}$ and 
$\rho^*$ was defined in %\eqref{eqn:defn_rho} and 
\eqref{eqn:defn_rho_star}.
\end{Lemma}
\begin{proof}
By Lemma \ref{lemma:nonzero_almost_smoothness}, we have
\begin{eqnarray}
C^i(\pmb{K}^{i},\pmb{K}^{-i*}) - C^i(\pmb{K}^{*}) 
 &=&\sum_{t=0}^{T-1}\Big[\Tr\big(\Sigma_t^{\pmb{K}^{i},\pmb{K}^{-i*}} (K_t^{i}-K_t^{i*})^{\top}(R_{t}^{i}+(B_t^i)^{\top}P_{t+1,i} ^{\pmb{K}^{*}}B_t^i)(K_t^{i}-K_t^{i*})\big)\Big]\nonumber\\
 &\geq& \sx\srmin\sum_{t=0}^{T-1}\|K_t^i-K_t^{i*}\|^2\geq \frac{\sx\srmin}{T} \vertiii{\pmb{K}^i-\pmb{K}^{i*}}^2,\label{eqn:Ci_diff_Ki}
\end{eqnarray}
where \eqref{eqn:Ci_diff_Ki} holds by the Cauchy-Schwarz inequality.
Then we have
\begin{eqnarray}
\left\|A_t-\sum_{j=1}^NB_t^jK_t^j\right\| &\leq& \left\|A_t-\sum_{j=1}^NB_t^jK_t^{j*}\right\|+\sum_{j=1}^N\|B_t^j\|\,\|K_t^j-K_t^{j*}\|\nonumber\\
&\leq& \left\|A_t-\sum_{j=1}^NB_t^jK_t^{j*}\right\|+\gamma_B\sum_{j=1}^N \vertiii{\pmb{K}^j-\pmb{K}^{j*}}\nonumber\\
&\leq& \left\|A_t-\sum_{j=1}^NB_t^jK_t^{j*}\right\|+N\gamma_B\sqrt{\frac{T}{\sx\srmin}\max_j\{ C^j(\pmb{K}^{j},\pmb{K}^{-j*}) - C^j(\pmb{K}^{*})\}},\label{eqn:rho_bd_inte1}
\end{eqnarray}
where \eqref{eqn:rho_bd_inte1} holds by \eqref{eqn:Ci_diff_Ki}.
Also, by the triangle inequality we have
\begin{eqnarray}\label{eqn:rho_bd_inte2}
\left\|A_t-\sum_{j=1,j\neq i}^N B_t^jK_t^{j}-B_t^iK_t^{i\prime}\right\|\leq \left\|A_t-\sum_{j=1}^NB_t^jK_t^j\right\| + \gamma_B\|K_t^{i\prime}-K_t^i\|.
\end{eqnarray}
and
\begin{eqnarray}\label{eqn:rho_bd_inte4}
\left\|A_t-\sum_{j=1,j\neq i}^N B_t^jK_t^{j*}-B_t^iK_t^{i\prime}\right\|\leq \left\|A_t-\sum_{j=1,j\neq i}^N B_t^jK_t^{j*}-B_t^iK_t^{i}\right\| + \gamma_B\|K_t^{i\prime}-K_t^i\|.
\end{eqnarray}
Finally, applying \eqref{eqn:Ci_diff_Ki}, 
\begin{eqnarray}
\left\|A_t-\sum_{j=1,j\neq i}^N B_t^jK_t^{j*}-B_t^iK_t^{i}\right\|&\leq& \left\|A_t-\sum_{j=1}^NB_t^jK_t^{j*}\right\|+\gamma_B\|K_t^{i*}-K_t^i\|\nonumber\\
&\leq& \left\|A_t-\sum_{j=1}^NB_t^jK_t^{j*}\right\|\nonumber\\
&&+\gamma_B\sqrt{\frac{T}{\sx\srmin} \big(C^i(\pmb{K}^{i},\pmb{K}^{-i*}) - C^i(\pmb{K}^{*})\big)}.\label{eqn:rho_bd_inte3}
\end{eqnarray}

Therefore combining \eqref{eqn:rho_bd_inte1}-\eqref{eqn:rho_bd_inte3}, we obtain the statement \eqref{eqn:rho_upper_bd}.
\end{proof}

Recall the definition of $\mathcal{F}_{t}=\mathcal{F}_{K_t}$ and $\mathcal{G}_t=\mathcal{G}_{t}^{\pmb{K}}$ in \eqref{eqn:defn_F_t} and \eqref{eqn:defn_Gt} associated with $\pmb{K}$, similarly let us define $\mathcal{G}^{i\prime}_t =\mathcal{F}_{K^{i\prime}_t,K_t^{-i}}\circ\mathcal{F}_{K^{i\prime}_{t-1},K_{t-1}^{-i}}\circ\cdots\circ \mathcal{F}_{K^{i\prime}_0,K_0^{-i}}$ for the set of policies $(\pmb{K}^{i\prime},\pmb{K}^{-i})$ and write $\mathcal{F}_t^{i\prime}=\mathcal{F}_{K_t^{i\prime},K_t^{-i}}$. We now establish a perturbation analysis for $\mathcal{F}_t$ and $\mathcal{G}_t$.
%in Lemmas \ref{lemma:bd_F_K1K2} and \ref{lemma:G_t_pertub}. 

\begin{Lemma}\label{lemma:bd_F_K1K2}
For all $t=0,1,\cdots,T-1$, we have
\begin{eqnarray}\label{eqn:F_K1pri_K2_bound}
\|\mathcal{F}_{t}-\mathcal{F}_{t}^{i\prime}\| \leq 2\rho_{\pmb{K},\pmb{K}^\prime}\gamma_B\|K_t^i-K_t^{i\prime}\|.
\end{eqnarray}
\end{Lemma}

\begin{proof}
For player $i$,
\begin{eqnarray*}
    &&(\mathcal{F}_{t}-\mathcal{F}_{t}^{i\prime})(X)\\
    &=&\left(A_t-\sum_{j=1}^NB_t^jK_t^j\right)X\left(A_t-\sum_{j=1}^NB_t^jK_t^j\right)^{\top} \\
     & & \qquad -
    \left(A_t-B_t^iK_t^{i\prime}-\sum_{j=1,j\neq i}^NB_t^jK_t^j\right)X\left(A_t-B_t^iK_t^{i\prime}-\sum_{j=1,j\neq i}^NB_t^jK_t^j\right)^{\top},
\end{eqnarray*}
By \eqref{defn_operator_norm} the operator norm of $\mathcal{F}_{t}-\mathcal{F}_{t}^{i\prime}$ is the maximum possible ratio of $\|(\mathcal{F}_{t}-\mathcal{F}_{t}^{i\prime})(X)\|$ and $\|X\|$. Then letting $Y=A_t-\sum_{j=1}^NB_t^jK_t^j$ and $Z=A_t-B_t^iK_t^{i\prime}-\sum_{j=1,j\neq i}^NB_t^jK_t^j$ in 
\begin{equation}\label{eqn:diff_YXY}
    Y XY^\top -Z XZ^\top =\frac{(Y+Z) X(Y-Z)^\top + (Y-Z) X(Y+Z)^\top}{2}
\end{equation}
and using the norm bound $\|AX\|\leq \|A\|\,\|X\|$, we have 
\begin{eqnarray*}
    &&\Big\|\left(A_t-\sum_{j=1}^NB_t^jK_t^j\right)X\left(A_t-\sum_{j=1}^NB_t^jK_t^j\right)^{\top} \\
     && -
    \left(A_t-B_t^iK_t^{i\prime}-\sum_{j=1,j\neq i}^NB_t^jK_t^j\right)X\left(A_t-B_t^iK_t^{i\prime}-\sum_{j=1,j\neq i}^NB_t^jK_t^j\right)^{\top}\Big\|\\
    &\leq& 2\rho_{\pmb{K},\pmb{K}^\prime} \|X\|\,\|B_t^i(K_t^i-K_t^{i'})\| \leq 2\rho_{\pmb{K},\pmb{K}^\prime} \|X\|\gamma_B\|(K_t^i-K_t^{i'})\|.
\end{eqnarray*}
Therefore we obtain the statement $\|\mathcal{F}_{t}-\mathcal{F}_{t}^{i\prime}\|\leq 2\rho_{\pmb{K},\pmb{K}^\prime}\gamma_B\|(K_t^i-K_t^{i'})\|$
% For player $i$,
% \begin{eqnarray*}
%     &&(\mathcal{F}_{t}-\mathcal{F}_{t}^{i\prime})(X)\\
%     &=&\left(A_t-\sum_{j=1}^NB_t^jK_t^j\right)X\left(A_t-\sum_{j=1}^NB_t^jK_t^j\right)^{\top} \\
%      & & \qquad -
%     \left(A_t-B_t^iK_t^{i\prime}-\sum_{j=1,j\neq i}^NB_t^jK_t^j\right)X\left(A_t-B_t^iK_t^{i\prime}-\sum_{j=1,j\neq i}^NB_t^jK_t^j\right)^{\top},
% \end{eqnarray*}
% Then \eqref{eqn:F_K1pri_K2_bound} follows by letting $Y=A_t-\sum_{j=1}^NB_t^jK_t^j$ and $Z=A_t-B_t^iK_t^{i\prime}-\sum_{j=1,j\neq i}^NB_t^jK_t^j$ in 
% \begin{equation}\label{eqn:diff_YXY}
%     Y XY^\top -Z XZ^\top =\frac{(Y+Z) X(Y-Z)^\top + (Y-Z) X(Y+Z)^\top}{2}
% \end{equation}
% and using the norm bound $\|AX\|\leq \|A\|\,\|X\|$.
\end{proof}

Recall that $\mathcal{F}_t$ and $\mathcal{G}_t$ are defined in equations \eqref{eqn:defn_F_t} and \eqref{eqn:defn_Gt}. Then we have the following lemma on perturbation analysis for $\mathcal{G}_t$.

\begin{Lemma}[Perturbation Analysis for  $\mathcal{G}_t$]\label{lemma:G_t_pertub}
For any symmetric matrix $\Sigma\in\mathbb{R}^{d\times d}$ and $i=1,\cdots,N$, we have that 
\[
\sum_{t=0}^{T-1}\Big\|(\mathcal{G}_{t}-\mathcal{G}_{t}^{i\prime})(\Sigma)\Big\| \leq \frac{\rho_{\pmb{K},\pmb{K}^\prime}^{2T}-1}{\rho_{\pmb{K},\pmb{K}^\prime}^2-1} \Big( \sum_{t=0}^{T-1}\|\mathcal{F}_{t}-\mathcal{F}^{i\prime}_{t}\|\Big)\|\Sigma\|.
\]
%where $\rho_{\pmb{K},\pmb{K}^\prime}$ is defined in \eqref{eqn:defn_rho}.
\end{Lemma}

\begin{proof}
By direct calculation,
\begin{eqnarray}\label{eqn:G_t_bound}
\|\mathcal{G}_t^{i\prime}\| \leq \rho_{\pmb{K},\pmb{K}^\prime}^{2(t+1)},\quad i=1,\cdots,N.
\end{eqnarray}
Then for any symmetric matrix $\Sigma\in\mathbb{R}^{d\times d}$ and $t \ge 0$,
\begin{eqnarray*}
\|(\mathcal{G}_{t+1}^{i\prime}-\mathcal{G}_{t+1})(\Sigma)\| &=& \|\mathcal{F}_{t+1}^{i\prime}\circ\mathcal{G}_t^{i\prime}(\Sigma)-\mathcal{F}_{t+1}\circ\mathcal{G}_t(\Sigma)\|\\
&=& \|\mathcal{F}_{t+1}^{i\prime}\circ\mathcal{G}_{t}^{i\prime}(\Sigma)-\mathcal{F}_{t+1}^{i\prime}\circ\mathcal{G}_{t}(\Sigma)+\mathcal{F}_{t+1}^{i\prime}\circ\mathcal{G}_{t}(\Sigma)-\mathcal{F}_{t+1}\circ\mathcal{G}_{t}(\Sigma)\| \\
&\leq &\|\mathcal{F}_{t+1}^{i\prime}\circ\mathcal{G}_{t}^{i\prime}(\Sigma)-\mathcal{F}_{t+1}^{i\prime}\circ\mathcal{G}_{t}(\Sigma)\|+\|\mathcal{F}_{t+1}^{i\prime}\circ\mathcal{G}_{t}(\Sigma)-\mathcal{F}_{t+1}\circ\mathcal{G}_{t}(\Sigma)\| \\
& =  & \|\mathcal{F}_{t+1}^{i\prime}\circ(\mathcal{G}_{t}^{i\prime}-\mathcal{G}_{t})(\Sigma)\|+\|(\mathcal{F}_{t+1}^{i\prime}-\mathcal{F}_{t+1})\circ\mathcal{G}_{t}(\Sigma)\|\\
&\leq & \|\mathcal{F}_{t+1}^{i\prime}\| \,\|(\mathcal{G}_{t}^{i\prime}-\mathcal{G}_{t})(\Sigma)\|+\|\mathcal{G}_{t}\| \,\|\mathcal{F}_{t+1}^{i\prime}-\mathcal{F}_{t+1}\|\,\|\Sigma\|\\
&\leq& \rho_{\pmb{K},\pmb{K}^\prime}^2\|(\mathcal{G}_{t}^{i\prime}-\mathcal{G}_{t})(\Sigma)\| +\rho_{\pmb{K},\pmb{K}^\prime}^{2(t+1)}\|\mathcal{F}_{t+1}^{i\prime}-\mathcal{F}_{t+1}\| \|\Sigma\|.
\end{eqnarray*}
Therefore,
\begin{eqnarray}\label{tn1}
\|(\mathcal{G}_{t+1}^{i\prime}-\mathcal{G}_{t+1})(\Sigma)\|
\leq \rho_{\pmb{K},\pmb{K}^\prime}^2\|(\mathcal{G}_{t}^{i\prime}-\mathcal{G}_{t})(\Sigma)\| +\rho_{\pmb{K},\pmb{K}^\prime}^{2(t+1)}\|\mathcal{F}_{t+1}^{i\prime}-\mathcal{F}_{t+1}\| \|\Sigma\|.
\end{eqnarray}
% {\color{blue}
% Adding \eqref{tn1} when $t=0$ and $\|\mathcal{G}^{i\prime}_0-\mathcal{G}_0\|=\|\mathcal{F}_0^{i\prime}-\mathcal{F}_0\|$ together, we have 
% \begin{eqnarray*}
% \sum_{t=0}^{1}\Big\|(\mathcal{G}_{t}-\mathcal{G}_{t}^{i\prime})(\Sigma)\Big\| &\leq& \left(1+\rho_{\pmb{K},\pmb{K}^\prime}^2\right)\|\mathcal{F}_0-\mathcal{F}_0^{i'}\|\,\|\Sigma\| + \rho_{\pmb{K},\pmb{K}^\prime}^{2}\|\mathcal{F}_{1}^{i\prime}-\mathcal{F}_{1}\| \|\Sigma\|\\
% &\leq& \left(1+\rho_{\pmb{K},\pmb{K}^\prime}^2\right)\Big( \sum_{t=0}^{1}\|\mathcal{F}_{t}-\mathcal{F}^{i\prime}_{t}\|\Big)\|\Sigma\|.
% \end{eqnarray*}
% Similarly,
% \begin{eqnarray*}
% \sum_{t=0}^{2}\Big\|(\mathcal{G}_{t}-\mathcal{G}_{t}^{i\prime})(\Sigma)\Big\| &\leq& \rho_{\pmb{K},\pmb{K}^\prime}^2\left(1+\rho_{\pmb{K},\pmb{K}^\prime}^2\right)\|\mathcal{G}_1-\mathcal{G}_1^{i'}\|\,\|\Sigma\|+\|\mathcal{F}_0-\mathcal{F}_0^{i'}\|\|\Sigma\| + \rho_{\pmb{K},\pmb{K}^\prime}^{4}\|\mathcal{F}_{2}^{i\prime}-\mathcal{F}_{2}\| \|\Sigma\|\\
% &\leq& \left(1+\rho_{\pmb{K},\pmb{K}^\prime}^2+\rho_{\pmb{K},\pmb{K}^\prime}^4\right)\Big( \sum_{t=0}^{2}\|\mathcal{F}_{t}-\mathcal{F}^{i\prime}_{t}\|\Big)\|\Sigma\|.
% \end{eqnarray*}
% }
As it is a geometric series, summing \eqref{tn1} over $t \in \{0,1,2,\cdots,T-2\}$ with $\|\mathcal{G}^{i\prime}_0-\mathcal{G}_0\|=\|\mathcal{F}_0^{i\prime}-\mathcal{F}_0\|$, gives
\[
\sum_{t=0}^{T-1}\Big\|(\mathcal{G}_{t}-\mathcal{G}_{t}^{i\prime})(\Sigma)\Big\| \leq \frac{\rho_{\pmb{K},\pmb{K}^\prime}^{2T}-1}{\rho_{\pmb{K},\pmb{K}^\prime}^2-1} \Big( \sum_{t=0}^{T-1}\|\mathcal{F}_{t}-\mathcal{F}^{i\prime}_{t}\|\Big)\|\Sigma\|.
\]
\end{proof}

Recall that $\gamma_A$, $\gamma_B$, and $\gamma_R$ are defined in \eqref{eqn:defn_gammaAgammaBgammaR}. Then we have the following perturbation analysis of $\Sigma_{\pmb{K}}$.

\begin{Lemma}[Perturbation Analysis of $\Sigma_{\pmb{K}}$]\label{lemma:perturb_Sigma_K1K2}
Assume Assumption \ref{ass:nonzero_cost} holds. Then
\begin{equation*}
\begin{split}
  \Big\|\Sigma_{\pmb{K}}-\Sigma_{\pmb{K}^{i\prime},\pmb{K}^{-i}}\Big\|
   & \leq 2\gamma_B\,\frac{\rho_{\pmb{K},\pmb{K}^\prime}( \rho_{\pmb{K},\pmb{K}^\prime}^{2T}-1)}{\rho_{\pmb{K},\pmb{K}^\prime}^2-1} \left(\frac{C^i(\pmb{K}^i,\pmb{K}^{-i*})}{\sqmin}+T\|W\|\right)\vertiii{\pmb{K}^i-\pmb{K}^{i\prime}}.
\end{split}
\end{equation*}
\end{Lemma}

\begin{proof}
Using Lemma \ref{lemma:bd_F_K1K2},
\begin{eqnarray*}
\sum_{t=0}^{T-1}\|\mathcal{F}_{t}-\mathcal{F}_{t}^{i\prime}\| 
&\leq& 2\rho_{\pmb{K},\pmb{K}^\prime}\gamma_B \sum_{t=0}^{T-1}\|K_t^i-K_t^{i\prime}\|.
\end{eqnarray*}
Define $D_{t,s}^{i\prime} = \Pi_{u=s}^{t} (A_u-B_u^iK_{u}^{i\prime}-\sum_{j=1,j\neq i}^N B_u^jK_u^j)$ (for $s=1,2,\cdots,t$). Then, in a similar way to the proof of Lemma \ref{lemma:G_t_pertub},  we have, $\forall\,t=1,\cdots,T-1$,
\begin{eqnarray}\label{eq:sum_D_bound}
 \sum_{s=1}^{t}\left\|D_{t,s}WD_{t,s}^{\top}-D^{i\prime}_{t,s}W(D^{i\prime}_{t,s})^{\top}\right\| \leq \frac{\rho_{\pmb{K},\pmb{K}^\prime}^{2T}-1}{\rho_{\pmb{K},\pmb{K}^\prime}^2-1} \left(\sum_{s=0}^{t}\|\mathcal{F}_{s}-\mathcal{F}_{s}^{i\prime}\| \right) \|W\|.
\end{eqnarray}
By Proposition \ref{prop:Sigma_Gamma_relation},  \eqref{eqn:defn_T_K1K1} and \eqref{eq:sum_D_bound}, we have
\begin{equation}\label{eq:Colla_intem}
\begin{split}
  \Big\|\Sigma_{\pmb{K}}-\Sigma_{\pmb{K}^{i\prime},\pmb{K}^{-i}}\Big\| & \leq 
  \Big\|(\mathcal{T}_{\pmb{K}}-\mathcal{T}_{\pmb{K}^{i\prime},\pmb{K}^{-i}})(\Sigma_0)\Big\| + \sum_{t=1}^{T-1} \sum_{s=1}^{t}\,\Big\|D_{t,s} W D_{t,s}^{\top}-D_{t,s}^{i\prime} W (D_{t,s}^{i\prime})^{\top}\Big\|\\
     & \leq \frac{\rho_{\pmb{K},\pmb{K}^\prime}^{2T}-1}{\rho_{\pmb{K},\pmb{K}^\prime}^2-1} \Big( \sum_{t=0}^{T-1}\|\mathcal{F}_{t}-\mathcal{F}_{t}^{i\prime}\|\Big)\left(\|\Sigma_0\|+T\|W\|\right) \\
     & \leq {\frac{ \rho_{\pmb{K},\pmb{K}^\prime}^{2T}-1}{\rho_{\pmb{K},\pmb{K}^\prime}^2-1} \left(\frac{C^i(\pmb{K}^i,\pmb{K}^{-i*})}{\sqmin}+T\|W\|\right)\left(2\rho_{\pmb{K},\pmb{K}^\prime}\gamma_B\,\vertiii{\pmb{K}^i-\pmb{K}^{i\prime}}\right)}.
\end{split}
\end{equation}
The last inequality holds since $\|\Sigma_0\| \leq \|\Sigma_{\pmb{K}^i,\pmb{K}^{-i*}}\|\leq \frac{C^i(\pmb{K}^i,\pmb{K}^{-i*})}{\sqmin}$ by Lemma \ref{lemma:nonzero_bds_P_Sigma}.
\end{proof}

\subsection{Convergence and Complexity Analysis}

We are now in a position to provide the proof of our main Theorem \ref{thm:local_conv_NPG}. This will follow from two important Lemmas. 
First, define
\begin{equation}\label{defn_rho_K}
    \rho_{\pmb{K}} := \rho^* +N\gamma_B\sqrt{\frac{T}{\sx\srmin} \max_i\big\{C^i(\pmb{K}^{i},\pmb{K}^{-i*})- C^i(\pmb{K}^{*})\big\}}+\frac{1}{20T^2}.
\end{equation}
Further define $g_1$ and $g_2$ as follows:
\begin{equation}\label{defn_g1}
    g_1 : = \frac{\srmin}{\|\Sigma_{\pmb{K}^{*}}\|},
\end{equation}
and
\begin{equation}\label{defn_g2}
    g_2 := 20(N-1)^2\,T^2\,d\, \frac{(\gamma_B)^4\max_i\{C^i(\pmb{K}^{i},\pmb{K}^{-i*})\}^4}{\sqmin^2\srmin}\left(\frac{\rho_{\pmb{K}}^{2T}-1}{\rho_{\pmb{K}}^2-1}\right)^2.
\end{equation}
We also write 
$C^{i,-i*}=C^i(\pmb{K}^i,\pmb{K}^{-i*})$, $C^{i*} = C^i(\pmb{K}^{*})$ and $C^{i\prime,-i*}=C^i(\pmb{K}^{i\prime},\pmb{K}^{-i*})$ to simplify notation.
%{\color{blue}$\rho_{\pmb{K}}$ is defined in \eqref{defn_rho_K}},  
%$\srmin$ and $\sqmin$ are defined in \eqref{eqn:defn_srmin_sqmin}, and $\gamma_R$ and $\gamma_B$ are defined in \eqref{eqn:defn_gammaAgammaBgammaR}.
\begin{Lemma}[One-step contraction]\label{lemma:one_step_loc_conv}
Assume Assumptions \ref{ass:nonzero_cost}, \ref{ass:nonzero_initial_noise}, and  \ref{ass:nonzero_exist_sol} hold, and that 
% \begin{equation}\label{eqn:one_step_noise_cond}
% \frac{(\sx)^5}{\|\Sigma_{\pmb{K}^{*}}\| } >\frac{5\,d\,T(N-1)^2\rho^2(\rho^{2T}-1)^2}{\sqmin^2 \srmin(1-\rho^2)^2}\Big(1+\max_i\{C^{i}(\pmb{K}^*)\}\Big)^4\max\{4(\gamma_B)^4,(\gamma_B)^6\}\max\left\{{\frac{T}{\srmin}},\frac{T^2(N-1)}{\sx \srmin^2}\right\},
% \end{equation}
\begin{equation}\label{eqn:one_step_noise_cond}
\sx^5 > \frac{g_2}{g_1}.
\end{equation}
%where $g_1$ and $g_2$ are defined in \eqref{defn_g1} and \eqref{defn_g2},
Also assume the policy update step for player $i$ at time $t$ is given by
\begin{equation}\label{eqn:one_step_update}
    K_t^{i\prime} = K_t^i - \eta\nabla_{K_t^i}C^i(\pmb{K})(\Sigma_{t}^{\pmb{K}})^{-1},
\end{equation}
where
\begin{equation}\label{eqn:step_size}
    \eta \leq \min\left\{I_1\,,I_2\,,\frac{1}{\srmin}\right\}
\end{equation}
with
 \begin{equation*}
 \begin{split}
 I_1 & = \left\{20 T\frac{\rho_{\pmb{K}}(\rho_{\pmb{K}}^{2T}-1)}{\rho_{\pmb{K}}^2-1}\left(\sum_{i=1}^NC^{i,-i*} +\sqmin T\|W\|\right)\gamma_B\max_i\{\max_{t}\{\|\nabla_{K_t^i}C^i(\pmb{K})\|\}\}+\sqmin(\sx)^2\right.\\
 &\quad \left.+4\big(\gamma_R\sx+(\gamma_B)^2\sum_{i=1}^NC^{i,-i*}\big)\sum_{i=1}^N\{C^{i,-i*}\}\right\}^{-1}\cdot\sqmin(\sx)^2,\\
 I_2 &= \left\{\big(\max_i\{k_i\}\big)\frac{10 T\sum_{i=1}^NC^{i,-i*}}{(10 T-1)\sqmin}\left(\gamma_R+(\gamma_B)^2\frac{\sum_{i=1}^NC^{i,-i*}}{\sx}\right)+\frac{2d}{\sx}\left(\frac{10 T\sum_{i=1}^NC^{i,-i*}}{(10 T-1)\sqmin}\right)^2\right.\cdot\\
 &\quad\quad\left.\left(\gamma_R+(\gamma_B)^2\frac{\sum_{i=1}^NC^{i,-i*}}{\sx}\right)^2\right\}^{-1}\cdot\frac{d}{80\sx} \left(\frac{10 T\min_i\{C^{i*}\}}{(10 T-1)\sqmin}\right)^2.
 \end{split}
 \end{equation*}
Let $\alpha:=\sx g_1- g_2/\sx^4>0$. Then, we have 
\begin{enumerate}
    \item  $\eta \in (0,\frac{1}{\alpha})$; and
    \item the following inequality holds
\begin{equation}\label{eq:c_sum_onestep_final} \sum_{i=1}^N \left(C^{i}(\pmb{
K}^{i\prime},\pmb{K}^{-i*})-C^{i}(\pmb{K}^*)\right)
\leq (1-\alpha\eta)\left(\sum_{i=1}^N \big(C^{i}(\pmb{K}^i,\pmb{K}^{-i*}) - C^{i}(\pmb{K}^*)\big)\right).
\end{equation}
\end{enumerate}
%\begin{equation}\label{eq:c_sum_onestep_final}
%\sum_{i=1}^N \left(C^{i}(\pmb{
%K}^{i\prime},\pmb{K}^{-i*})-C^{i}(\pmb{K}^*)\right)
%\leq (1-\alpha)\left(\sum_{i=1}^N \big(C^{i}(\pmb{K}^i,\pmb{K}^{-i*}) - C^{i}(\pmb{K}^*)\big)\right),
%\end{equation}
% with {\color{blue}$\alpha=\eta(\sx g_1- g_2/\sx^4)$}
%such that $0<\alpha <1$.
\end{Lemma}

\begin{Remark}
%[Comparison between condition \eqref{eqn:one_step_noise_cond} and condition \eqref{eq:system_noise} in Assumption \ref{ass:noise_cond}]
{\rm
(1) In the one-step contraction analysis (Lemma \ref{lemma:one_step_loc_conv}), \eqref{eqn:one_step_noise_cond} imposes a condition on $C(\pmb{K}^{i},\pmb{K}^{-i*})$ which is associated with the current policy $\pmb{K}$. In the analysis of the global convergence result (Theorem \ref{thm:local_conv_NPG}),   \eqref{eq:system_noise} imposes a similar condition on $C(\pmb{K}^{i,(0)},\pmb{K}^{-i*})$ which is associated with the initial policy $\pmb{K}^{(0)}$. Condition \eqref{eq:system_noise} in Assumption \ref{ass:noise_cond} and the step size condition in Theorem \ref{thm:local_conv_NPG} guarantee that, condition \eqref{eqn:one_step_noise_cond} holds for any $\pmb{K}=\pmb{K}^{(m)}$ throughout the training process $(m=1,2,\cdots,M)$. This further ensures that the one-step contraction analysis in Lemma \ref{lemma:one_step_loc_conv} can be applied iteratively which leads to the global convergence result as stated in Theorem \ref{thm:local_conv_NPG}.\\
%This condition is guaranteed throughout the training procedure $\{\pmb{K}^{(m)}\}_{i=1}^{M}$ when the step size $\eta$ is chosen properly according to the assumption in Theorem \ref{thm:local_conv_NPG}. 
(2) Note that the numbers such as 20 and 80 that appear in $I_1, I_2$ are not arbitrary and, although some minor improvements can be made by optimizing at various stages, they enable us to obtain reasonable bounds.}
%$I_1$ can be sufficiently small such that $\gamma_B\|K_t^{i\prime}-K_t^i\|  \leq \frac{1}{\kappa T^2}$ for some sufficiently large $\kappa\in\mathbb{R}$. Here we choose $I_1$ such that $\kappa = 20$ which leads to reasonable bounds in \eqref{eqn:I1_conseq_3}-\eqref{eqn:I1_conseq_2}. Similarly, $I_2$ can also be sufficiently small such that \eqref{eqn:I2_conseq} $\leq \frac{d\delta}{\sx} \left(\frac{10 T\max_i\{C^{i,-i*}\}}{(10 T-1)\sqmin}\right)^2$ for some sufficiently small $\delta\in\mathbb{R}$. We choose $I_2$ such that $\delta=1/20$ which gives a reasonable bound in \eqref{eqn:I2_conseq}.
\end{Remark}

% \begin{Remark}{\rm
% In Lemma \ref{lemma:one_step_loc_conv}\label{remark:vani} we provide a one-step contraction result for the natural policy gradient method. This result can be developed in the case of the vanilla policy gradient method, where the updating rule is $K_t^{i\prime} = K_t^i - \eta\nabla_{K_t^i}C^i(\pmb{K})$. We provide the lemma for the vanilla method, Lemma \ref{lemma:one_step_cont_vani}, and its proof in Appendix \ref{app:vani}.
% }
% \end{Remark}

\begin{proof}We break this proof up into a series of steps. \\ [.1in]
% Given \eqref{eqn:one_step_update} and $\eta\leq I_1$ from \eqref{eqn:step_size}, we have
% \[
% \|K_t^{i\prime}-K_t^i\| = \eta \|\nabla_{K_t^i} C^i(\pmb{K})\Sigma_t^{-1}\| \leq \frac{\sqmin\sx}{2\gamma_B C^{i,-i*}}.
% \]
% Therefore, 
% \begin{equation}\label{eqn:BdiffK_bd}
%     \|B_t^i\|\|K_t^{i\prime}-K_t^i\|  \leq \frac{\sqmin\sx}{2C^{i,-i*}} \leq \frac{1}{2}.
% \end{equation}
% The last inequality holds since $\sx\leq \frac{C^{i,-i*}}{\sqmin}$ given by Lemma \ref{lemma:nonzero_bds_P_Sigma}. Therefore, 
{\it Step 1:} We first consider the consequences of the condition  $\eta\leq \min\{I_1,I_2\}$. Straightforward calculations show that when condition $\eta\leq I_1$ is satisfied, the following inequalities hold:
\begin{enumerate}
    \item $\forall i=1,\cdots,N$,
    \begin{eqnarray}\label{eqn:I1_conseq_3}
    \|K_t^{i\prime}-K_t^i\| = \eta \|\nabla_{K_t^i} C^i(\pmb{K})\Sigma_t^{-1}\| \leq \frac{\sqmin\sx}{20 T \gamma_B C^{i,-i*}}.    
    \end{eqnarray}
    \item $\forall i=1,\cdots,N$, \begin{eqnarray}\label{eqn:I1_conseq_1}
    &&\eta \left(\frac{ \rho_{\pmb{K}}^{2T}-1}{\rho_{\pmb{K}}^2-1} \left(\frac{C^{i,-i*}}{\sqmin}+T\|W\|\right)2\rho_{\pmb{K}}\,\frac{\gamma_B}{\sx}\sum_{t=0}^{T-1}\|\nabla_{K_t^i}C^i(\pmb{K})\|\right)\nonumber\\
    &\leq& I_1 \left(\frac{ \rho_{\pmb{K}}^{2T}-1}{\rho_{\pmb{K}}^2-1} \left(\frac{\sum_{i=1}^N C^{i,-i*}}{\sqmin}+T\|W\|\right)2\rho_{\pmb{K}}\,\frac{\gamma_B}{\sx}T\max_t\{\|\nabla_{K_t^i}C^i(\pmb{K})\|\}\right)\nonumber \\
    &\leq& \frac{\sx}{10}.
\end{eqnarray}
\item $\forall i=1,\cdots,N$, \begin{eqnarray}\label{eqn:I1_conseq_2}
\eta\leq I_1&\leq& \frac{\sx}{\sx+2\frac{2C^{i,-i*}}{\sqmin}(\gamma_R+\gamma_B^2\frac{C^{i,-i*}}{\sx})}.
\end{eqnarray}
\end{enumerate}
In the case where  $\eta\leq I_2$, we have  $\forall i=1,\cdots,N$
\begin{eqnarray}\label{eqn:I2_conseq}
&&4\eta\,k_i\frac{10 TC^{i,-i*}}{(10 T-1)\sqmin}\left(\gamma_R+(\gamma_B)^2\frac{C^{i,-i*}}{\sx}\right)+8\eta\frac{d}{\sx}\left(\frac{10 TC^{i,-i*}}{(10 T-1)\sqmin}\right)^2\left(\gamma_R+(\gamma_B)^2\frac{C^{i,-i*}}{\sx}\right)^2\nonumber\\ 
&\leq& 4 I_2 \left(\,k_i\frac{10 TC^{i,-i*}}{(10 T-1)\sqmin}\left(\gamma_R+(\gamma_B)^2\frac{C^{i,-i*}}{\sx}\right)+2\frac{d}{\sx}\left(\frac{10 TC^{i,-i*}}{(10 T-1)\sqmin}\right)^2\left(\gamma_R+(\gamma_B)^2\frac{C^{i,-i*}}{\sx}\right)^2\right)\nonumber\\
&\leq& \frac{d}{20\sx} \left(\frac{10 T\min_i\{C^{i*}\}}{(10 T-1)\sqmin}\right)^2 \leq \frac{d}{20\sx} \left(\frac{10 T\max_i\{C^{i,-i*}\}}{(10 T-1)\sqmin}\right)^2.
\end{eqnarray}
By \eqref{eqn:I1_conseq_3} and Lemma~\ref{lemma:nonzero_bds_P_Sigma} we have
\begin{equation*}
    \gamma_B\|K_t^{i\prime}-K_t^i\|  \leq \frac{\sqmin\sx}{20 TC^{i,-i*}} \leq \frac{1}{20 T^2}.
\end{equation*}
%The last inequality holds since $\sx\leq \frac{C^{i,-i*}}{T\sqmin}$ by Lemma \ref{lemma:nonzero_bds_P_Sigma}. 
Therefore, we have $\rho_{\pmb{K},\pmb{K}^{\prime}}\leq \rho_{\pmb{K}}$ by Lemma \ref{lemma:rho_upper_bd}. \\[.1in]
%where $\pmb{K}^\prime$ is defined in \eqref{eqn:one_step_update}. 
{\it Step 2:} We bound the norm of the state covariance matrix.
By Lemma \ref{lemma:bd_F_K1K2},
\[
\sum_{t=0}^{T-1}\|\mathcal{F}_{K_t^i,K_t^{-i*}}-\mathcal{F}_{K^{i\prime}_t,K_t^{-i*}}\| \leq 2\rho_{\pmb{K}}\gamma_B\left(\sum_{t=0}^{T-1}\|K_t^i-K_t^{i\prime}\|\right),
\] 
and hence by Lemma \ref{lemma:perturb_Sigma_K1K2}, we have
\begin{eqnarray}
  \Big\|\Sigma_{\pmb{K}^i,\pmb{K}^{-i*}}-\Sigma_{\pmb{K}^{i\prime},\pmb{K}^{-i*}}\Big\| 
     & \leq & \frac{\rho_{\pmb{K}}^{2T}-1}{\rho_{\pmb{K}}^2-1} \left( \sum_{t=0}^{T-1}\|\mathcal{F}_{K_t^i,K_t^{-i*}}-\mathcal{F}_{K^{i\prime}_t,K_t^{-i*}}\|\right)\left(\|\Sigma_0\|+T\|W\|\right) \nonumber \\
     & \leq & {\frac{ \rho_{\pmb{K}}^{2T}-1}{\rho_{\pmb{K}}^2-1} \left(\frac{C^{i,-i*}}{\sqmin}+T\|W\|\right)2\rho_{\pmb{K}}\gamma_B\,\vertiii{\pmb{K}^i-\pmb{K}^{i\prime}}} \nonumber \\
       & \leq & \frac{ \rho_{\pmb{K}}^{2T}-1}{\rho_{\pmb{K}}^2-1} \left(\frac{C^{i,-i*}}{\sqmin}+T\|W\|\right)2\rho_{\pmb{K}}\gamma_B\,\frac{\eta}{\sx}\sum_{t=0}^{T-1}\|\nabla_{K_t^i}C^i(\pmb{K})\|
       \nonumber \\
       &\leq & \frac{\sx}{10}, \label{eqn:one_step_inte_SigDiff}
\end{eqnarray}
where the last inequality holds  by \eqref{eqn:I1_conseq_1} when $\sum_{t=0}^{T-1}\|\nabla_{K_t^i}C^i(\pmb{K}\|>0$. Note that when \\
$\sum_{t=0}^{T-1}\|\nabla_{K_t^i}C^i(\pmb{K})\|=0$, we have $\Big\|\Sigma_{\pmb{K}^i,\pmb{K}^{-i*}}-\Sigma_{\pmb{K}^{i\prime},\pmb{K}^{-i*}}\Big\| =0< \frac{\sx}{10}$ and hence \eqref{eqn:one_step_inte_SigDiff} still holds in this case. Therefore, by Lemma \ref{lemma:nonzero_bds_P_Sigma}, and noting $\sx\leq\|\Sigma_{\pmb{K}^{i\prime},\pmb{K}^{-i*}}\|/T$,
\begin{eqnarray}\label{eqn:bd_SigmaKprimeK}
 \big\|\Sigma_{\pmb{K}^{i\prime},\pmb{K}^{-i*}}\big\|  &\leq& \big\|\Sigma_{\pmb{K}^{i\prime},\pmb{K}^{-i*}}-\Sigma_{\pmb{K}^i,\pmb{K}^{-i*}}\big\| + \big\|\Sigma_{\pmb{K}^i,\pmb{K}^{-i*}}\big\|\leq \frac{\sx}{10} + \frac{C^{i,-i*}}{\sqmin}\nonumber\\
 &\leq& \frac{\big\|\Sigma_{\pmb{K}^{i\prime},\pmb{K}^{-i*}}\big\| }{10 T} + \frac{C^{i,-i*}}{\sqmin},
\end{eqnarray}
which leads to 
\begin{equation}\label{eqn:bd_pertSigma}
   \big\|\Sigma_{\pmb{K}^{i\prime},\pmb{K}^{-i*}}\big\| \leq\frac{10 T\,C^{i,-i*}}{(10 T-1)\sqmin}.
\end{equation}
%Recall that $\gamma_B$ and $\gamma_R$ are defined in \eqref{eqn:defn_gammaAgammaBgammaR}. 
{\it Step 3:} We now bound $\|P_{t,i}^{\pmb{K}}-P_{t,i}^{\pmb{K}^i,\pmb{K}^{-i*}}\|$ by $\sum_{j=1,j\neq i}^N\vertiii{\pmb{K}^j-\pmb{K}^{j*}}$ {where $\pmb{K} = (\pmb{K}^i,\pmb{K}^{-i})$.}
\begin{eqnarray}\label{eqn:diff_P}
 &&\left\|P_{t,i}^{\pmb{K}}-P_{t,i}^{\pmb{K}^i,\pmb{K}^{-i*}}\right\|\nonumber \\
 &=& \left\|\left(A_t - \sum_{j=1}^NB_t^jK_t^{j}\right)^\top P_{t+1,i}^{\pmb{K}}\left(A_t - \sum_{j=1}^NB_t^jK_t^{j}\right)\right.\nonumber\\
&& \left.-\left(A_t - B_t^iK_t^{i} -\sum_{j=1,j\neq i}^N B_t^jK_t^{j*}\right)^\top P_{t+1,i}^{\pmb{K}^i,\pmb{K}^{-i*}}\left(A_t - B_t^iK_t^{i} -\sum_{j=1,j\neq i}^N B_t^jK_t^{j*}\right)\right\|\nonumber\\
% &\leq& \left\|\left(A_t -
% \sum_{j=1}^NB_t^jK_t^{j}\right)^\top \left(P_{t+1,i}^{\pmb{K}}-P_{t+1,i}^{\pmb{K}^i,\pmb{K}^{-i*}}\right)\left(A_t - \sum_{j=1}^NB_t^jK_t^{j}\right)\right\|\\
% && +\left\|\left(\sum_{j=1,j\neq i}^N B_t^j(K_t^{j}-K_t^{j*})\right)^\top P_{t+1,i}^{\pmb{K}^i,\pmb{K}^{-i*}}\left(\sum_{j=1,j\neq i}^N B_t^j(K_t^{j}-K_t^{j*})\right)\right\|\\
% && + 2\left\|\left(A_t - \sum_{j=1}^NB_t^jK_t^{j}\right)^\top  P_{t+1,i}^{\pmb{K}^i,\pmb{K}^{-i*}}\sum_{j=1,j\neq i}^N B_t^j(K_t^{j}-K_t^{j*})\right\|\\
&\leq& \left\|\left(A_t -
\sum_{j=1}^NB_t^jK_t^{j}\right)^\top \left(P_{t+1,i}^{\pmb{K}}-P_{t+1,i}^{\pmb{K}^i,\pmb{K}^{-i*}}\right)\left(A_t - \sum_{j=1}^NB_t^jK_t^{j}\right)\right\|\nonumber\\
&& +\left\|\left(A_t - \sum_{j=1}^NB_t^jK_t^{j}\right)^\top  P_{t+1,i}^{\pmb{K}^i,\pmb{K}^{-i*}}\left(A_t - \sum_{j=1}^NB_t^jK_t^{j}\right)\right.\nonumber\\
&& \left.-\left(A_t - B_t^iK_t^{i} -\sum_{j=1,j\neq i}^N B_t^jK_t^{j*}\right)^\top P_{t+1,i}^{\pmb{K}^i,\pmb{K}^{-i*}}\left(A_t - B_t^iK_t^{i} -\sum_{j=1,j\neq i}^N B_t^jK_t^{j*}\right)\right\|\nonumber\\
&\leq& \rho_{\pmb{K}}^2 \left\|P_{t+1,i}^{\pmb{K}}-P_{t+1,i}^{\pmb{K}^i,\pmb{K}^{-i*}}\right\|+2\rho_{\pmb{K}}\gamma_B\left\|P_{t+1,i}^{\pmb{K}^i,\pmb{K}^{-i*}}\right\|\sum_{j=1,j\neq i}^N\left\|K_t^{j}-K_t^{j*}\right\|,
% &\leq& \rho^2 \left\|P_{t+1,i}^{\pmb{K}}-P_{t+1,i}^{\pmb{K}^i,\pmb{K}^{-i*}}\right\|+2\rho\gamma_B\left\|P_{t+1,i}^{\pmb{K}^i,\pmb{K}^{-i*}}\right\|\sum_{j=1,j\neq i}^N\left\|K_t^{j}-K_t^{j*}\right\|+(\gamma_B)^2\left\|P_{t+1,i}^{\pmb{K}^i,\pmb{K}^{-i*}}\right\|\left(\sum_{j=1,j\neq i}^N\left\|K_t^{j}-K_t^{j*}\right\|\right)^2\\
% &\leq& \rho^2 \left\|P_{t+1,i}^{\pmb{K}}-P_{t+1,i}^{\pmb{K}^i,\pmb{K}^{-i*}}\right\|+2\rho\gamma_B\frac{C^{i,-i*}}{\sx}\sum_{j=1,j\neq i}^N\left\|K_t^{j}-K_t^{j*}\right\|+(\gamma_B)^2\frac{C^{i,-i*}}{\sx}\left(\sum_{j=1,j\neq i}^N\left\|K_t^{j}-K_t^{j*}\right\|\right)^2,
\end{eqnarray}
where the last inequality holds by letting $X=P_{t+1,i}^{\pmb{K}^i,\pmb{K}^{-i*}}$, $Y=A_t -
\sum_{j=1}^NB_t^jK_t^{j}$, and $Z=A_t - B_t^iK_t^{i} -\sum_{j=1,j\neq i}^N B_t^jK_t^{j*}$ in \eqref{eqn:diff_YXY}.
Since $\left\|P_{T,i}^{\pmb{K}}-P_{T,i}^{\pmb{K}^i,\pmb{K}^{-i*}}\right\|=0$ holds at terminal time $T$, we obtain $\forall t=0,\cdots, T-1$,
\begin{equation}\label{eqn:diffP_t}
    \begin{split}
         \left\|P_{t,i}^{\pmb{K}}-P_{t,i}^{\pmb{K}^i,\pmb{K}^{-i*}}\right\|
 &\leq \frac{2\gamma_BC^{i,-i*}}{\sx}\frac{\rho_{\pmb{K}}(\rho_{\pmb{K}}^{2(T-t)}-1)}{\rho_{\pmb{K}}^2-1}\sum_{s=t}^{T-1}\left(\sum_{j=1,j\neq i}^N\left\|K_s^{j}-K_s^{j*}\right\|\right).
    \end{split}
\end{equation}
Therefore, $\forall t=0,1,\cdots,T-1$,
\begin{eqnarray}
    && \left\|E_{t,i}^{\pmb{K}}-E_{t,i}^{\pmb{K}^i,\pmb{K}^{-i*}}\right\|\nonumber\\
 &=& \left\|(B_t^i)^\top \left(P_{t+1,i}^{\pmb{K}}-P_{t+1,i}^{\pmb{K}^i,\pmb{K}^{-i*}}\right)\left(A_t-\sum_{j=1}^NB_t^jK_t^j\right)-(B_t^i)^\top P_{t+1,i}^{\pmb{K}^i,\pmb{K}^{-i*}}\sum_{j=1,j\neq i}^NB_t^j(K_t^{j}-K_t^{j*})
 \right\|\nonumber\\
&&\leq \rho_{\pmb{K}}\gamma_B\|P_{t+1,i}^{\pmb{K}}-P_{t+1,i}^{\pmb{K}^i,\pmb{K}^{-i*}}\|+(\gamma_B)^2\frac{C^{i,-i*}}{\sx}\sum_{j=1,j\neq i}^N\|K_t^{j}-K_t^{j*}\|\label{eqn:diffE_inte1}\\
 &&\leq \frac{(\gamma_B)^2C^{i,-i*}}{\sx}\left(\frac{2\rho_{\pmb{K}}^2(\rho_{\pmb{K}}^{2(T-t-1)}-1)}{\rho_{\pmb{K}}^2-1}\sum_{s=t+1}^{T-1}\sum_{j=1,j\neq i}^N\|K_s^{j}-K_s^{j*}\|+\sum_{j=1,j\neq i}^N\|K_t^{j}-K_t^{j*}\|\right)\label{eqn:diffE_inte2}\\
 &&\leq \frac{(\gamma_B)^2C^{i,-i*}}{\sx}\left(\frac{2(\rho_{\pmb{K}}^{2(T-t)}-1)}{\rho_{\pmb{K}}^2-1}\sum_{s=t}^{T-1}\left(\sum_{j=1,j\neq i}^N\|K_s^{j}-K_s^{j*}\|\right)\right)\nonumber\\
 &&\leq \frac{(\gamma_B)^2C^{i,-i*}}{\sx}\left(\frac{2(\rho_{\pmb{K}}^{2T}-1)}{\rho_{\pmb{K}}^2-1}\sum_{j=1,j\neq i}^N\vertiii{\pmb{K}^{j}-\pmb{K}^{j*}}\right)\label{eqn:diff_E_bd},
\end{eqnarray}
where \eqref{eqn:diffE_inte1} holds by \eqref{eqn:defn_Ei}, and \eqref{eqn:diffE_inte2} holds by \eqref{eqn:diffP_t}. \\[.1in]
{\it Step 4:} We can now estimate the cost difference between using $\pmb{K}^i$ and the update $\pmb{K}^{i\prime}$. 
By Lemma~\ref{lemma:nonzero_almost_smoothness} we have
\begin{eqnarray}
     &&C^{i\prime,-i*}-C^{i,-i*}\nonumber\\
     && = \sum_{t=0}^{T-1}\Big[\Tr\big(\Sigma_t^{\pmb{K}^{i\prime},\pmb{K}^{-i*}} (K_t^{i\prime}-K_t^i)^{\top}(R_t^{i}+(B_t^i)^{\top}P_{t +1,i}^{\pmb{K}^{i},\pmb{K}^{-i*}}B_t^i)(K_t^{i\prime}-K_t^i)\big)\nonumber\\
     &&+2\Tr\big(\Sigma_t^{\pmb{K}^{i\prime},\pmb{K}^{-i*}} (K_t^{i\prime}-K_t^i)^{\top}E_{t,i}^{\pmb{K}^i,\pmb{K}^{-i*}}\big)\Big]. \nonumber
     %\label{eqn:Cdiff_inte1}
\end{eqnarray}     
Using the updating rule $K_t^{i\prime} =K_t^i - \eta \nabla_{K_t^i}C(\pmb{K})(\Sigma_t^{\pmb{K}})^{-1}$  and the expression for the gradient from Lemma~\ref{lemma:nonzero_policygrad}
\begin{eqnarray}
%&&C^i(\pmb{K}^{i\prime},\pmb{K}^{-i*})-C^i(\pmb{K}^i,\pmb{K}^{-i*})\nonumber\\
&& C^{i\prime,-i*}-C^{i,-i*} \nonumber \\
&& = \sum_{t=0}^{T-1}\Big[4\eta^2\Tr\Big(\Sigma_t^{\pmb{K}^{i\prime},\pmb{K}^{-i*}} (E_{t,i}^{\pmb{K}})^\top(R_t^{i}+(B_t^i)^{\top}P_{t +1,i}^{\pmb{K}^{i},\pmb{K}^{-i*}}B_t^i)(E_{t,i}^{\pmb{K}})\Big) -4\eta\Tr\big(\Sigma_t^{\pmb{K}^{i\prime},\pmb{K}^{-i*}} (E_{t,i}^{\pmb{K}})^\top E_{t,i}^{\pmb{K}^i,\pmb{K}^{-i*}}\big)\Big] \nonumber\\ 
&& = \sum_{t=0}^{T-1}\Big[4\eta^2\Tr\big(\Sigma_t^{\pmb{K}^{i\prime},\pmb{K}^{-i*}} (E_{t,i}^{\pmb{K}}-E_{t,i}^{\pmb{K}^i,\pmb{K}^{-i*}}+E_{t,i}^{\pmb{K}^i,\pmb{K}^{-i*}})^\top(R_t^{i}+(B_t^i)^{\top}P_{t +1,i}^{\pmb{K}^{i},\pmb{K}^{-i*}}B_t^i)\nonumber\\
&&\quad\quad(E_{t,i}^{\pmb{K}}-E_{t,i}^{\pmb{K}^i,\pmb{K}^{-i*}}+E_{t,i}^{\pmb{K}^i,\pmb{K}^{-i*}})\big) -4\eta\Tr\big(\Sigma_t^{\pmb{K}^{i\prime},\pmb{K}^{-i*}} (E_{t,i}^{\pmb{K}}-E_{t,i}^{\pmb{K}^i,\pmb{K}^{-i*}}+E_{t,i}^{\pmb{K}^i,\pmb{K}^{-i*}})^\top E_{t,i}^{\pmb{K}^i,\pmb{K}^{-i*}}\big)\Big]\nonumber\\ 
&& = \sum_{t=0}^{T-1}\Big[4\eta^2\Tr\big(\Sigma_t^{\pmb{K}^{i\prime},\pmb{K}^{-i*}} (E_{t,i}^{\pmb{K}}-E_{t,i}^{\pmb{K}^i,\pmb{K}^{-i*}})^\top(R_t^{i}+(B_t^i)^{\top}P_{t +1,i}^{\pmb{K}^{i},\pmb{K}^{-i*}}B_t^i)(E_{t,i}^{\pmb{K}}-E_{t,i}^{\pmb{K}^i,\pmb{K}^{-i*}})\big)\nonumber\\
&&\quad\quad+ 8\eta^2\Tr\big(\Sigma_t^{\pmb{K}^{i\prime},\pmb{K}^{-i*}} (E_{t,i}^{\pmb{K}}-E_{t,i}^{\pmb{K}^i,\pmb{K}^{-i*}})^\top(R_t^{i}+(B_t^i)^{\top}P_{t +1,i}^{\pmb{K}^{i},\pmb{K}^{-i*}}B_t^i)E_{t,i}^{\pmb{K}^i,\pmb{K}^{-i*}}\big) \nonumber\\
&& \quad\quad+4\eta^2\Tr\big(\Sigma_t^{\pmb{K}^{i\prime},\pmb{K}^{-i*}} (E_{t,i}^{\pmb{K}^i,\pmb{K}^{-i*}})^\top(R_t^{i}+(B_t^i)^{\top}P_{t +1,i}^{\pmb{K}^{i},\pmb{K}^{-i*}}B_t^i)E_{t,i}^{\pmb{K}^i,\pmb{K}^{-i*}}\big)\nonumber\\
&&\quad\quad-4\eta\Tr\big(\Sigma_t^{\pmb{K}^{i\prime},\pmb{K}^{-i*}} (E_{t,i}^{\pmb{K}}-E_{t,i}^{\pmb{K}^i,\pmb{K}^{-i*}})^\top E_{t,i}^{\pmb{K}^i,\pmb{K}^{-i*}}\big)-4\eta\Tr\big(\Sigma_t^{\pmb{K}^{i\prime},\pmb{K}^{-i*}} (E_{t,i}^{\pmb{K}^i,\pmb{K}^{-i*}})^\top E_{t,i}^{\pmb{K}^i,\pmb{K}^{-i*}}\big)\Big].\nonumber
\end{eqnarray}
     
Now, letting $\omega^2=\frac{2}{\sx}$ in
\[
2\Tr(A^\top B)=\Tr(A^\top B + B^\top A)\leq \omega^2 \Tr(A^\top A)+\frac{1}{\omega^2}\Tr(B^\top B),
\]
(which holds for any matrices $A$ and $B$ of the same dimension) we have
\begin{eqnarray}     
&& C^{i\prime,-i*}-C^{i,-i*} \nonumber \\    
&& \leq \sum_{t=0}^{T-1}\Big[4\eta^2\Tr\big(\Sigma_t^{\pmb{K}^{i\prime},\pmb{K}^{-i*}} (E_{t,i}^{\pmb{K}}-E_{t,i}^{\pmb{K}^i,\pmb{K}^{-i*}})^\top(R_t^{i}+(B_t^i)^{\top}P_{t +1,i}^{\pmb{K}^{i},\pmb{K}^{-i*}}B_t^i)(E_{t,i}^{\pmb{K}}-E_{t,i}^{\pmb{K}^i,\pmb{K}^{-i*}})\big) \nonumber \\
&&\quad\quad+ 8\eta^2\frac{\sx}{4}\Tr\big((E_{t,i}^{\pmb{K}^i,\pmb{K}^{-i*}})^\top E_{t,i}^{\pmb{K}^i,\pmb{K}^{-i*}}\big)+ 8\eta^2\frac{1}{\sx}\Tr\big(\Sigma_t^{\pmb{K}^{i\prime},\pmb{K}^{-i*}} (E_{t,i}^{\pmb{K}}-E_{t,i}^{\pmb{K}^i,\pmb{K}^{-i*}})^\top\cdot \nonumber\\
&&\qquad\,\,\quad(R_t^{i}+(B_t^i)^{\top}P_{t +1,i}^{\pmb{K}^{i},\pmb{K}^{-i*}}B_t^i)(R_t^{i}+(B_t^i)^{\top}P_{t +1,i}^{\pmb{K}^{i},\pmb{K}^{-i*}}B_t^i)(E_{t,i}^{\pmb{K}}-E_{t,i}^{\pmb{K}^i,\pmb{K}^{-i*}})\Sigma_t^{\pmb{K}^{i\prime},\pmb{K}^{-i*}}\big) \nonumber\\
&& \quad\quad+4\eta^2\Tr\big(\Sigma_t^{\pmb{K}^{i\prime},\pmb{K}^{-i*}} (E_{t,i}^{\pmb{K}^i,\pmb{K}^{-i*}})^\top(R_t^{i}+(B_t^i)^{\top}P_{t +1,i}^{\pmb{K}^{i},\pmb{K}^{-i*}}B_t^i)E_{t,i}^{\pmb{K}^i,\pmb{K}^{-i*}}\big)\nonumber\\
&&\quad\quad+4\eta\frac{1}{\sx}\Tr\big(\Sigma_t^{\pmb{K}^{i\prime},\pmb{K}^{-i*}} (E_{t,i}^{\pmb{K}}-E_{t,i}^{\pmb{K}^i,\pmb{K}^{-i*}})^\top (E_{t,i}^{\pmb{K}}-E_{t,i}^{\pmb{K}^i,\pmb{K}^{-i*}})\Sigma_t^{\pmb{K}^{i\prime},\pmb{K}^{-i*}}\big)\nonumber\\
&&\quad\quad +4\eta\frac{\sx}{4}\Tr\big( (E_{t,i}^{\pmb{K}^i,\pmb{K}^{-i*}})^\top E_{t,i}^{\pmb{K}^i,\pmb{K}^{-i*}}\big)-4\eta\sx\Tr\big( (E_{t,i}^{\pmb{K}^i,\pmb{K}^{-i*}})^\top E_{t,i}^{\pmb{K}^i,\pmb{K}^{-i*}}\big)\Big]\nonumber \\
&& \leq \sum_{t=0}^{T-1}\Big[\Big(4\eta^2\|\Sigma_t^{\pmb{K}^{i\prime},\pmb{K}^{-i*}}\|\Tr\big(R_t^{i}+(B_t^i)^{\top}P_{t +1,i}^{\pmb{K}^{i},\pmb{K}^{-i*}}B_t^i\big)+ \frac{8\eta^2}{\sx}\|R_t^{i}+(B_t^i)^{\top}P_{t +1,i}^{\pmb{K}^{i},\pmb{K}^{-i*}}B_t^i\|^2\cdot\nonumber\\
&&\quad\quad\Tr\big(\Sigma_t^{\pmb{K}^{i\prime},\pmb{K}^{-i*}} \Sigma_t^{\pmb{K}^{i\prime},\pmb{K}^{-i*}}\big)+ \frac{4\eta}{\sx}\Tr\big(\Sigma_t^{\pmb{K}^{i\prime},\pmb{K}^{-i*}} \Sigma_t^{\pmb{K}^{i\prime},\pmb{K}^{-i*}}\big)\Big)\|E_{t,i}^{\pmb{K}}-E_{t,i}^{\pmb{K}^i,\pmb{K}^{-i*}}\|^2 \nonumber\\
&& \quad\quad  +\big(2\eta^2\sx+4\eta^2\|\Sigma_t^{\pmb{K}^{i\prime},\pmb{K}^{-i*}}\|\|R_t^{i}+(B_t^i)^{\top}P_{t +1,i}^{\pmb{K}^{i},\pmb{K}^{-i*}}B_t^i\|+\eta\sx-4\eta\sx\big)\cdot\nonumber\\
&&\quad\quad \Tr\big( (E_{t,i}^{\pmb{K}^i,\pmb{K}^{-i*}})^\top E_{t,i}^{\pmb{K}^i,\pmb{K}^{-i*}}\big)\Big].\label{eqn:Cdiff_inte3}
\end{eqnarray}
     Now, using $\big\|\Sigma_{\pmb{K}^{i\prime},\pmb{K}^{-i*}}\big\| \leq\frac{2\,C^{i,-i*}}{\sqmin}$ by \eqref{eqn:bd_pertSigma} (this is a loose approximation in order to ease the analysis and smooth out the presentation), we can bound the step size condition in \eqref{eqn:I1_conseq_2} by
 \begin{equation*}
     \eta \leq \frac{\sx}{\sx+2\|\Sigma_t^{\pmb{K}^{i\prime},\pmb{K}^{-i*}}\|\|R_t^{i}+(B_t^i)^{\top}P_{t +1,i}^{\pmb{K}^{i},\pmb{K}^{-i*}}B_t^i\|}.
 \end{equation*}    
This gives 
\[
2\eta^2\sx+4\eta^2\left\|\Sigma_t^{\pmb{K}^{i\prime},\pmb{K}^{-i*}}\right\|\left\|R_t^{i}+(B_t^i)^{\top}P_{t +1,i}^{\pmb{K}^{i},\pmb{K}^{-i*}}B_t^i\right\|+\eta\sx-4\eta\sx \leq -\eta\sx.
\]    
Hence, using this in \eqref{eqn:Cdiff_inte3}, we have
     \begin{equation}
C^{i\prime,-i*}-C^{i,-i*} \leq \eta\,h_{\rm diff}^i\sum_{t=0}^{T-1}\|E_{t,i}^{\pmb{K}}-E_{t,i}^{\pmb{K}^i,\pmb{K}^{-i*}}\|^2  -\eta\sx \sum_{t=0}^{T-1}\Tr\big( (E_{t,i}^{\pmb{K}^i,\pmb{K}^{-i*}})^\top E_{t,i}^{\pmb{K}^i,\pmb{K}^{-i*}}\big)\label{eqn:onestep_inte},
\end{equation}
where 
\begin{equation*}
\begin{split}
 h_{\rm diff}^i &:= 4\eta\,k_i\frac{10 TC^{i,-i*}}{(10 T-1)\sqmin}\left(\gamma_R+(\gamma_B)^2\frac{C^{i,-i*}}{\sx}\right)+8\eta\frac{d}{\sx}\left(\frac{10 TC^{i,-i*}}{(10 T-1)\sqmin}\right)^2\left(\gamma_R+(\gamma_B)^2\frac{C^{i,-i*}}{\sx}\right)^2\\
 &\qquad +4\frac{d}{\sx} \left(\frac{10 TC^{i,-i*}}{(10 T-1)\sqmin}\right)^2.
 \end{split}
\end{equation*}
%Here \eqref{eqn:Cdiff_inte1} holds by Lemma ; \eqref{eqn:Cdiff_inte2} holds by \eqref{eqn:one_step_update} and Lemma ;
%\eqref{eqn:Cdiff_inte3} holds by  \eqref{eqn:onestep_inte} holds since, by step size condition \eqref{eqn:step_size},we have 
%by , which leads to

Therefore, by \eqref{eqn:bd_pertSigma}, \eqref{eqn:diff_E_bd}, Lemma \ref{lemma:grad_domi}, and Lemma \ref{lemma:nonzero_bds_P_Sigma},
\begin{eqnarray}
 C^{i\prime,-i*}-C^{i,-i*} &\leq & \eta\,h_{\rm diff}^i\,T\,\left[\frac{(\gamma_B)^2C^{i,-i*}}{\sx}\frac{2(\rho_{\pmb{K}}^{2T}-1)}{\rho_{\pmb{K}}^2-1}\sum_{j=1,j\neq i}^N\vertiii{\pmb{K}^{j}-\pmb{K}^{j*}}\right]^2 \nonumber \\ &&\quad-\eta\sx\sum_{t=0}^{T-1}\Tr\big( (E_{t,i}^{\pmb{K}^i,\pmb{K}^{-i*}})^\top E_{t,i}^{\pmb{K}^i,\pmb{K}^{-i*}}\big) \nonumber \\
     &\leq & \eta\,h_{\rm glob}\left(\sum_{j=1,j\neq i}^N\vertiii{\pmb{K}^{j}-\pmb{K}^{j*}}\right)^2-\eta\frac{\sx\srmin}{\|\Sigma_{\pmb{K}^{*}}\|}\big(C^{i,-i*}-C^{i*}\big), \label{eqn:ind_cost_diff}
\end{eqnarray}
where 
\begin{equation}\label{eqn:defn_hloc}
 \begin{split}
      h_{\rm glob} &= 4T\left[\eta\,k_i\frac{10 T\max_i\{C^{i,-i*}\}}{(10 T-1)\sqmin}\left(\gamma_R+(\gamma_B)^2\frac{\max_i\{C^{i,-i*}\}}{\sx}\right)+2\eta\frac{d}{\sx}\left(\frac{10 T\max_i\{C^{i,-i*}\}}{(10 T-1)\sqmin}\right)^2\cdot\right.\\
      &\qquad\left.\left(\gamma_R+(\gamma_B)^2\frac{\max_i\{C^{i,-i*}\}}{\sx}\right)^2+\,\frac{d}{\sx} \left(\frac{10 T\max_i\{C^{i,-i*}\}}{(10 T-1)\sqmin}\right)^2\right] \cdot\\
      &\qquad\left[\frac{(\gamma_B)^2\max_i\{C^{i,-i*}\}}{\sx}\frac{2(\rho_{\pmb{K}}^{2T}-1)}{\rho_{\pmb{K}}^2-1}\right]^2.
 \end{split}
\end{equation}
{\it Step 5:}
Finally we can establish the one-step contraction. Using \eqref{eqn:ind_cost_diff}, we have
\begin{eqnarray}
 C^{i\prime,-i*}-C^{i*} &=& C^{i\prime,-i*}-C^{i,-i*}+C^{i,-i*}-C^{i*}\nonumber\\
 &\leq& \left(1-\eta\frac{\sx\srmin}{\|\Sigma_{\pmb{K}^{*}}\|}\right)\big(C^{i,-i*}-C^{i*}\big) + \eta\,h_{\rm glob}\left(\sum_{j=1,j\neq i}^N\vertiii{\pmb{K}^{j}-\pmb{K}^{j*}}\right)^2.\label{eq:C_diff}
\end{eqnarray}
Hence by Lemma \ref{lemma:nonzero_almost_smoothness} and \eqref{eqn:Ci_diff_Ki}, we have
\begin{eqnarray*}
\sum_{j=1,j\neq i}^N\left(C^{j,-j*} - C^{j*}\right) 
 \geq \frac{\sx\srmin}{T} \sum_{j=1,j\neq i}^N\vertiii{\pmb{K}^j-\pmb{K}^{j*}}^2 \geq \frac{\sx\srmin}{T(N-1)} \left(\sum_{j=1,j\neq i}^N\vertiii{\pmb{K}^j-\pmb{K}^{j*}}\right)^2,
\end{eqnarray*}
and thus
\begin{equation}
    \begin{split}
        C^{i\prime,-i*}-C^{i*}& \leq \big(1-\eta\frac{\sx\srmin}{\|\Sigma_{\pmb{K}^{*}}\|}\big)\big(C^{i,-i*}-C^{i*}\big) + \,\eta\,\, h_{\rm glob}\frac{T(N-1)}{\sx \srmin}\left(\sum_{j=1,j\neq i}^N (C^{j,-j*} - C^{j*}) \right).\label{eq:C_diff1}
    \end{split}
\end{equation}
Summing up \eqref{eq:C_diff1} for $i=1,\cdots,N$, we have
\begin{equation}\label{eq:c_sum_onestep}
\sum_{i=1}^N \left(C^{i\prime,-i*}-C^{i*}\right)
\leq \left(1-\eta\frac{\sx\srmin}{\|\Sigma_{\pmb{K}^{*}}\|}+\eta(N-1)h_{\rm glob}\frac{T(N-1)}{\sx \srmin}\right)\left(\sum_{i=1}^N (C^{i,-i*} - C^{i*})\right).
\end{equation}
Since $\eta\leq I_2$, we have \eqref{eqn:I2_conseq} and
% \begin{eqnarray*}
% \eta \leq \left(k_i\frac{\kappa T C^{i,-i*}}{(\kappa T-1)\sqmin}\left(\gamma_R+(\gamma_B)^2\frac{ C^{i,-i*}}{\sx}\right)+\frac{2d}{\sx}\left(\frac{\kappa T C^{i,-i*}}{(\kappa T-1)\sqmin}\right)^2\left(\gamma_R+(\gamma_B)^2\frac{C^{i,-i*}}{\sx}\right)^2\right)^{-1}\frac{d}{4\sx} \left(\frac{\kappa T  C^{i,-i*}}{(\kappa T-1)\sqmin}\right)^2.
%  \end{eqnarray*}
then
 \begin{eqnarray*}
  h_{\rm diff}^i \leq (4+\frac{1}{20})\frac{d}{\sx} \left(\frac{10 T \max_i\{C^{i,-i*}\}}{(10 T-1)\sqmin}\right)^2\leq 5\frac{d}{\sx} \left(\frac{ \max_i\{C^{i,-i*}\}}{\sqmin}\right)^2,\quad\text{and}\quad h_{\rm glob} \leq \bar{h}_{\rm glob},
 \end{eqnarray*}
where $\bar{h}_{\rm glob}$ is given by
\begin{equation}\label{eqn:defn_hloc_bar}
  \begin{split}
      \bar{h}_{\rm glob} &= 5\,T\,d\, \frac{\max_i\{C^{i,-i*}\}^4(\gamma_B)^4}{\sqmin^2\sx^3}\left[\frac{2(\rho_{\pmb{K}}^{2T}-1)}{\rho_{\pmb{K}}^2-1}\right]^2.
  \end{split}
 \end{equation}
Under condition \eqref{eqn:one_step_noise_cond}, we have
\begin{eqnarray*}
\sx g_1 - \frac{g_2}{\sx^4} = \frac{\sx\srmin}{\|\Sigma_{\pmb{K}^{*}}\|}-(N-1)^2\bar{h}_{\rm glob}\frac{T}{\sx\srmin}>0,
 \end{eqnarray*}
which indicates that $\alpha\eta>0$. Since $\eta\leq \frac{1}{\srmin}$, we have
 \begin{eqnarray*}
 \eta < \frac{\|\Sigma_{\pmb{K}^{*}}\|}{\sx\srmin}<\left( \frac{\sx\srmin}{\|\Sigma_{\pmb{K}^{*}}\|}-(N-1)^2\bar{h}_{\rm glob}\frac{T}{\sx\srmin}\right)^{-1} = \left(\sx g_1 - \frac{g_2}{\sx^4}\right)^{-1}.
 \end{eqnarray*}
Recall that in the statement we define $\alpha=\sx g_1-g_2/\sx^4$. Therefore we have $\alpha\eta<1$. Along with \eqref{eq:c_sum_onestep}, we obtain the one-step contraction \eqref{eq:c_sum_onestep_final}.
\end{proof}

\begin{Lemma}\label{lemma:grad_bound}
Assume Assumptions \ref{ass:nonzero_cost},  \ref{ass:nonzero_initial_noise}, and \ref{ass:nonzero_exist_sol} hold. Then we have that for player $i$,
\begin{eqnarray}
\sum_{t=0}^{T-1} \|\nabla_{K_t^i}C^i(\pmb{K})\|^2 &\leq& 8\left\{\frac{\rho_{\pmb{K}}^{2(T+1)}-1}{\rho_{\pmb{K}}^2-1}\|\Sigma_0\|+\frac{\rho_{\pmb{K}}^{2T}-1}{\rho_{\pmb{K}}^2-1}\,T\,\|W\|\right\}^2\cdot\nonumber\\
&&\left\{d\frac{T^2\,(N-1)}{\sx\srmin}\left[\frac{(\gamma_B)^2\sum_{i=1}^NC^i(\pmb{K}^i,\pmb{K}^{-i*})}{\sx}\frac{2(\rho_{\pmb{K}}^{2T}-1)}{\rho_{\pmb{K}}^2-1}\right]^2\cdot\right.\nonumber\\
&& \quad\Big(\sum_{j=1,j\neq i}^N(C^j(\pmb{K}^{j},\pmb{K}^{-j*}) - C^j(\pmb{K}^{*}))\Big)\nonumber\\
&&\left.+\frac{\sx\gamma_R+(\gamma_B)^2\sum_{i=1}^NC^i(\pmb{K}^i,\pmb{K}^{-i*})}{\sx^2}\left(C^i(\pmb{K}^i,\pmb{K}^{-i*})-C^i(\pmb{K}^{*})\right)
\right\}.\label{eqn:grad_bound_lemma_state}
\end{eqnarray}
%{\color{blue}where $\rho_{\pmb{K}}$ is defined in \eqref{defn_rho_K}}.
\end{Lemma}

\begin{proof}
Using Lemma \ref{lemma:nonzero_policygrad}, we have
\begin{equation}\label{eqn:grad_sq_bd}
    \sum_{t=0}^{T-1} \|\nabla_{K_t^i}C^i(\pmb{K})\|^2 \leq 4\sum_{t=0}^{T-1} \Tr\left(\Sigma_t^{\pmb{K}}(E_{t,i}^{\pmb{K}})^\top E_{t,i}^{\pmb{K}}\Sigma_t^{\pmb{K}}\right)\leq   4\left(\|\Sigma_{\pmb{K}}\|\right)^2\sum_{t=0}^{T-1} \Tr\left((E_{t,i}^{\pmb{K}})^\top E_{t,i}^{\pmb{K}}\right),
\end{equation}
and
\begin{eqnarray}
\sum_{t=0}^{T-1} \Tr\left((E_{t,i}^{\pmb{K}})^\top E_{t,i}^{\pmb{K}}\right) &=& \sum_{t=0}^{T-1} \Tr\left((E_{t,i}^{\pmb{K}}-E_{t,i}^{\pmb{K}^i,\pmb{K}^{-i*}}+E_{t,i}^{\pmb{K}^i,\pmb{K}^{-i*}})^\top (E_{t,i}^{\pmb{K}}-E_{t,i}^{\pmb{K}^i,\pmb{K}^{-i*}}+E_{t,i}^{\pmb{K}^i,\pmb{K}^{-i*}})\right)\nonumber\\
&\leq& 2\sum_{t=0}^{T-1} \Tr\left((E_{t,i}^{\pmb{K}}-E_{t,i}^{\pmb{K}^i,\pmb{K}^{-i*}})^\top (E_{t,i}^{\pmb{K}}-E_{t,i}^{\pmb{K}^i,\pmb{K}^{-i*}})\right)  +2\sum_{t=0}^{T-1} \Tr\left((E_{t,i}^{\pmb{K}^i,\pmb{K}^{-i*}})^\top E_{t,i}^{\pmb{K}^i,\pmb{K}^{-i*}}\right)\nonumber\\
&\leq& 2d\sum_{t=0}^{T-1} \left\|E_{t,i}^{\pmb{K}}-E_{t,i}^{\pmb{K}^i,\pmb{K}^{-i*}}\right\|^2 +2\sum_{t=0}^{T-1} \Tr\left((E_{t,i}^{\pmb{K}^i,\pmb{K}^{-i*}})^\top E_{t,i}^{\pmb{K}^i,\pmb{K}^{-i*}}\right).\label{eqn:grad_bound}
\end{eqnarray}
By \eqref{eqn:Ci_diff_Ki}, \eqref{eqn:diff_E_bd}, we have
\begin{eqnarray}
    &&\sum_{t=0}^{T-1} \left\|E_{t,i}^{\pmb{K}}-E_{t,i}^{\pmb{K}^i,\pmb{K}^{-i*}}\right\|^2\nonumber\\
    &\leq& T\left[\frac{(\gamma_B)^2C^i(\pmb{K}^i,\pmb{K}^{-i*})}{\sx}\frac{2(\rho_{\pmb{K}}^{2T}-1)}{\rho_{\pmb{K}}^2-1}\left(\sum_{j=1,j\neq i}^N\vertiii{\pmb{K}^{j}-\pmb{K}^{j*}}\right)\right]^2\nonumber\\
    &\leq& \frac{T^2\,(N-1)}{\sx\srmin}\left[\frac{(\gamma_B)^2C^i(\pmb{K}^i,\pmb{K}^{-i*})}{\sx}\frac{2(\rho_{\pmb{K}}^{2T}-1)}{\rho_{\pmb{K}}^2-1}\right]^2\Big(\sum_{j=1,j\neq i}^N(C^j(\pmb{K}^{j},\pmb{K}^{-j*}) - C^j(\pmb{K}^{*}))\Big)
    \label{eqn:grad_bound_1}.
\end{eqnarray}
By Lemma \ref{lemma:grad_domi} and Lemma \ref{lemma:nonzero_bds_P_Sigma}, we have
\begin{equation}\label{eqn:grad_bound_2}
\begin{split}
    \sum_{t=0}^{T-1}\Tr\big((E_{t,i}^{\pmb{K}^i,\pmb{K}^{-i*}})^\top E_{t,i}^{\pmb{K}^i,\pmb{K}^{-i*}}\big)&\leq\frac{\sx\gamma_R+(\gamma_B)^2C^i(\pmb{K}^i,\pmb{K}^{-i*})}{\sx^2} \left(C^i(\pmb{K}^i,\pmb{K}^{-i*})-C^i(\pmb{K}^{*})\right).
      \end{split}
\end{equation}
By Proposition \ref{prop:Sigma_Gamma_relation},
\begin{eqnarray}\label{eqn:grad_bound_3}
\|\Sigma_{\pmb{K}}\| &\leq& \|\Sigma_0\|
+\sum_{t=0}^{T-1}\|\mathcal{G}_t(\Sigma_0)\| +\sum_{t=1}^{T-1}\sum_{s=1}^t \|D_{t,s}WD_{t,s}^{\top}\| + T\|W\|\nonumber\\
&\leq& \|\Sigma_0\|
+\|\Sigma_0\|\sum_{t=0}^{T-1}\rho_{\pmb{K}}^{2(t+1)} +T\|W\|\sum_{s=1}^{T-1} \rho_{\pmb{K}}^{2(T-s)}+T\|W\|\nonumber\\
&\leq& \frac{\rho_{\pmb{K}}^{2(T+1)}-1}{\rho_{\pmb{K}}^2-1}\|\Sigma_0\|+\frac{\rho_{\pmb{K}}^{2T}-1}{\rho_{\pmb{K}}^2-1}\,T\,\|W\|.
\end{eqnarray}
Note that here we use $\sum_{t=1}^{T-1}\sum_{s=1}^t \rho_{\pmb{K}}^{2(t-s)}\leq T\sum_{s=1}^{T-1} \rho_{\pmb{K}}^{2(T-s)}$, which is a loose bound in order to simplify the presentation.
Therefore, combining \eqref{eqn:grad_sq_bd}-\eqref{eqn:grad_bound_3}, we obtain the statement \eqref{eqn:grad_bound_lemma_state}.
\end{proof}

Finally, we are ready to provide the proof for the main result based on Lemmas \ref{lemma:one_step_loc_conv} and \ref{lemma:grad_bound}.
\begin{proof}[Proof of Theorem \ref{thm:local_conv_NPG}]
We first show that the total cost for the $N$ players decreases at round $m=1$. Take $\pmb{K}=\pmb{K}^{(0)}$ and $\pmb{K}^{\prime}=\pmb{K}^{(1)}$ in Lemma \ref{lemma:one_step_loc_conv}, we begin by showing that there exists a positive lower bound on the RHS of \eqref{eqn:step_size}. 
%Denote $C^{i,-i*}=C^i(\pmb{K}^{i,(0)},\pmb{K}^{-i*})$ and $C^{i*}=C^i(\pmb{K}^*)$. 
By the Cauchy-Schwarz inequality and Lemma \ref{lemma:grad_bound},  
\begin{eqnarray*}
&&\sum_{t=0}^{T-1} \|\nabla_{K_t^i}C^i(\pmb{K}^{(0)})\| \leq\sqrt{T\cdot\sum_{t=0}^{T-1} \|\nabla_{K_t^i}C^i(\pmb{K}^{(0)})\|^2}\\
&\leq &2\sqrt{2}\left(\frac{\rho_{\pmb{K}^{(0)}}^{2(T+1)}-1}{\rho_{\pmb{K}^{(0)}}^2-1}\|\Sigma_0\|+\frac{\rho_{\pmb{K}^{(0)}}^{2T}-1}{\rho_{\pmb{K}^{(0)}}^2-1}\,T\,\|W\|\right)\cdot\left\{d\frac{T^2\,(N-1)}{\sx\srmin}\left[\frac{(\gamma_B)^2\sum_{i=1}^NC^{i,-i*}}{\sx}\cdot\right.\right.\\
&&\left.\left.\frac{2(\rho_{\pmb{K}^{(0)}}^{2T}-1)}{\rho_{\pmb{K}^{(0)}}^2-1}\right]^2+\frac{\sx\gamma_R+(\gamma_B)^2\sum_{i=1}^NC^{i,-i*}}{\sx^2}\right\}^{1/2}\sqrt{T\left(\sum_{i=1}^N(C^{i,-i*} - C^{i*})\right)}\\
&\leq &2\sqrt{2}\left(\frac{\bar{\rho}^{2(T+1)}-1}{\bar{\rho}^2-1}\|\Sigma_0\|+\frac{\bar{\rho}^{2T}-1}{\bar{\rho}^2-1}\,T\,\|W\|\right)\cdot\left\{d\frac{T^2\,(N-1)}{\sx\srmin}\left[\frac{(\gamma_B)^2\sum_{i=1}^NC^{i,-i*}}{\sx}\cdot\right.\right.\\
&&\left.\left.\frac{2(\bar{\rho}^{2T}-1)}{\bar{\rho}^2-1}\right]^2+\frac{\sx\gamma_R+(\gamma_B)^2\sum_{i=1}^NC^{i,-i*}}{\sx^2}\right\}^{1/2}\sqrt{T\left(\sum_{i=1}^N(C^{i,-i*} - C^{i*})\right)},
\end{eqnarray*} 
where the last inequality holds since, after performing one-step natural policy gradient ${K}_t^{i,(1)} = {K}_t^{i,(0)} - \eta \nabla_{K_t^i}C^i(\pmb{K}^{(0)})(\Sigma_t^{\pmb{K}^{(0)}})^{-1}$, we have 
\begin{eqnarray}
\rho_{\pmb{K}^{(0)}} &\leq& \rho^* +N\gamma_B\sqrt{\frac{T}{\sx\srmin} \psi}+\frac{1}{20 T^2} = \bar{\rho}  \label{eq:inter1}
\end{eqnarray}
where $\rho^*$ and $\bar{\rho}$ are defined in \eqref{eqn:defn_rho_star} and \eqref{eqn:defn_rho_bar}, and $\psi := \max_i\{C^i(\pmb{K}^{i,(0)},\pmb{K}^{-i*})-C^{i*}\}$.
\eqref{eq:inter1} holds by Lemma \ref{lemma:rho_upper_bd}.

Now we aim to show that $\frac{1}{\bar{\rho}}$ is bounded below by polynomials in some model parameters. 
%This is equivalent to showing that $\bar{\rho}$ is bounded above by polynomials in $T$, $N$, $\gamma_A$, $\gamma_B$, $\frac{1}{\sx}$, $\frac{1}{\srmin}$, and $\vertiii{\pmb{K}^*}$, or a constant $1+\xi$.  
Given that $\left\|A_t-\sum_{i=1}^NB_t^iK_t^{i*}\right\|\leq \gamma_A + \gamma_B\vertiii{\pmb{K}^*}$ and that Lemma \ref{lemma:rho_upper_bd}, $\bar{\rho}$ can be bounded above by polynomials in $T$, $N$, $\gamma_A$,  $\gamma_B$, $\frac{1}{\sx}$, $\frac{1}{\srmin}$, and $\vertiii{\pmb{K}^*}$, or a constant $1+\xi$. Therefore, along with the fact that $\frac{1}{d}>\frac{1}{(a+1)(b+1)(c+1)}$ for $d<ab+c$ with some $a,b,c, d>0$ and  $\frac{1}{a^n+1}>\frac{1}{(a+1)^n}$ for $a>0$ and $n\in \mathbb{N}^+$, $\frac{1}{\bar{\rho}}$ is bounded below by polynomials in $\frac{1}{T+1}$, $\frac{1}{N+1}$, $\frac{1}{\gamma_A+1}$, $\frac{1}{\gamma_B+1}$, $\sx$, $\frac{1}{\sx+1}$, $\srmin$, $\frac{1}{\srmin+1}$, and $\frac{1}{\vertiii{\pmb{K}^*}+1}$, or a constant $\frac{1}{1+\xi}$.
Similarly, $I_1$ can be bounded below by polynomials in $\frac{1}{d+1}$, $\frac{1}{N+1}$, $\frac{1}{T+1}$, $\frac{1}{\sum_{i=1}^NC^{i,-i*}+1}$, $\frac{1}{\|W\|+1}$,$\frac{1}{\|\Sigma_0\|+1}$,$\frac{1}{\gamma_A+1}$, $\frac{1}{\gamma_B+1}$,$\frac{1}{\gamma_R+1}$, $\frac{1}{\srmin+1}$, $\srmin$, $\frac{1}{\sqmin+1}$, $\sqmin$, $\frac{1}{\sx+1}$, $\sx$, and $\frac{1}{\vertiii{\pmb{K}^*}+1}$; and $I_2$ can be bounded below by polynomials in $\frac{1}{\sum_ik_i+1}$, $\frac{1}{d+1}$, $d$, $\frac{1}{\sum_{i=1}^NC^{i,-i*}+1}$,  $\frac{1}{\gamma_B+1}$, $\frac{1}{\sqmin+1}$, $\sqmin$, $\frac{1}{\gamma_R+1}$, $\frac{1}{\sx+1}$, $\sx$.

% Hence, by choosing $\eta\in\HH\left(\frac{1}{\sum_{i=1}^NC^i(\pmb{K}^{i,(0)},\pmb{K}^{-i*})+1}\right)$ as an appropriate polynomial in $\frac{1}{\sum_{i=1}^NC^i(\pmb{K}^{i,(0)},\pmb{K}^{-i*})+1}$, $\frac{1}{\sum_i k_i+1}$, $\frac{1}{d+1}$, $d$, $\frac{1}{N+1}$, $\frac{1}{T+1}$, $\frac{1}{\|W\|+1}$,$\frac{1}{\|\Sigma_0\|+1}$, $\frac{1}{\gamma_A+1}$, $\frac{1}{\gamma_B+1}$, $\frac{1}{\gamma_R+1}$, $\frac{1}{\sx+1}$, $\sx$,  $\frac{1}{\srmin+1}$, $\srmin$, $\frac{1}{\sqmin+1}$, $\sqmin$, and $\frac{1}{\vertiii{\pmb{K}^*}+1}$, the step size condition \eqref{eqn:step_size} is satisfied.
Hence, there exists  $\eta_0\in\HH\left(\frac{1}{\sum_{i=1}^NC^i(\pmb{K}^{i,(0)},\pmb{K}^{-i*})+1}\right)$ as an appropriate polynomial in $\frac{1}{\sum_{i=1}^NC^i(\pmb{K}^{i,(0)},\pmb{K}^{-i*})+1}$, $\frac{1}{\sum_i k_i+1}$, $\frac{1}{d+1}$, $d$, $\frac{1}{N+1}$, $\frac{1}{T+1}$, $\frac{1}{\|W\|+1}$,$\frac{1}{\|\Sigma_0\|+1}$, $\frac{1}{\gamma_A+1}$, $\frac{1}{\gamma_B+1}$, $\frac{1}{\gamma_R+1}$, $\frac{1}{\sx+1}$, $\sx$,  $\frac{1}{\srmin+1}$, $\srmin$, $\frac{1}{\sqmin+1}$, $\sqmin$, and $\frac{1}{\vertiii{\pmb{K}^*}+1}$, such that when $\eta<\eta_0$, the step size condition \eqref{eqn:step_size} is satisfied. Therefore, by Lemma \ref{lemma:one_step_loc_conv}, we have
\begin{eqnarray*}
 \sum_{i=1}^N\Big( C^i(\pmb{K}^{i,(1)},\pmb{K}^{-i*})-C^{i*}) \Big) \leq (1-\widehat{\alpha}\eta) \sum_{i=1}^N\Big( C^i(\pmb{K}^{i,(0)},\pmb{K}^{-i*})-C^{i*}) \Big),
\end{eqnarray*}
with $\widehat{\alpha}$ defined in \eqref{eqn:defn_hat_alpha}, which implies the total cost of $N$ players decreases at $m=1$. Proceeding inductively, assume the following facts hold at round $m$:
\begin{itemize}
    \item $C^i(\pmb{K}^{i,(m-1)},\pmb{K}^{-i*})-C^{i*}\leq \psi$ for $i=1,2,\cdots,N$;
    \item $\rho_{\pmb{K}^{(m-1)}}\leq \bar{\rho}$;
    \item $\sum_{i=1}^N\Big( C^i(\pmb{K}^{i,(m)},\pmb{K}^{-i*})-C^{i*} \Big) \leq (1-\widehat{\alpha}\eta )\sum_{i=1}^N\Big( C^i(\pmb{K}^{i,(m-1)},\pmb{K}^{-i*})-C^{i*}\Big)$.
\end{itemize}
Now we prove the above facts also hold in round $m+1$.
% Note that the step size condition \eqref{eqn:step_size} is still satisfied by Lemma \ref{lemma:grad_bound} since $\sum_{i=1}^N C^i(\pmb{K}^{i,(m-1)},\pmb{K}^{-i*}) \leq\sum_{i=1}^NC^i(\pmb{K}^{i,(0)},\pmb{K}^{-i*})$.} 
Taking $\pmb{K} = \pmb{K}^{(m-1)}$ and  $\pmb{K}^{\prime} =  \pmb{K}^{(m)}$ in \eqref{eq:C_diff1}, we have 
\begin{eqnarray*}
C^i(\pmb{K}^{i,(m)},\pmb{K}^{-i*})-C^{i*} &\leq& \big(1-\eta\frac{\sx\srmin}{\|\Sigma_{\pmb{K}^{*}}\|}\big)\big(C^i(\pmb{K}^{i,(m-1)},\pmb{K}^{-i*})-C^{i*}\big) \\
&&+ \,\eta\,\, h_{\rm glob}\frac{T(N-1)}{\sx \srmin}\left(\sum_{j=1,j\neq i}^N (C^j(\pmb{K}^{j,(m-1)},\pmb{K}^{-j*})-C^{j*}) \right)\\
&\leq & \big(1-\eta\frac{\sx\srmin}{\|\Sigma_{\pmb{K}^{*}}\|}\big) \psi + \eta\,\, h_{\rm glob}\frac{T(N-1)}{\sx \srmin}(N-1)\psi \leq \psi.
\end{eqnarray*}
The last inequality holds since $1-\eta\frac{\sx\srmin}{\|\Sigma_{\pmb{K}^{*}}\|}+ \eta\,\, h_{\rm glob}\frac{T(N-1)^2}{\sx \srmin}<1$ under Assumption \ref{ass:noise_cond}. Therefore by \eqref{defn_rho_K}, we have
\[
\rho_{\pmb{K}^{(m)}} = \rho^* +N\gamma_B\sqrt{\frac{T}{\sx\srmin} \max_i\big\{C^i(\pmb{K}^{i,(m)},\pmb{K}^{-i*})- C^{i*}\big\}}+\frac{1}{20 T^2} \leq \bar{\rho}.
\]
Thus the step size condition \eqref{eqn:step_size} is still satisfied for $\pmb{K} = \pmb{K}^{(m)}$ and $\pmb{K}^{\prime} = \pmb{K}^{(m+1)}$ with $\eta<\eta_0$,
%$\eta\in\HH\left(\frac{1}{\sum_{i=1}^NC^i(\pmb{K}^{i,(0)},\pmb{K}^{-i*})+1}\right)$, 
since $\rho_{\pmb{K}^{(m)}}\leq \bar{\rho}$
and $\sum_{i=1}^N C^i(\pmb{K}^{i,(m)},\pmb{K}^{-i*}) \leq\sum_{i=1}^NC^i(\pmb{K}^{i,(0)},\pmb{K}^{-i*})$. Therefore, Lemma \ref{lemma:one_step_loc_conv} can be applied again for the update at round $m+1$ {with $\pmb{K}=\pmb{K}^{(m)}$ and $\pmb{K}^{\prime}=\pmb{K}^{(m+1)}$} to obtain:
\begin{eqnarray*}
 \sum_{i=1}^N\Big( C^i(\pmb{K}^{i,(m+1)},\pmb{K}^{-i*})-C^{i*} \Big) \leq
 (1-\widehat{\alpha}\eta )\sum_{i=1}^N\Big( C^i(\pmb{K}^{i,(m)},\pmb{K}^{-i*})-C^{i*} \Big).
\end{eqnarray*}
For $\epsilon>0$, provided $M\geq \frac{1}{\widehat{\alpha}\eta}\log\left(\frac{\sum_{i=1}^N(C^i(\pmb{K}^{i,(0)},\pmb{K}^{-i*})-C^{i*})}{\epsilon}\right)$, we have
\[
\sum_{i=1}^N\Big( C^i(\pmb{K}^{i,(M)},\pmb{K}^{-i*})-C^{i*}\Big)\leq\epsilon.
\]
\end{proof}

\section{The Natural Policy Gradient Method with Unknown Parameters}\label{sec:mf_NPG}

Based on the update rule in \eqref{eqn:nonzero_grad_update_rule}, it is straightforward to develop a 
{\it model-free version} of the natural policy gradient algorithm using {\it sampled data}. 
See Algorithm \ref{alg:NPG_unknwon_est} for the natural policy gradient method with unknown parameters.
In contrast to the model-based case, where the gradient $\nabla C^i(\pmb{K})$ and covariance matrix $\Sigma_t^i$ can be calculated directly, these two terms can not be calculated in the model-free setting since the model parameters are unknown.  Therefore, we propose to use a zeroth-order optimization method to estimate the gradient and an empirical covariance matrix to estimate the covariance matrix (see \eqref{eqn:algo_est}). Building upon the theories in Section \ref{sec:exact_NPG}, high-probability convergence guarantees (linear convergence rate and polynomial sample complexity) for the model-free counterpart can be established in the same way as for the 
Linear Quadratic Regulator setting in \cite{hambly2020policy}.

\iffalse
{\color{blue}[TODO:]
\begin{itemize}
    \item Add a few comments on Algorithm \ref{alg:NPG_unknwon_est}.
    \item Explain the difference between Algorithm \ref{alg:NPG_unknwon_est} and  Algorithm \ref{alg:NPG_known}.
    \item Add some intuition on how to show the convergence of Algorithm \ref{alg:NPG_unknwon_est} using method in \cite{hambly2020policy}.
\end{itemize}
}
\fi

\begin{algorithm}[H]
\caption{\textbf{Natural Policy Gradient Method with Unknown Parameters}}
\label{alg:NPG_unknwon_est}
\begin{algorithmic}[1]
    \STATE \textbf{Input}: Number of iterations $M$, time horizon $T$, initial policy $\pmb{K}^{(0)}=(\pmb{K}^{1,(0)},\cdots,\pmb{K}^{N,(0)})$, step size $\eta$, number of trajectories $L$, smoothing parameters $r_i$, dimensions $D_i=k_i\times d$.
    \FOR {$m\in\{1, \ldots, M\}$}
        \FOR {$i\in\{1,\ldots,N\}$}
        \FOR {$l\in\{1, \ldots, L\}$}
        \FOR {$t\in\{0, \ldots, T-1\}$}
           \STATE Sample the (sub)-policy at time $t$: $\widehat{K}^{i,l}_t = K_t^{i,(m-1)} +U^{i,l}_t$ where $U_t^{i,l}$ is drawn uniformly at random over matrices such that $\|U_t^{i,l}\|_F=r_i$.
           \STATE Denote  $\widehat{c_t}^{i,l}$ as the single trajectory cost of player $i$ with policy $(\widehat{\pmb{K}}_{l,t}^{i,(m-1)},\pmb{K}^{-i,(m-1)})$ where $\widehat{\pmb{K}}_{l,t}^{i,(m-1)}:=({K}_0^{i,(m-1)},\cdots,{K}_{t-1}^{i,(m-1)},\widehat{K}^{i,l}_t,{K}_t^{i,(m-1)},\cdots,K_{T-1}^{i,(m-1)})$ starting from $x^l_0$.
           \STATE Denote $\widehat{\Sigma}_t^{i,l}$ as the state covariance matrix with $\widehat{\Sigma}_t^{i,l}=x_t^{i,l}(x_t^{i,l})^\top$.
            \ENDFOR
            \ENDFOR
            \ENDFOR
             \STATE Obtain the estimates of $\nabla_{K_t^i} C^i(\pmb{K}^{(m-1)})$ and $\Sigma_t^{\pmb{K}^{(m-1)}}$ for each $i$ and $t$:
 \begin{eqnarray}\label{eqn:algo_est}
\widehat{\nabla_{K_t^i} C^i(\pmb{K}^{(m-1)})} = \frac{1}{L}\sum_{l=1}^L \frac{D_i}{r_i^2}\,\widehat{c_t}^{i,l}\, U^{i,l}_t,\qquad \widehat{\Sigma}_t^i=\frac{1}{L}\sum_{l=1}^L\widehat{\Sigma}_t^{i,l}.
 \end{eqnarray}
 \STATE Update the policies using natural policy gradient updating rule: 
            \begin{equation}
                 K_t^{i,(m)} = K_t^{i,(m-1)} - \eta \widehat{\nabla_{K_t^i} C^i(\pmb{K}^{(m-1)})}(\widehat{\Sigma}_t^i)^{-1}.
             \end{equation}
        \ENDFOR
 \STATE Return the iterates $\pmb{K}^{(M)}=(\pmb{K}^{1,(M)},\cdots,\pmb{K}^{N,(M)})$.
\end{algorithmic}
\end{algorithm}

% \paragraph{Experimental Results for the Bargaining Model.} We perform the NPG algorithm with unknown parameters to the bargaining model introduced in Section \ref{subsec:bargaining} under the parameters given in Figure \ref{fig:optimal_state_traj}. The state trajectory under the learned policies from the algorithm (see Figure \ref{fig:state_mf_learned}) is very similar to the one under Nash equilibrium shown in Figure \ref{fig:optimal_state_traj}, and the normalized error falls below 0.05 within 10000 iterations (see Figure \ref{fig:conv_mf_bargaining}). Note that although Assumption \ref{ass:nonzero_cost} and \ref{ass:nonzero_initial_noise} are not satisfied under this set of parameters, the NPG will converge to the Nash equilibrium under proper conditions.
%     \begin{figure}[H]
%   \centering
%   \begin{subfigure}[b]{0.43\textwidth}
%     \includegraphics[width=\textwidth]{images/MF_learned_K.png}
%     \caption{\label{fig:state_mf_learned}State trajectory under learned policies.}
%   \end{subfigure}
%   \begin{subfigure}[b]{0.43\textwidth}
%     \includegraphics[width=\textwidth]{images/Model_free_bar_sig_5_eta_10_5_xlog.png}
%     \caption{\label{fig:conv_mf_bargaining}Convergence.}
%   \end{subfigure}
%   \caption{Performance of the NPG algorithm with unknown parameters for the bargaining model (Model parameters are the same as in Figure \ref{fig:optimal_state_traj}; for the NPG algorithm with unknown parameters, number of trajectory $L_1=L_2=500$, and smoothing parameters $r_1=r_2=0.5$.).}
% \end{figure}

\section{Numerical Experiments}\label{sec:numerical_experiments}

We demonstrate the performance of the natural policy gradient algorithms with three general-sum game examples.  We will specifically focus on the following questions:
\begin{itemize}
    \item In practice, how fast do the natural policy gradient algorithms %with known and unknown parameters 
    converge to the true solution?
   How sensitive is the natural policy gradient algorithm %with known parameters 
    to the step size and the initial policy?
    \item Do natural gradient methods converge when Assumption \ref{ass:noise_cond} is violated? How restrictive is Assumption \ref{ass:noise_cond} in practice?
   \item Can Theorem \ref{thm:local_conv_NPG} provide any guidance on hyper-parameter tuning? In particular, does adding system noise improve the convergence?
    %\item How does the system noise affect the performance of the natural policy gradient algorithms?
\end{itemize}

The first example is a modified example from \cite{mazumdar2019policy}, in which they show that in the setting of infinite time horizon and deterministic dynamics, the (vanilla) policy gradient algorithms have no guarantees of even local convergence to the Nash equilibria with known parameters. We will show in Section \ref{sec:compar_infi_paper} that, under the same experimental set-up but over a finite time horizon {with stochastic dynamics}, the natural policy gradient algorithm with known parameters  finds the Nash equilibrium with properly chosen initial policies and step sizes. The second example is a two-player LQ game with synthetic data (see Section \ref{sec:sys_noise}). We will show that the system noise helps the natural policy gradient algorithm with unknown parameters to converge to the Nash equilibrium. Finally, we investigate  the algorithm's performance with known and unknown parameters for a three-player game example in Section \ref{sec:three_player}.

%This section is organized as follows. We demonstrate the performance of the natural policy gradient algorithm with known parameters under the setting of a modified example of \cite{mazumdar2019policy} in Section \ref{sec:compar_infi_paper}. We then perform the natural policy gradient algorithm with unknown parameters to a two-player LQG example with synthetic data to show the impact of the system noise in Section \ref{sec:sys_noise}.

\paragraph{Performance Measure.} We use the following \textit{normalized error} to quantify the performance of a given pair of policies $(\pmb{K}^1,\pmb{K}^2)$: for $i,j=1,2$ and $j\neq i$,
\[
\text{Normalized error (of player} \,\,i) =  \frac{C^i(\pmb{K}^i,\pmb{K}^{j*})-C^i(\pmb{K}^{1*},\pmb{K}^{2*})}{C^i(\pmb{K}^{1*},\pmb{K}^{2*})}.
\]

\subsection{Convergence of the Natural Policy Gradient Algorithm}\label{sec:compar_infi_paper}

For policy gradient MARL under the setting of N-player general-sum LQ games, some difficulties in convergence have been identified in some empirical studies. For example, Mazumdar et al. in \cite{mazumdar2019policy} show by a counterexample that in the setting of infinite time horizon and deterministic dynamics, the (vanilla) policy gradient method avoids the Nash equilibria for a non-negligible subset of problems {(with known parameters)}. In this section, we illustrate that in our setting with finite time horizon and stochastic dynamics, the natural policy gradient algorithm (with known parameters) finds the (unique) Nash equilibrium under the same experimental set-up as in a modified example given in \cite{mazumdar2019policy}.

\paragraph{Set-up.} We set up the model parameters and initialize the policies in the same way as Section 5.1 of \cite{mazumdar2019policy}. Under the following set of parameters, there exists a unique Nash equilibrium since the sufficient condition in Remark \ref{remark:suff_unique} is satisfied.
\begin{enumerate}
    \item Parameters: for $t=1,\cdots,T-1$
    \[
A_t = 
\begin{bmatrix}
0.588 & 0.028\\
0.570 & 0.056
\end{bmatrix},\quad
B_t^1 = 
\begin{bmatrix}
1 \\
1
\end{bmatrix},\quad
B_t^2 = 
\begin{bmatrix}
0 \\
1
\end{bmatrix},\quad
W = 
\begin{bmatrix}
\sigma^2 & 0 \\
0 & \sigma^2
\end{bmatrix},
\]
    \[
Q_T^1=Q_t^1 = 
\begin{bmatrix}
0.01 & 0\\
0 & 1
\end{bmatrix},\quad
Q_T^2 =Q_t^2 = 
\begin{bmatrix}
1 & 0\\
0 & 0.147
\end{bmatrix},\quad
R_t^1(t)=R_t^2(t)=0.01,
\]
where $\sigma\in\mathbb{R}$ and $T=10$.
    \item Initialization: we assume the initial state distribution to be $[1,1]^\top$ or $[1,1.1]^\top$ with probability 0.5 each. We initialize both players' policies $K_t^{i,(0)}=(K_{t0}^{i,(0)},K_{t1}^{i,(0)})$ such that $(K_{t0}^{i,(0)}-K_{t0}^{i*})^2+(K_{t1}^{i,(0)}-K_{t1}^{i*})^2\leq r^2$, where  $K_t^{i*}=(K_{t0}^{i*},K_{t1}^{i*})$ denotes the Nash equilibrium, and $r$ is the radius of the ball centered at $K_t^{i*}$ in which we initialize the policies.
\end{enumerate}

\begin{figure}[htbp]
  \centering
  \begin{subfigure}[b]{0.43\textwidth}
    \includegraphics[width=\textwidth]{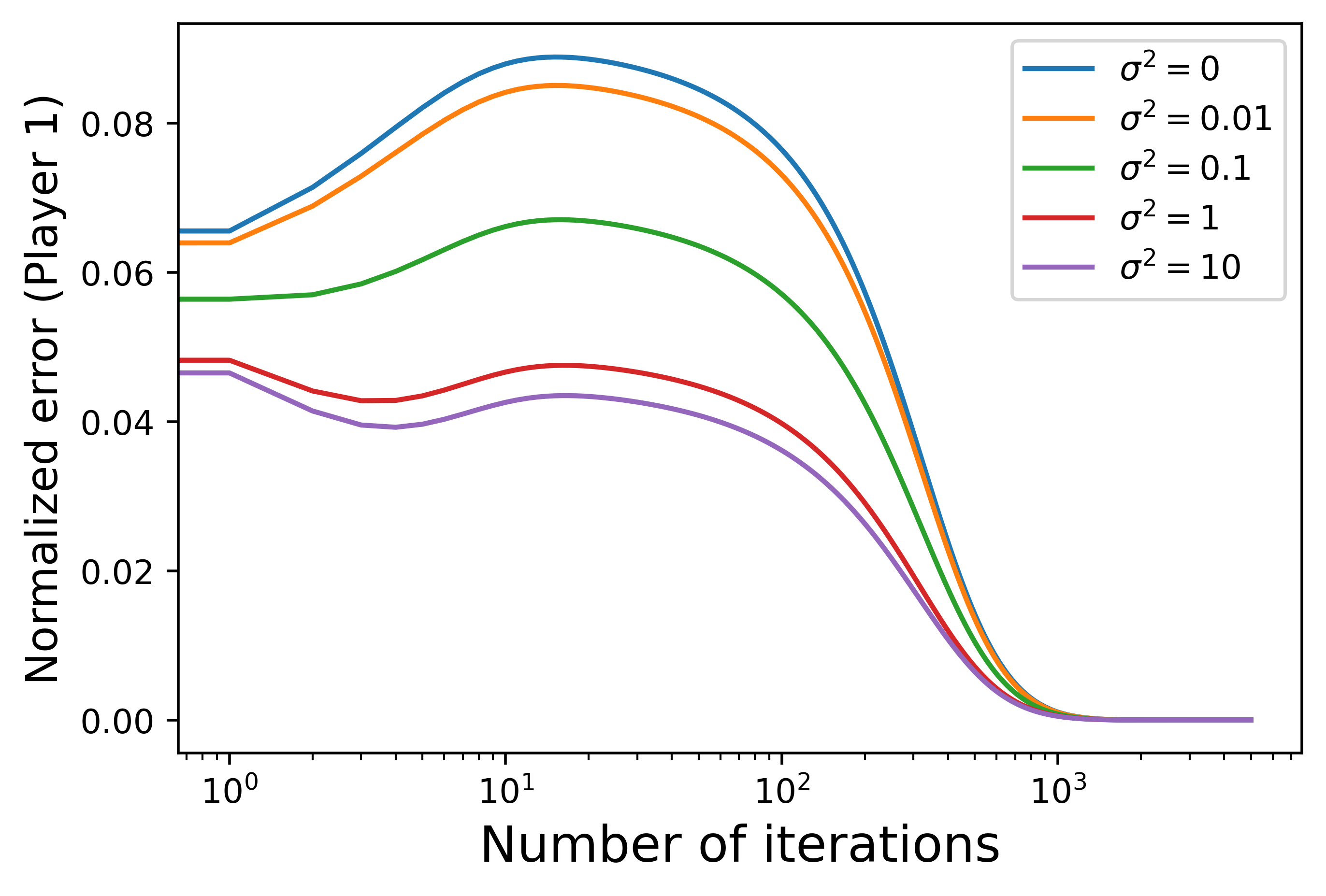}
    \caption{Player 1.}
  \end{subfigure}
  \begin{subfigure}[b]{0.43\textwidth}
    \includegraphics[width=\textwidth]{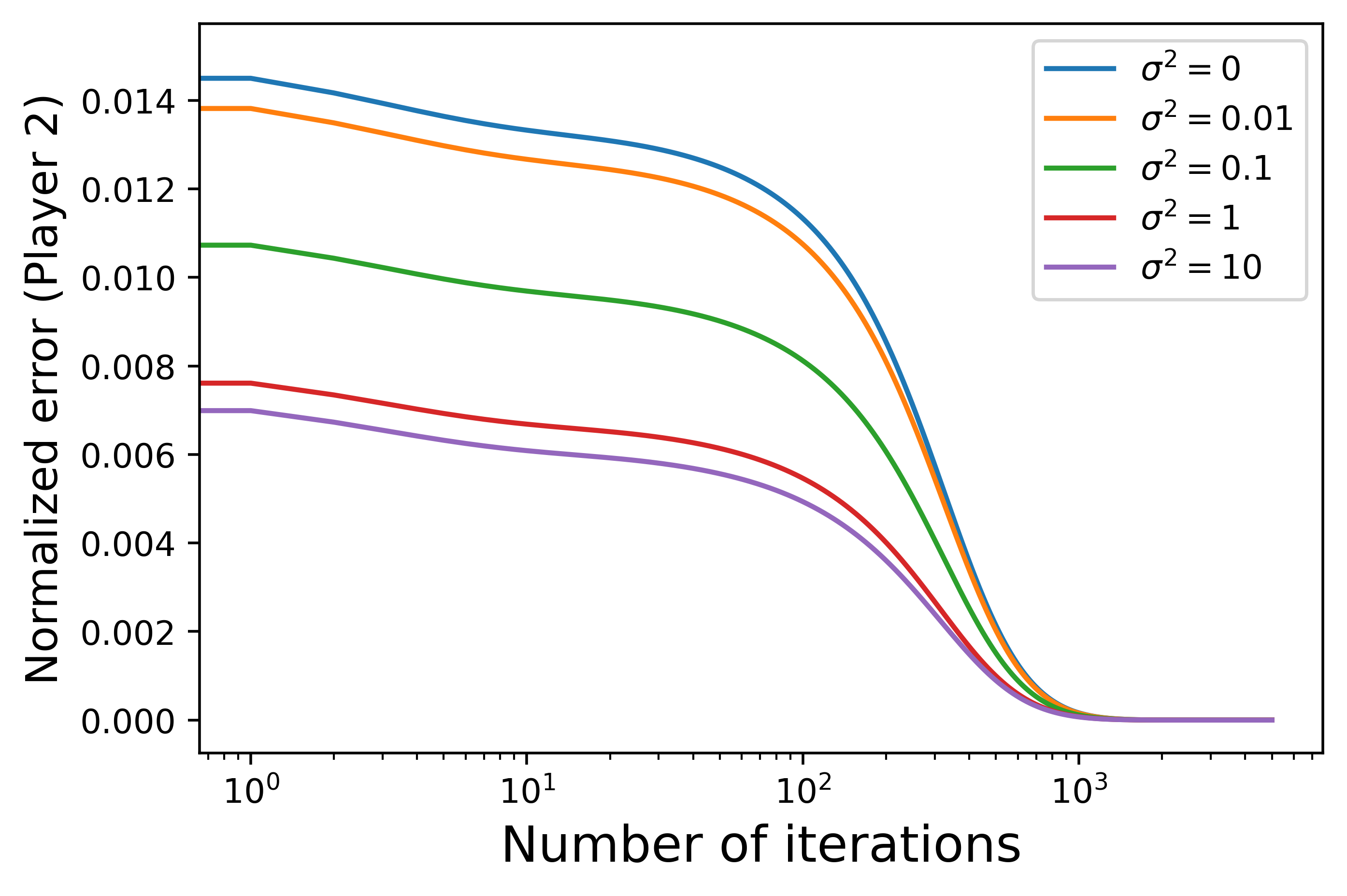}
    \caption{Player 2.}
  \end{subfigure}
  \caption{\label{fig:NPG_noise_small_r}Normalized error under different $\sigma^2$ when $r=0.25$ ($\eta_1=\eta_2=0.1$ and $M=5000$).}
\end{figure}

\begin{figure}[htbp]
  \centering
  \begin{subfigure}[b]{0.43\textwidth}
    \includegraphics[width=\textwidth]{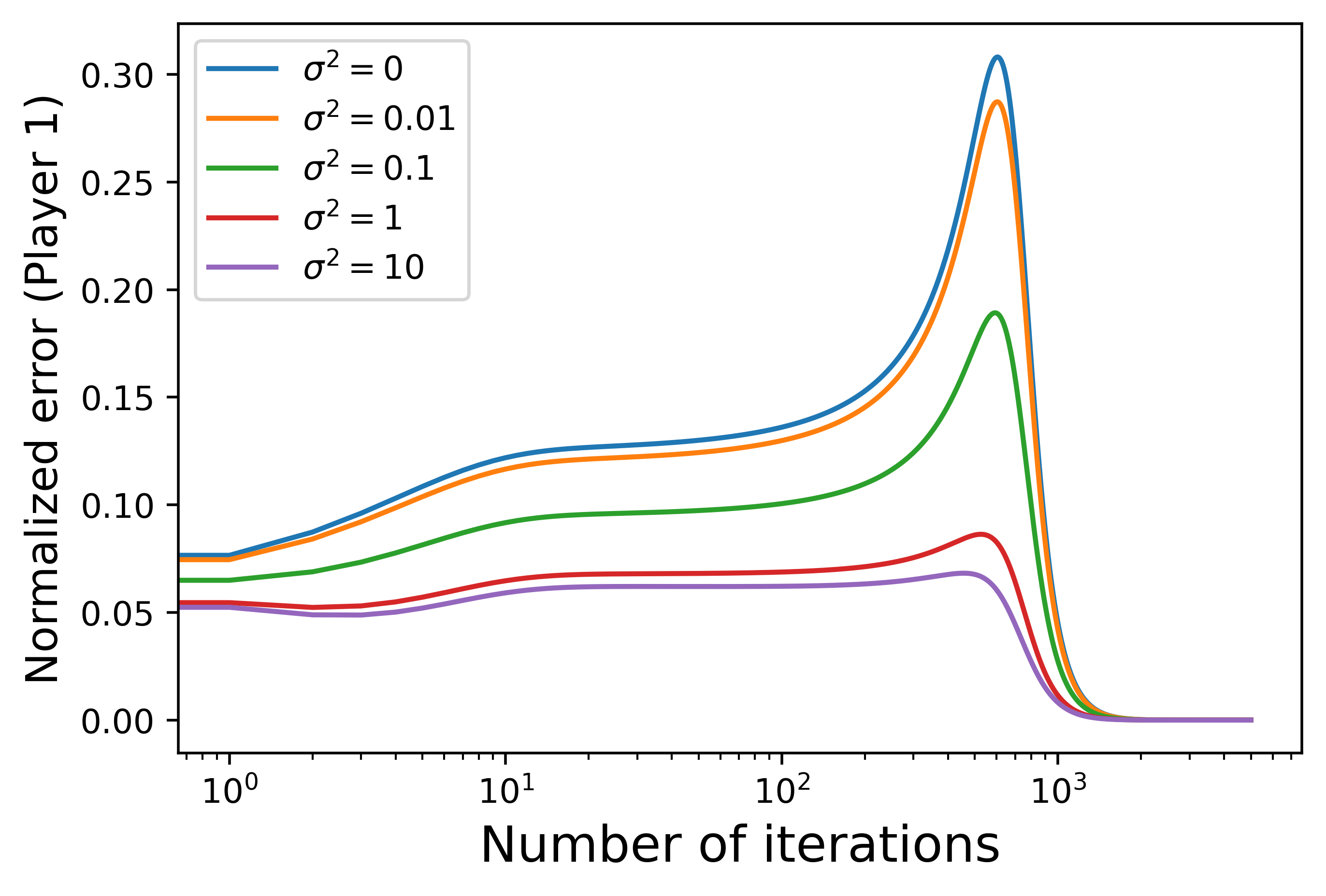}
    \caption{Player 1.}
  \end{subfigure}
  \begin{subfigure}[b]{0.43\textwidth}
    \includegraphics[width=\textwidth]{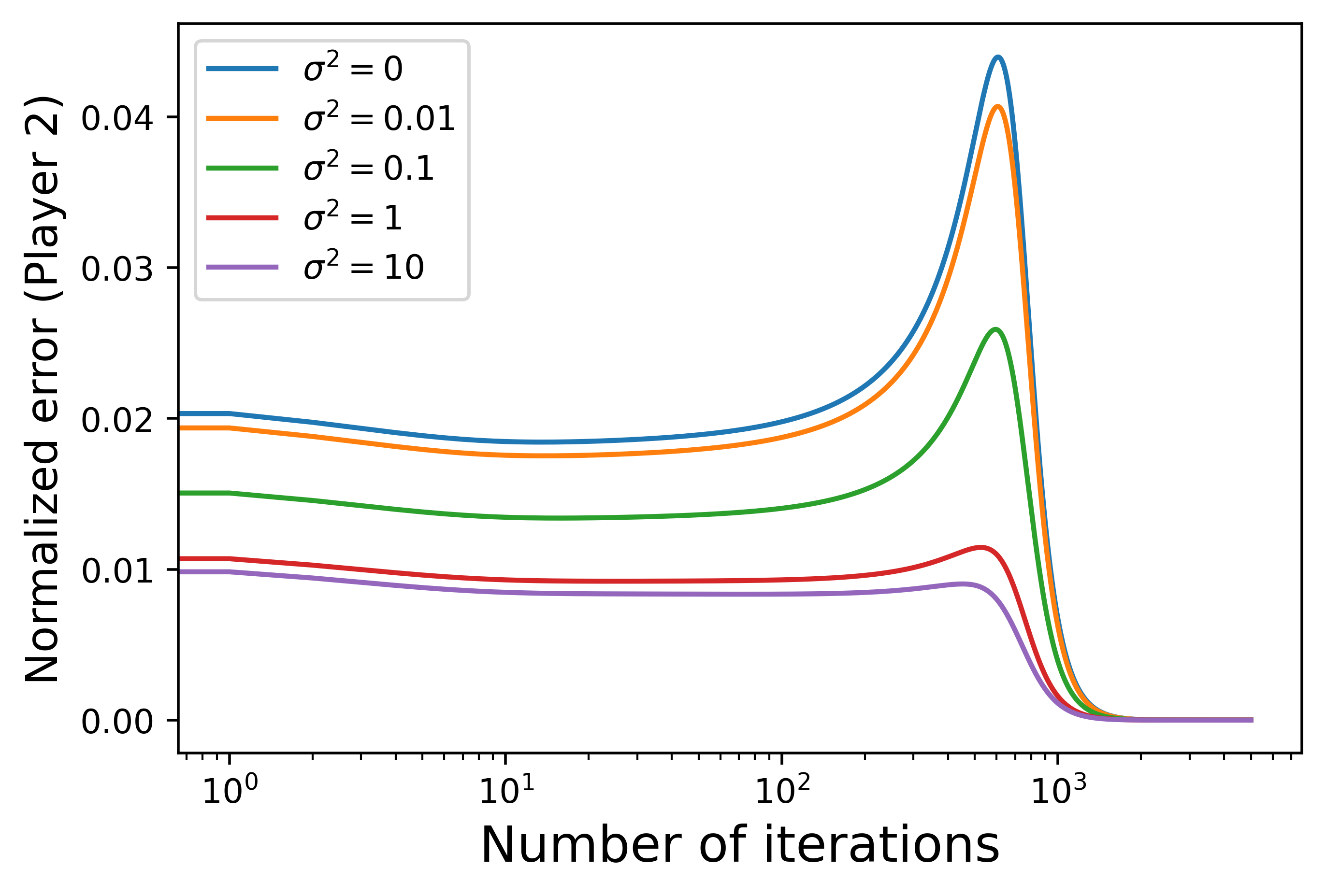}
    \caption{Player 2.}
  \end{subfigure}
  \caption{\label{fig:NPG_noise_large_r}Normalized error under different $\sigma^2$ when $r=0.3$ ($\eta_1=\eta_2=0.1$ and $M=5000$).}
\end{figure}

\paragraph{Convergence.} The natural policy gradient algorithm shows a reasonable level of accuracy within 1000 iterations (i.e., the normalized error is less than 0.5\%) for both players under different levels of system noise $\sigma^2$, which ranges from 0 (deterministic dynamics) to 10. See Figure \ref{fig:NPG_noise_small_r} for the case where $r=0.25$ and Figure \ref{fig:NPG_noise_large_r} for the case where $r=0.30$. We observe that when we initialize the policies in a larger neighborhood of the Nash equilibrium, it takes the natural policy gradient algorithm (with the same step size) more iterations to converge. {There is a sharp peak in the normalized error for player 2 and this peak diminishes when the noise level increases.}
%{\color{red}[add a summary of the results on different values of $\sigma$.]}

\begin{figure}[htbp]
  \centering
  \begin{subfigure}[b]{0.3\textwidth}
    \includegraphics[width=\textwidth]{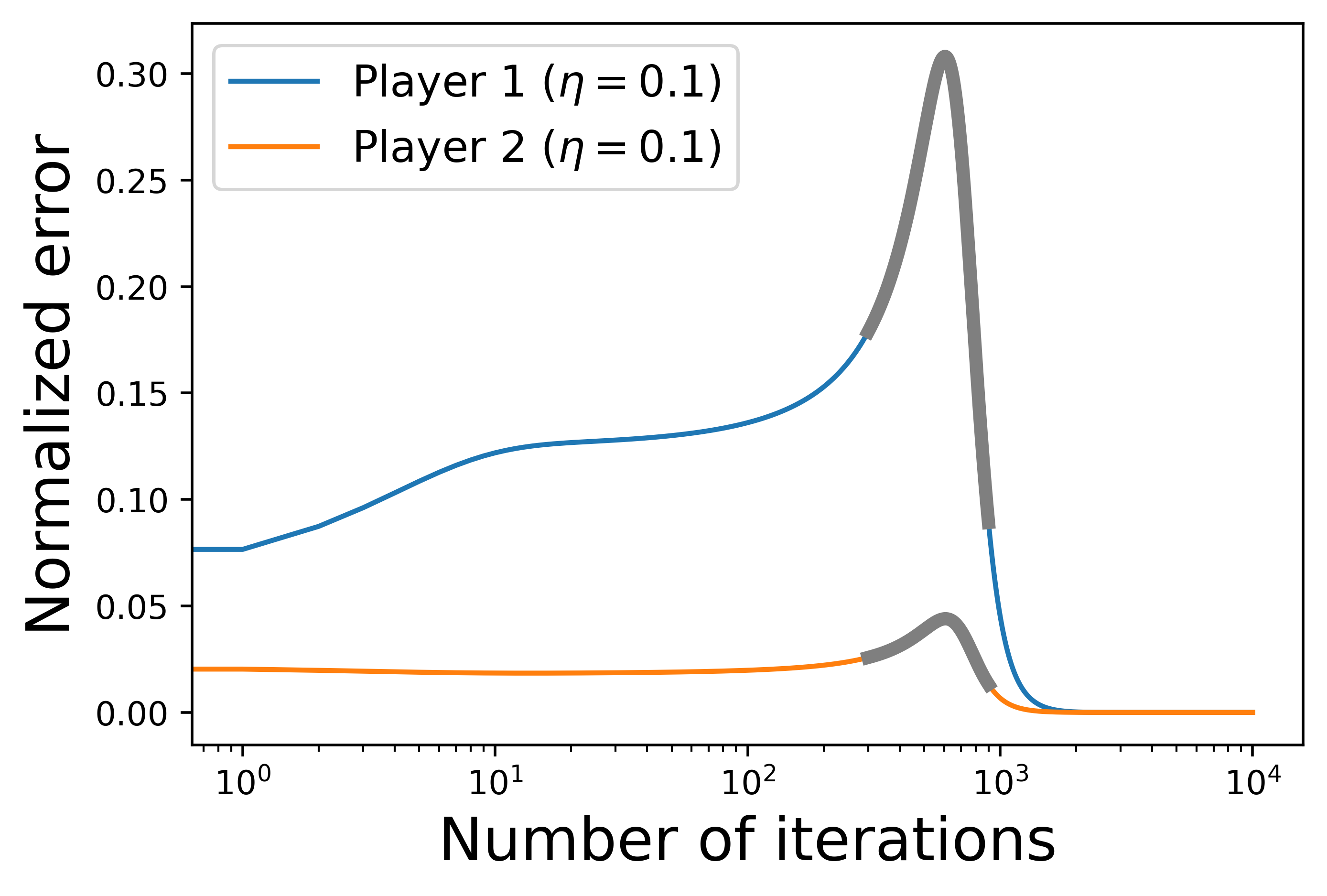}
    \caption{\label{fig:peak_nor_error}Normalized error.}
  \end{subfigure}
  \begin{subfigure}[b]{0.3\textwidth}
    \includegraphics[width=\textwidth]{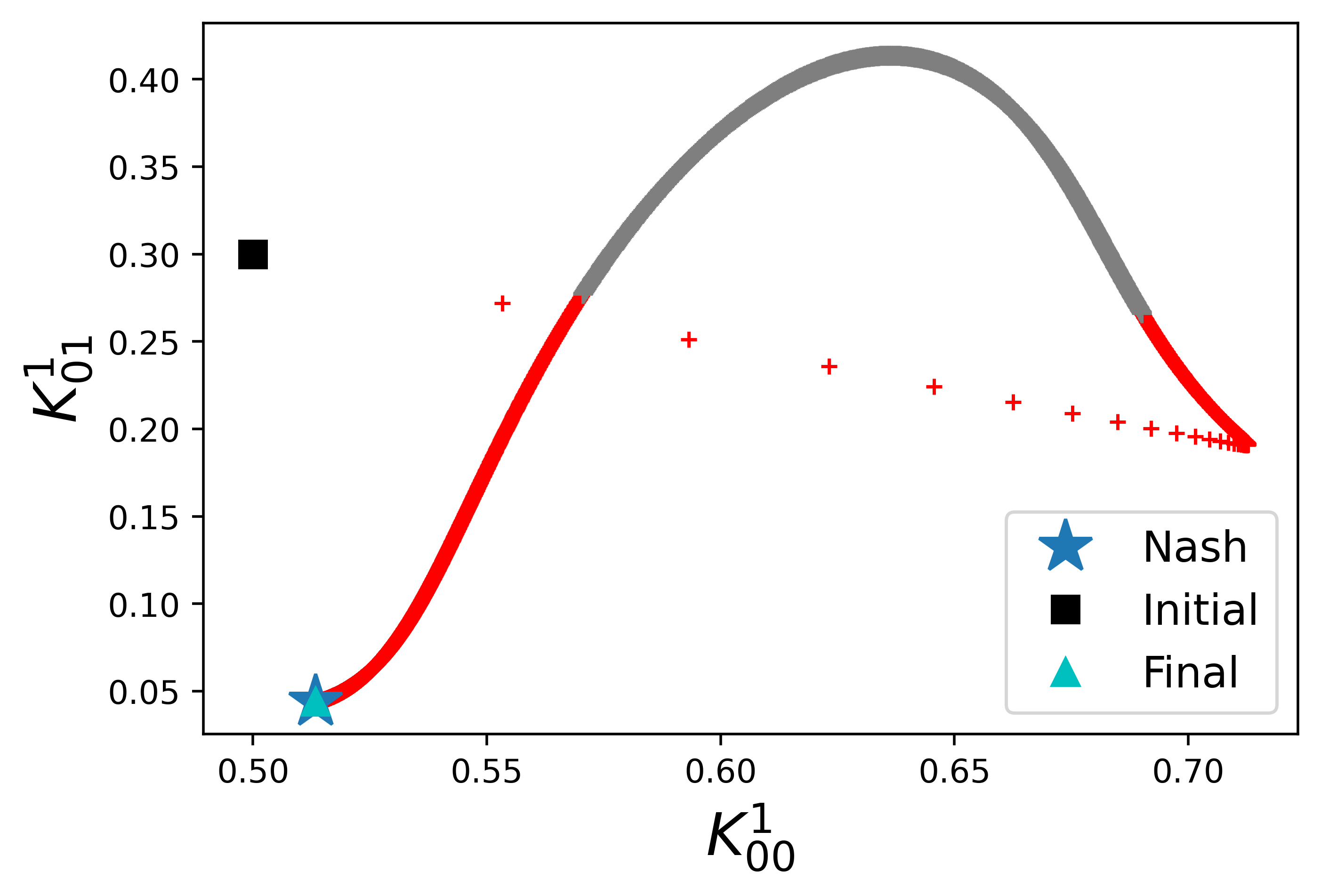}
    \caption{\label{fig:peak_K1_traj}Trajectory of $K_{0}^{1}$.}
  \end{subfigure}
   \begin{subfigure}[b]{0.3\textwidth}
    \includegraphics[width=\textwidth]{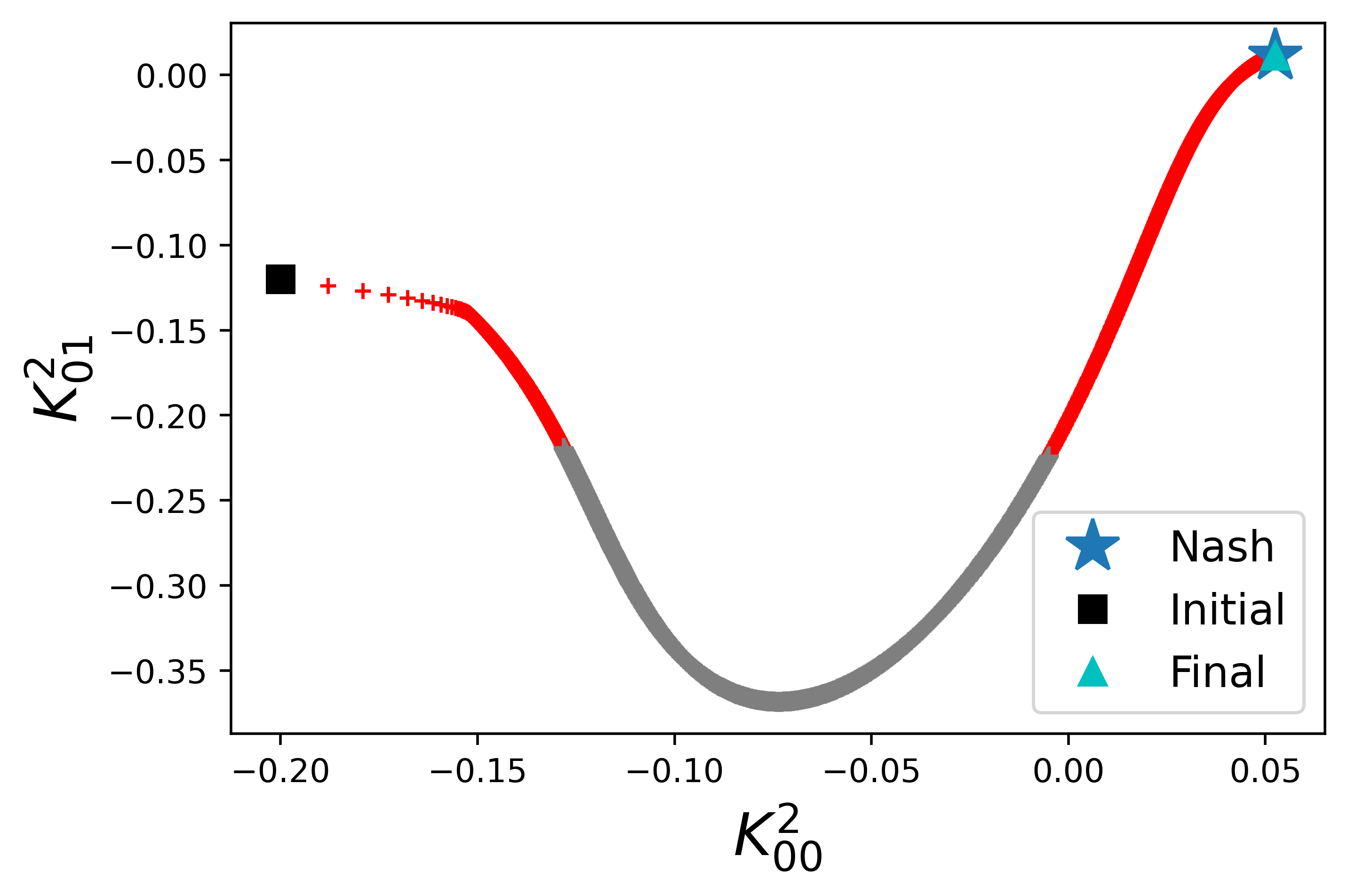}
    \caption{\label{fig:peak_K2_traj}Trajectory of $K_{0}^{2}$.}
  \end{subfigure}
  \caption{\label{fig:NPG_noise_peak}Normalized error and trajectories of learned policies with $\sigma^2=0$ and $r=0.3$ ($\eta_1=\eta_2=0.1$, $M=10000$). The peak period between iterations 300-900 is indicated in grey.}
\end{figure}

In Figure \ref{fig:NPG_noise_peak}, we show the normalized error and the corresponding trajectories of learned policies near the peak observed in Figure \ref{fig:NPG_noise_large_r} under $\sigma=0$. We denote by $K_0^i=(K_{00}^i,K_{01}^i)$ the learned policy at $t=0$ of player $i$. The peak period (iterations 300-900) is indicated in grey in Figures \ref{fig:peak_nor_error}-\ref{fig:peak_K2_traj}, and the trajectories of learned policy $K_0^i$ for the rest of the whole 10000 iterations are indicated in red in Figures \ref{fig:peak_K1_traj}-\ref{fig:peak_K2_traj}. The natural policy gradient algorithm {overshoots in the first few iterations but} detects the right direction after about 500 iterations and eventually converges to the Nash equilibrium (see the blue star in Figures \ref{fig:peak_K1_traj}-\ref{fig:peak_K2_traj}).

\begin{figure}[htbp]
  \centering
  \begin{subfigure}[b]{0.3\textwidth}
    \includegraphics[width=\textwidth]{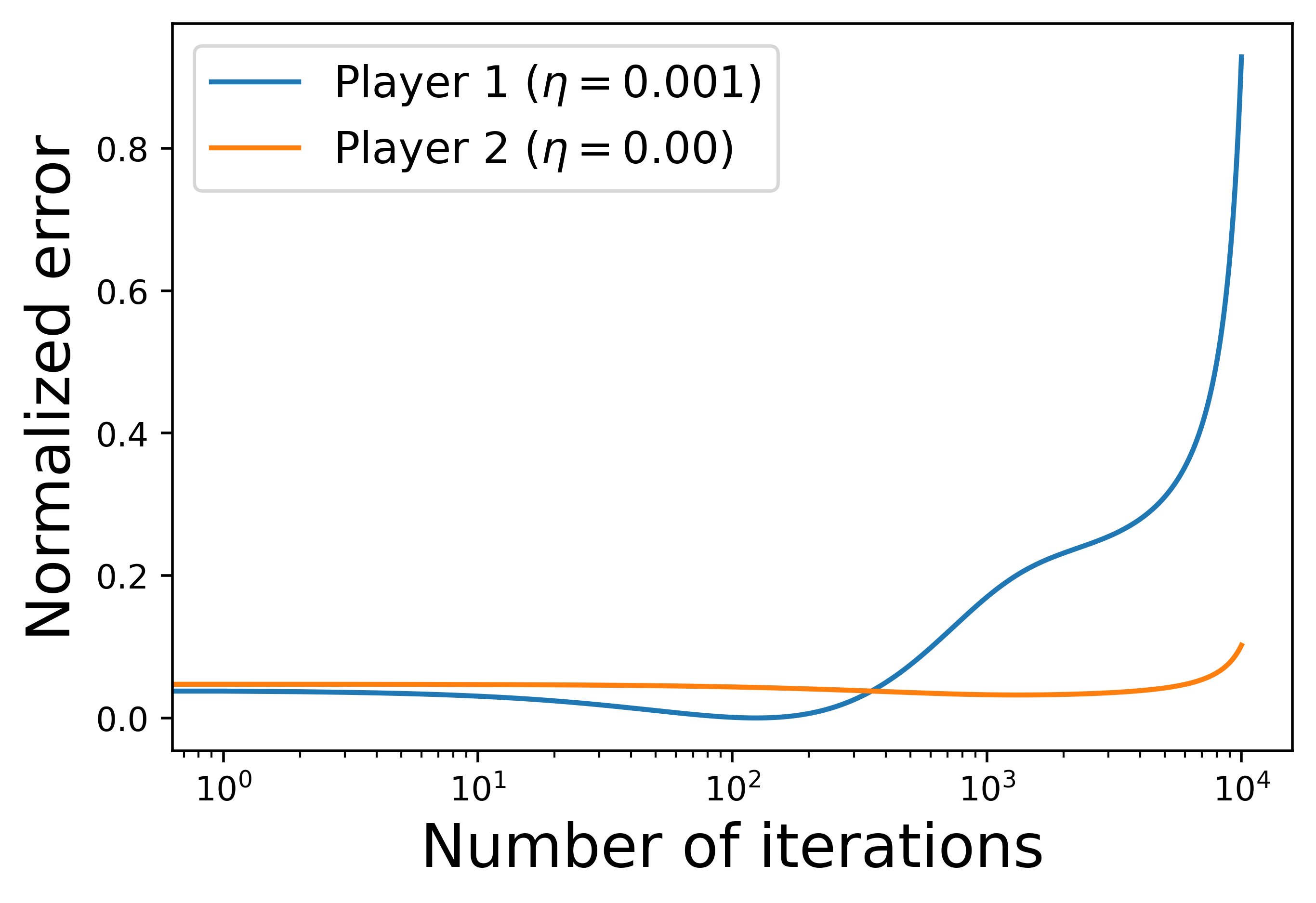}
    \caption{Normalized error.}
  \end{subfigure}
  \begin{subfigure}[b]{0.3\textwidth}
    \includegraphics[width=\textwidth]{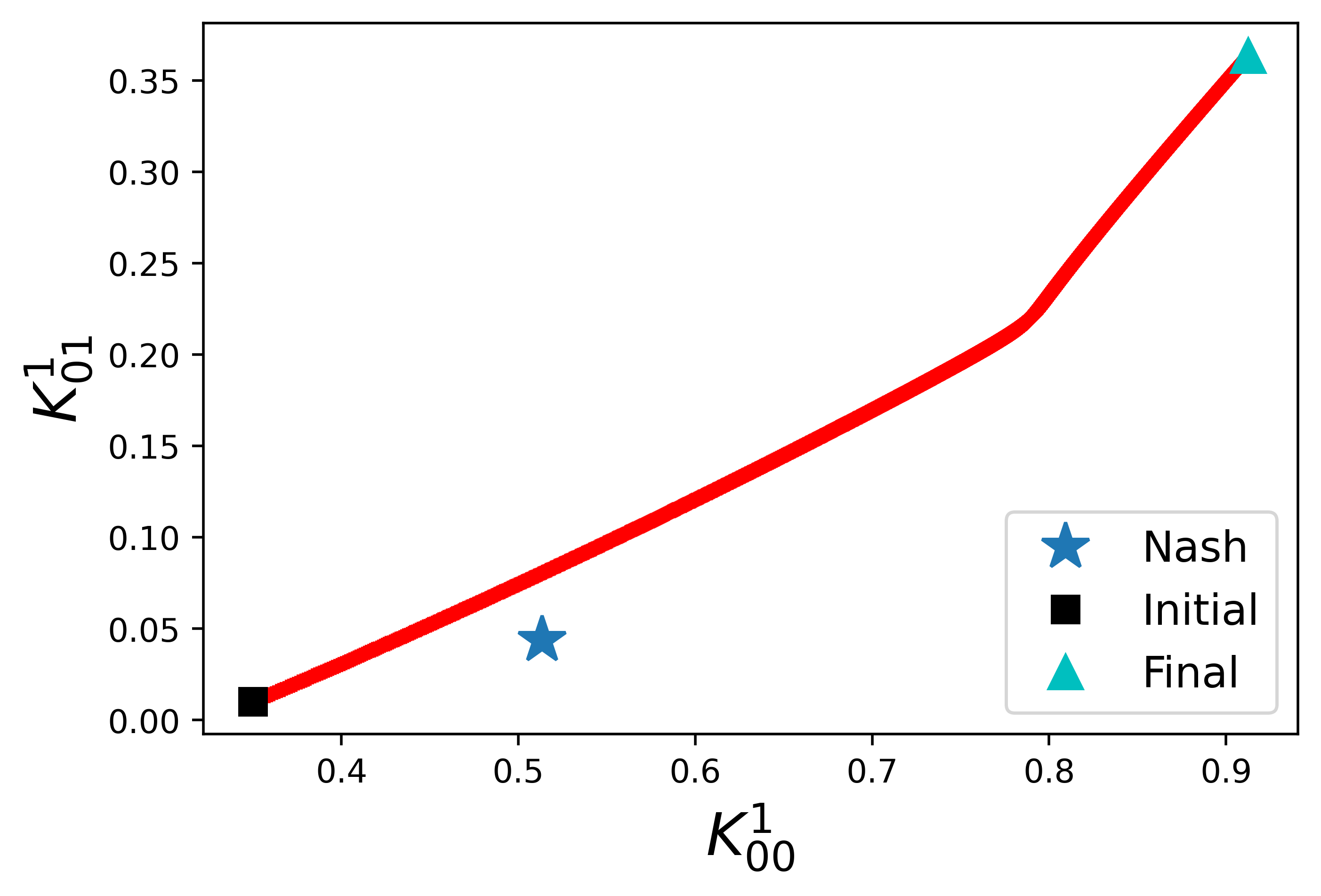}
    \caption{Trajectory of $K_{0}^{1}$.}
  \end{subfigure}
   \begin{subfigure}[b]{0.3\textwidth}
    \includegraphics[width=\textwidth]{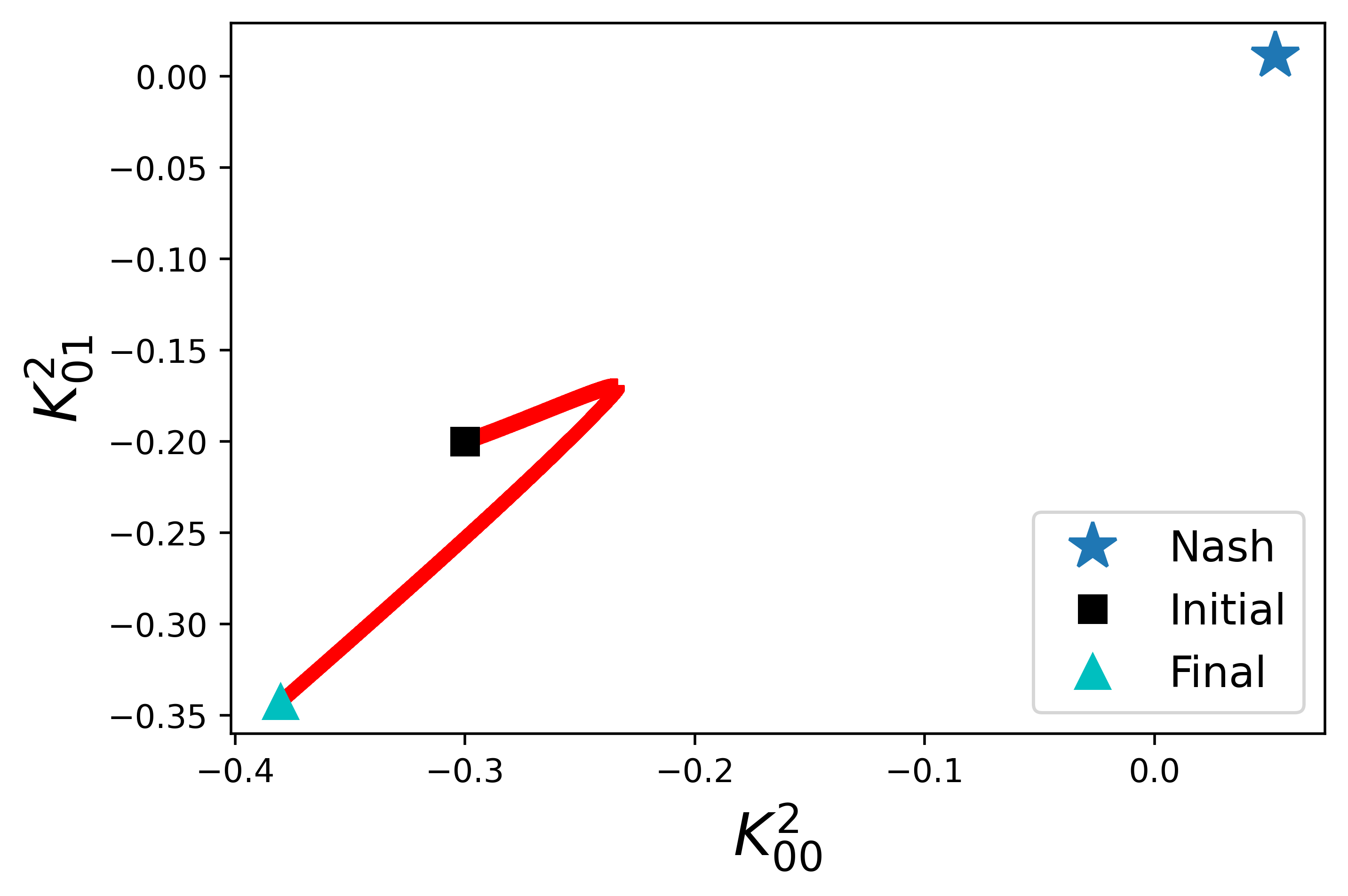}
    \caption{Trajectory of $K_{0}^{2}$.}
  \end{subfigure}
  \caption{\label{fig:NPG_div}Performance of the natural policy gradient algorithm with $r=0.42$ and $\eta_1=\eta_2=0.001$ ($M=10000$).}
\end{figure}

\begin{figure}[htbp]
  \centering
  \begin{subfigure}[b]{0.3\textwidth}
    \includegraphics[width=\textwidth]{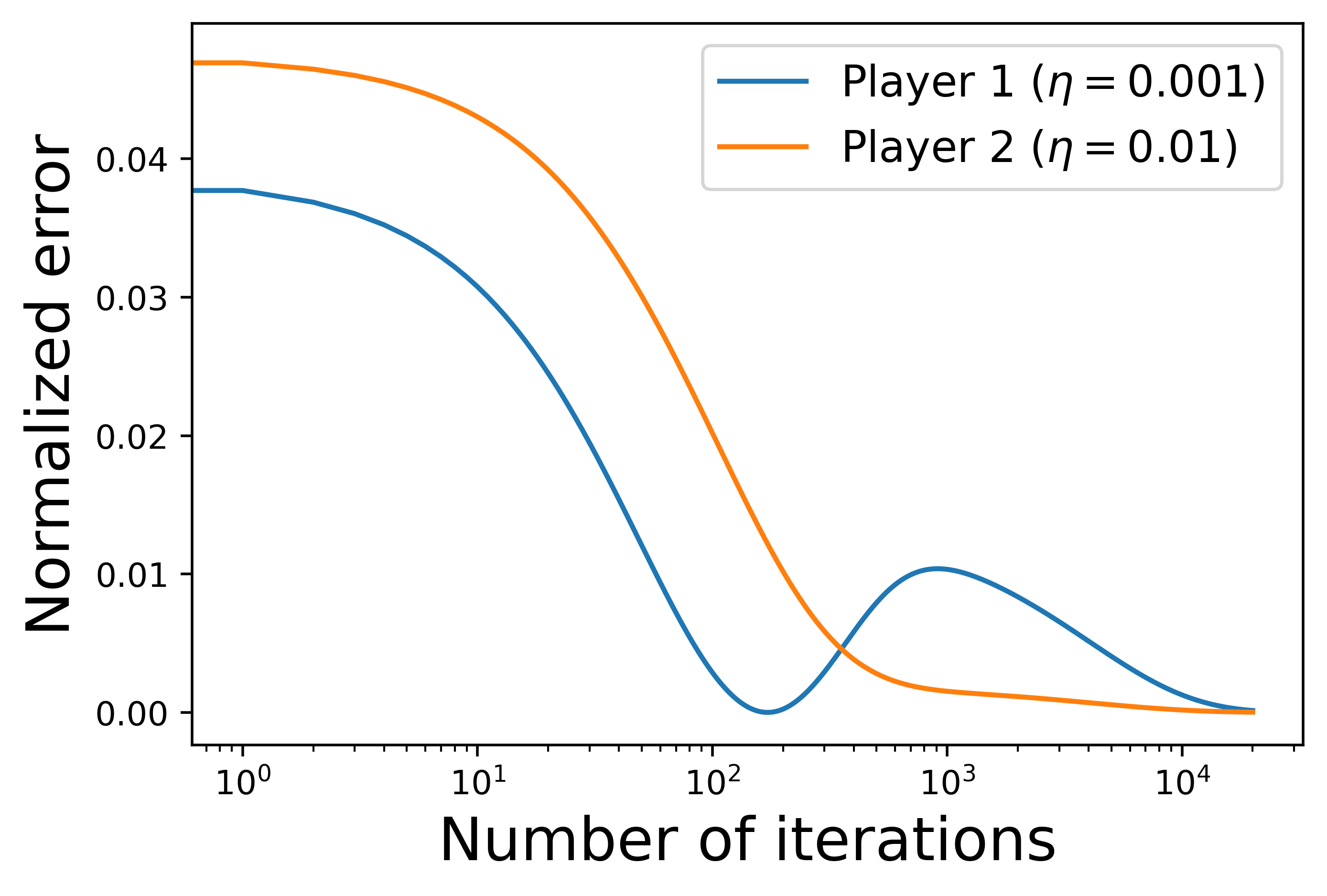}
    \caption{Normalized error.}
  \end{subfigure}
  \begin{subfigure}[b]{0.3\textwidth}
    \includegraphics[width=\textwidth]{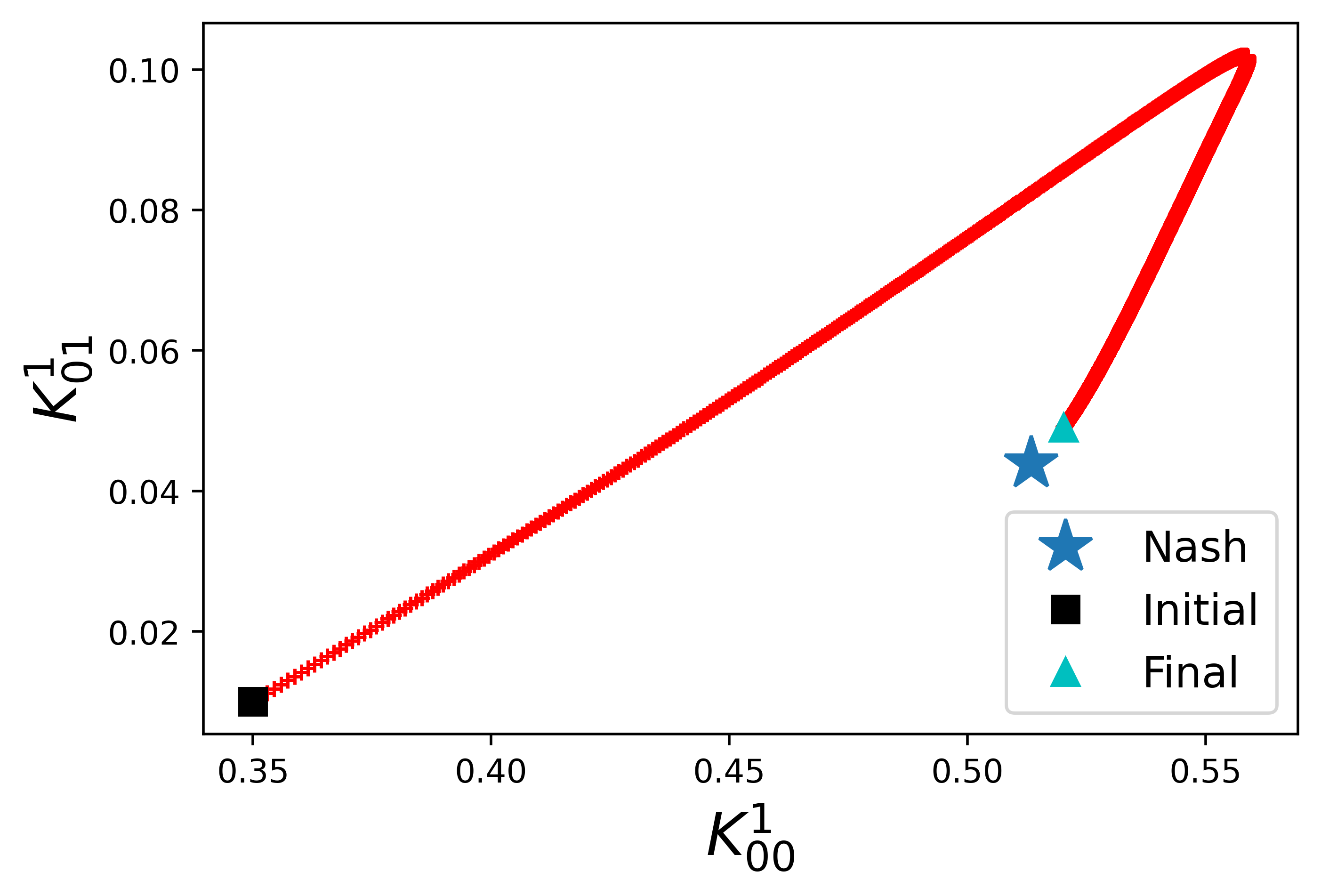}
    \caption{Trajectory of $K_{0}^{1}$.}
  \end{subfigure}
   \begin{subfigure}[b]{0.3\textwidth}
    \includegraphics[width=\textwidth]{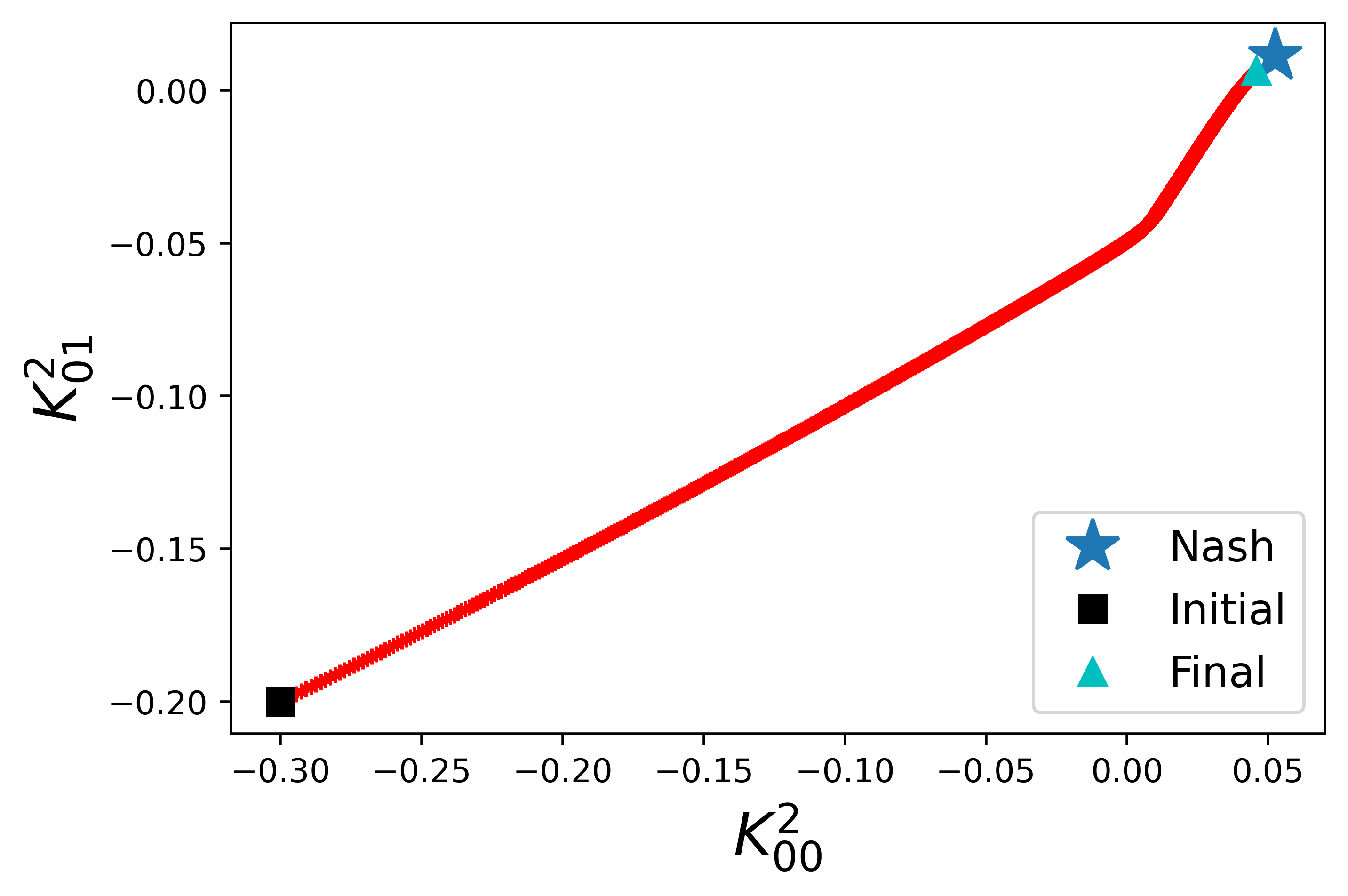}
    \caption{Trajectory of $K_{0}^{2}$.}
  \end{subfigure}
  \caption{\label{fig:NPG_div_eta}Performance of the natural policy gradient algorithm with $r=0.42$,  $\eta_1=0.001$, and $\eta_2=0.01$ ($M=20000$).}
\end{figure}

\begin{figure}[htbp]
  \centering
  \begin{subfigure}[b]{0.3\textwidth}
    \includegraphics[width=\textwidth]{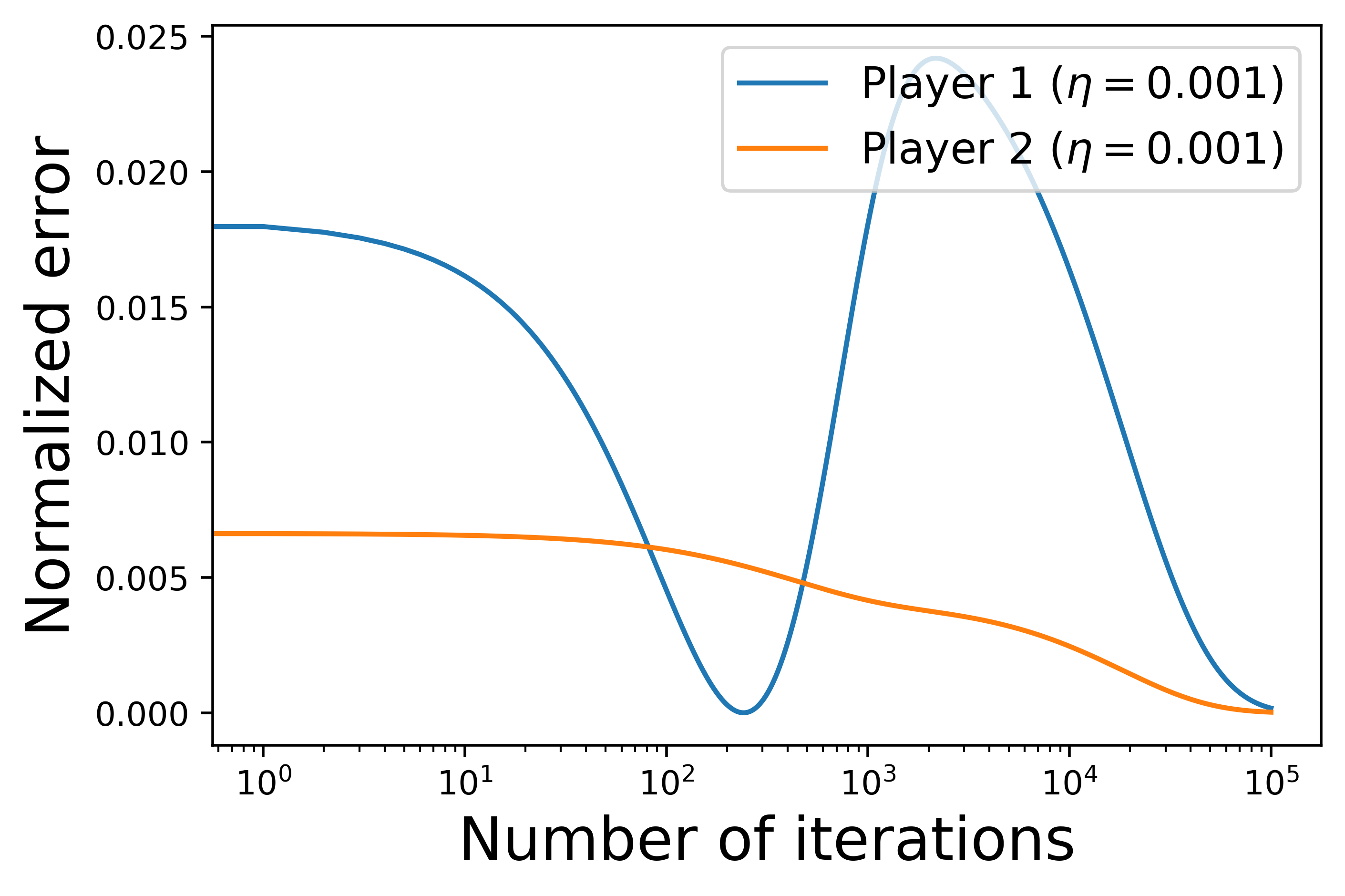}
    \caption{Normalized error.}
  \end{subfigure}
  \begin{subfigure}[b]{0.3\textwidth}
    \includegraphics[width=\textwidth]{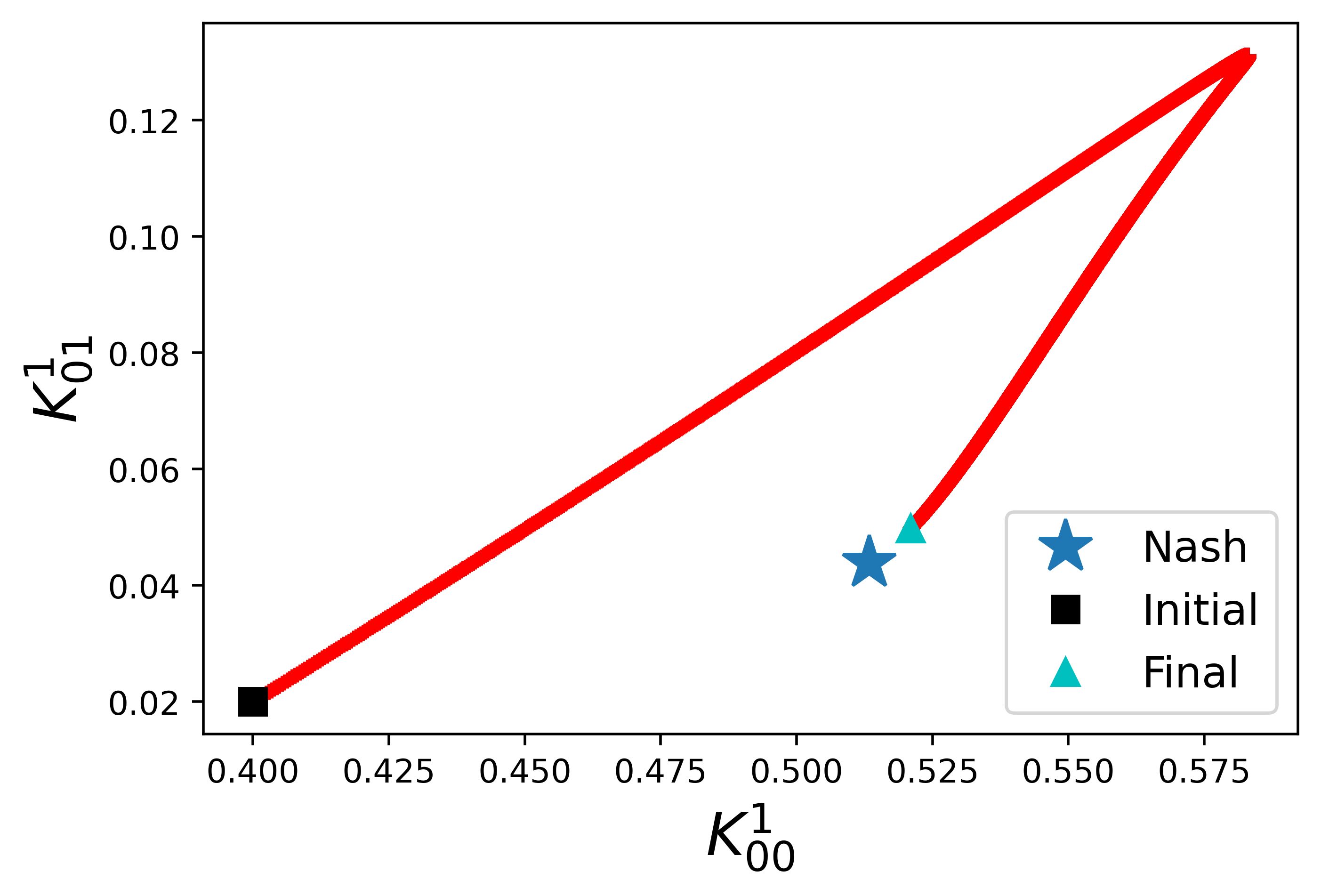}
    \caption{Trajectory of $K_{0}^{1}$.}
  \end{subfigure}
   \begin{subfigure}[b]{0.3\textwidth}
    \includegraphics[width=\textwidth]{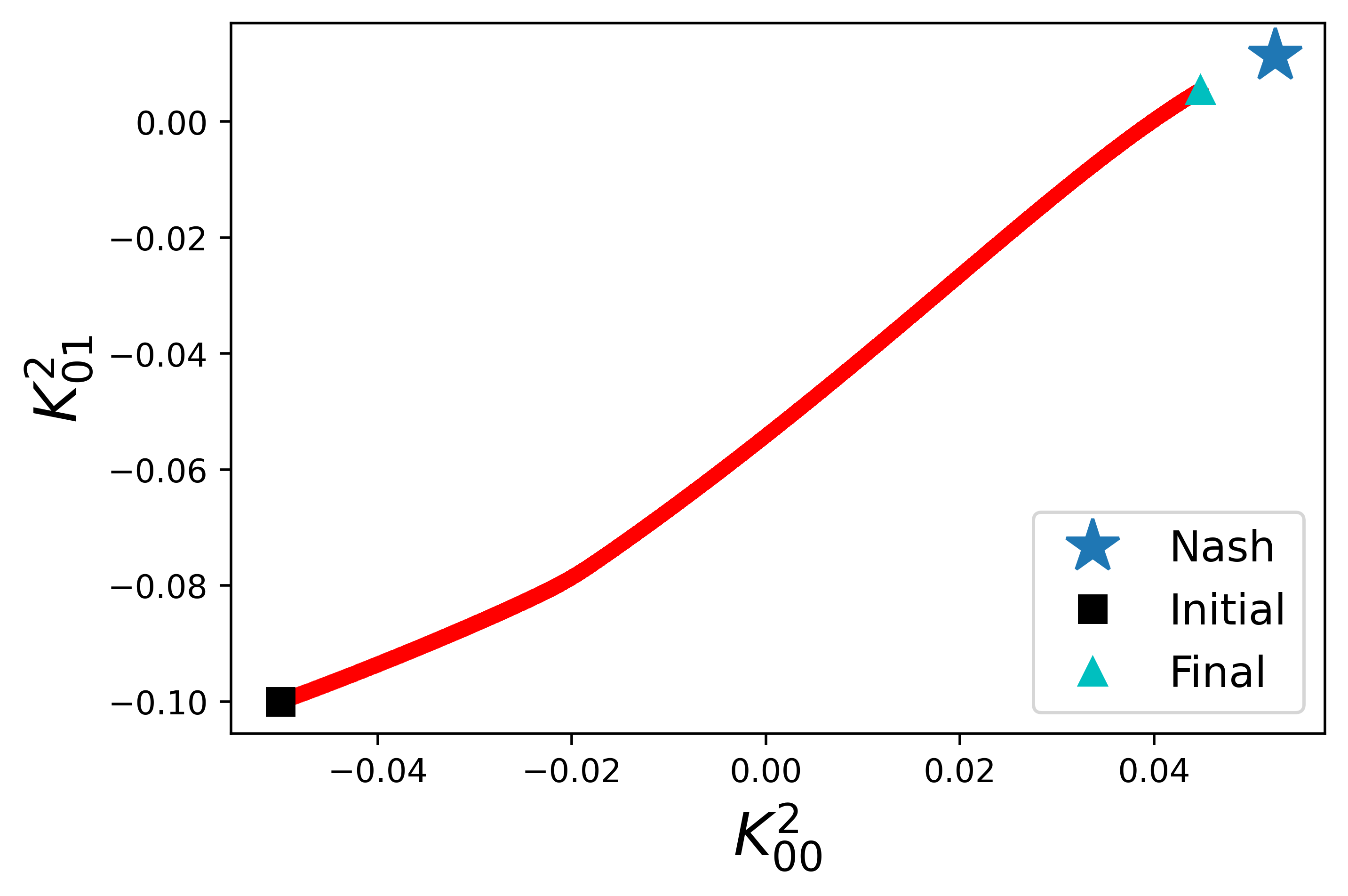}
    \caption{Trajectory of $K_{0}^{2}$.}
  \end{subfigure}
  \caption{\label{fig:NPG_div_ini}Performance of the natural policy gradient algorithm with $r=0.16$ and $\eta_1=\eta_2=0.001$ ($M=100000$).}
\end{figure}

\paragraph{Performance under Deterministic Dynamics.} We observe that under carefully chosen initial policies and step sizes, the natural policy gradient converges to the Nash equilibrium even with deterministic state dynamics ($\sigma=0$). We first show the case when the natural policy gradient algorithm diverges with $r=0.42$ and $\eta_1=\eta_2=0.001$ in Figure \ref{fig:NPG_div} (a trajectory of 10000 iterations is indicated in red). However, either by adjusting the step size to $\eta_2=0.01$ (see Figure \ref{fig:NPG_div_eta}), or by initializing the policies from a smaller neighbourhood around the Nash equilibrium (see Figure \ref{fig:NPG_div_ini}), the natural policy gradient method converges to the Nash equilibrium. {This further demonstrates that the theoretical result, along with its assumptions, in Theorem \ref{thm:local_conv_NPG} could provide 
insightful guidance on how to tune the hyper-parameters for practical examples.}

\subsection{Effect of the System Noise}\label{sec:sys_noise}

As illustrated in the theoretical analysis in Section \ref{sec:exact_NPG}, the system noise plays an important role in the convergence guarantee of the natural policy gradient algorithm. To test the sensitivity of the performance of this algorithm to the level of system noise, we apply the natural policy gradient algorithm with unknown parameters to a two-player LQ game example with synthetic data consisting of a two-dimensional state variable and a one-dimensional control variable. The model parameters (except the level of noise $\sigma^2$ in $W$ which we will discuss later) are randomly picked such that the conditions for our LQ game framework are satisfied.

\paragraph{Set-up.} We perform the natural policy gradient algorithm with synthetic data given as follows.
\begin{enumerate}
    \item Parameters: for $t=1,\cdots,T-1$,
    \[
A_t = 
\begin{bmatrix}
0.588 & 0.28\\
0.57 & 0.56
\end{bmatrix},\quad
B_t^1 = 
\begin{bmatrix}
1 \\
1
\end{bmatrix},\quad
B_t^2 = 
\begin{bmatrix}
0.5 \\
1
\end{bmatrix},\quad
W = 
\begin{bmatrix}
\sigma^2 & 0 \\
0 & \sigma^2
\end{bmatrix},
\]
    \[
Q_T^1 = 
\begin{bmatrix}
0.5 & 0\\
0 & 1
\end{bmatrix},\quad
Q_T^2 = 
\begin{bmatrix}
1 & 0\\
0 & 0.3
\end{bmatrix},\quad
Q_t^1=Q_t^2=
\begin{bmatrix}
0.001 & 0\\
0 & 0.001
\end{bmatrix},\quad
R_t^1=R_t^2=1,
\]
where $\sigma\in\mathbb{R}_+$ and $T=5$. The smoothing parameter is $r_1=r_2=0.5$ and the number of trajectories is $L_1=L_2=200$.
    \item Initialization: we assume $x_0=(x_0^1,x_0^2)$ with $x_0^2$ and $x_0^2$ independent and sampled from  $\mathcal{N}(10,2)$ and $\mathcal{N}(12,3)$ respectively. The initial policy for player 1: $\pmb{K}^1\in\mathbb{R}^{1\times 10}=(K_0^{1,(0)},\cdots,K_{T-1}^{1,(0)})$ with $K^{1,(0)}_{t}=[0.3,0.15]$ for all $t$. The initial policy for player 2: $\pmb{K}^2\in\mathbb{R}^{1\times 10}=(K_0^{2,(0)},\cdots,K_{T-1}^{2,(0)})$ with $K^{2,(0)}_{t}=[0.1,0.05]$ for all $t$.
\end{enumerate}

\begin{figure}[htbp]
  \centering
  \begin{subfigure}[h]{0.42\textwidth}
    \includegraphics[width=\textwidth]{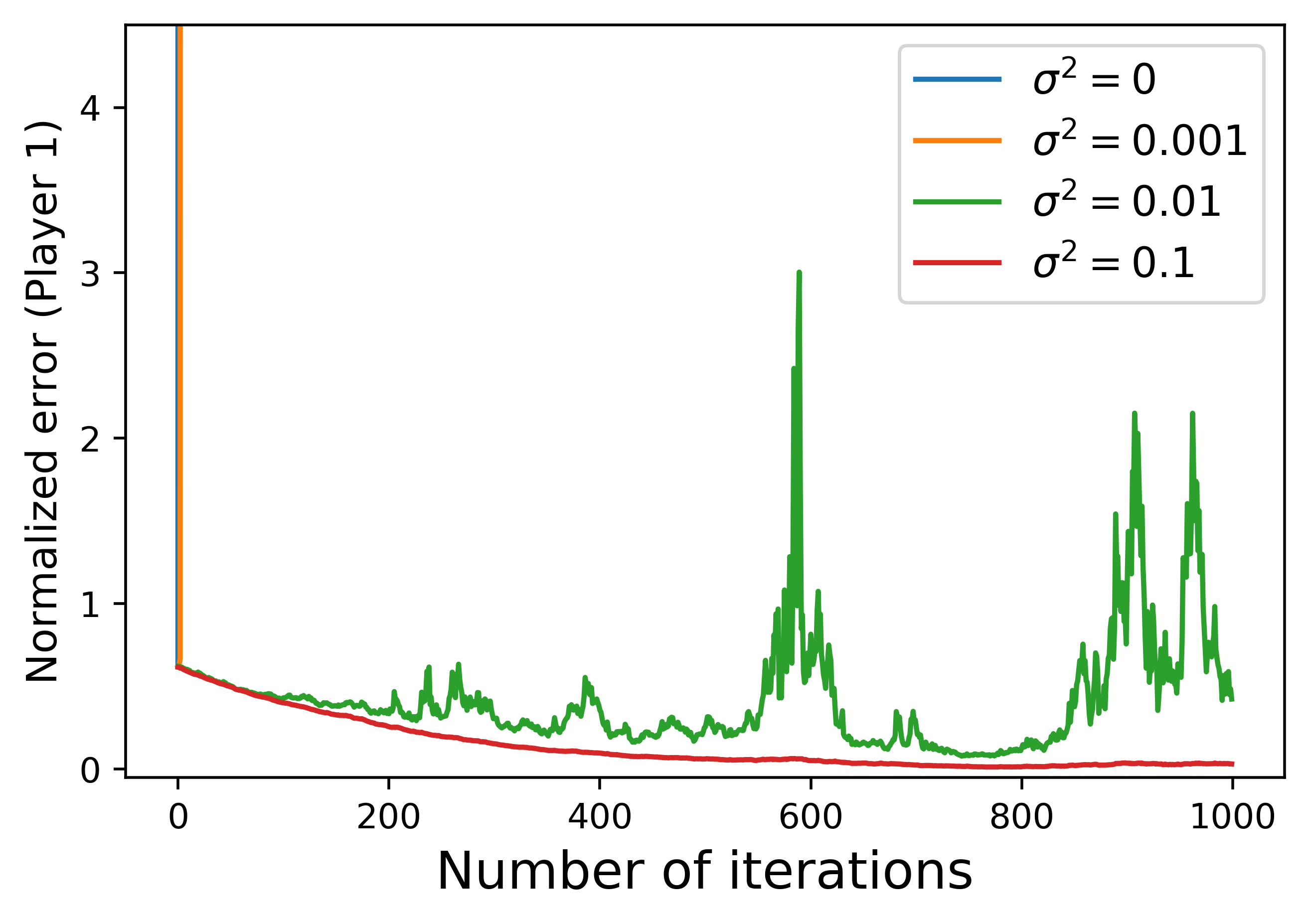}
    \caption{Player 1.}
  \end{subfigure}
  \begin{subfigure}[h]{0.42\textwidth}
    \includegraphics[width=\textwidth]{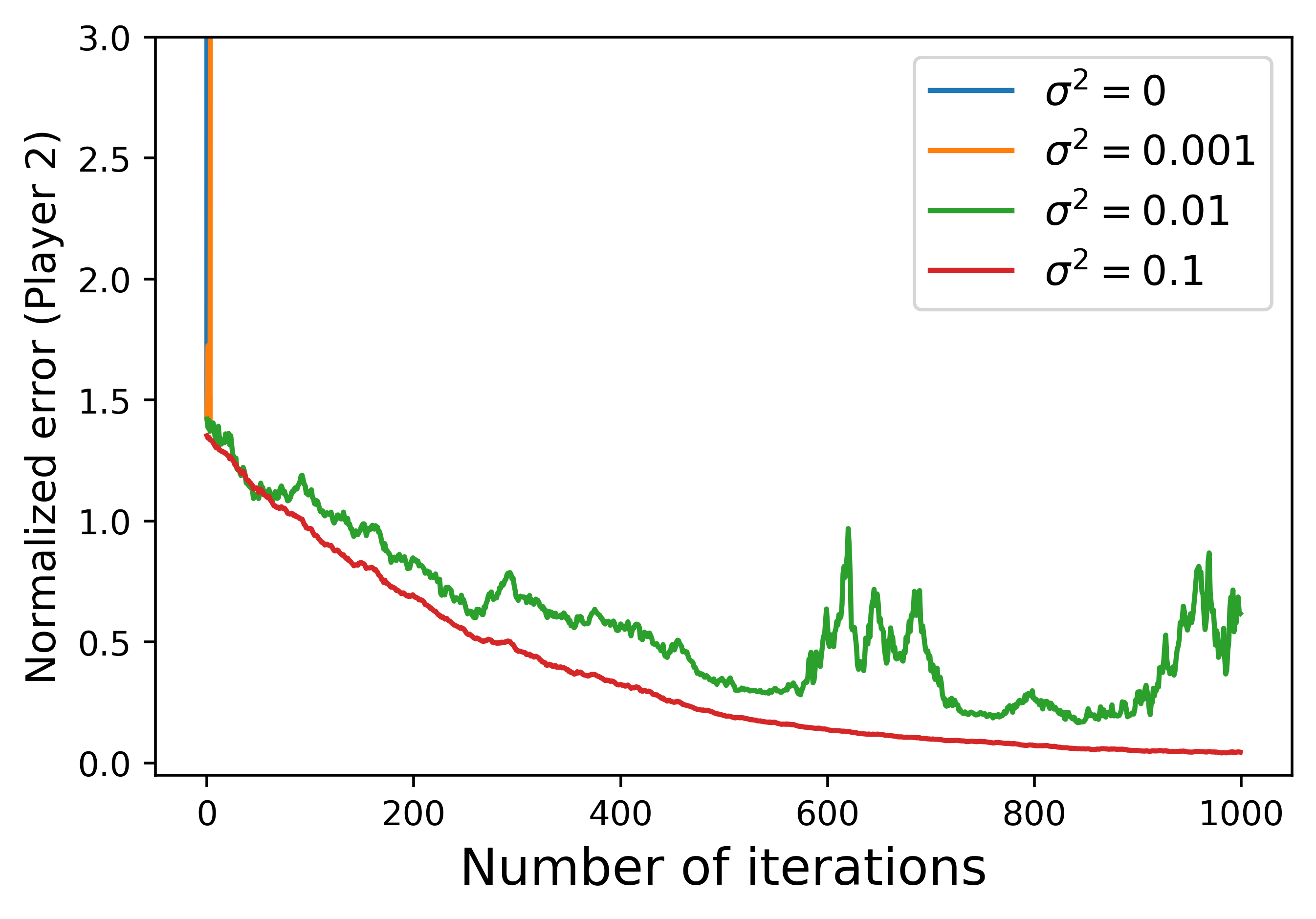}
    \caption{Player 2.}
  \end{subfigure}
\caption{\label{fig:noi_mf_conv}Performance of the natural policy gradient algorithm with unknown parameters ($\eta_1=\eta_2=0.001$). }
\end{figure}

\paragraph{Convergence.} To show that (even a low level of) the system noise can indeed help the natural policy gradient algorithm to find the Nash equilibrium, we vary the value of $\sigma^2$ from 0 to 0.1 and show the normalized error for different values of $\sigma^2$ in Figure \ref{fig:noi_mf_conv}. The natural policy gradient algorithm diverges when $\sigma^2\leq 0.001$ for both players, and it starts to converge with large fluctuations when $\sigma^2=0.01$. When $\sigma^2=0.1$, the algorithm shows a reasonable accuracy within 300 iterations (that is the normalized error is less than 5\%) for both players without fluctuations. It is worth pointing out that in fact, although Assumption \ref{ass:noise_cond} is violated when $\sigma^2=0.1$, the natural policy gradient algorithm converges to the Nash equilibrium. Finally, it is possible to make the algorithm converge when $\sigma^2\leq 0.001$ by adjusting the initial policy and the step size, as shown in Section \ref{sec:compar_infi_paper}. In practice we may not be able to change the variance of $w_t$ directly, however the level of system noise can also be increased by adding (Gaussian) explorations to the agents' policies. See \cite{houthooft2016vime,wang2020reinforcement} for more discussion of Gaussian exploration.
 
\begin{figure}[htbp]
  \centering
  \begin{subfigure}[h]{0.42\textwidth}
    \includegraphics[width=\textwidth]{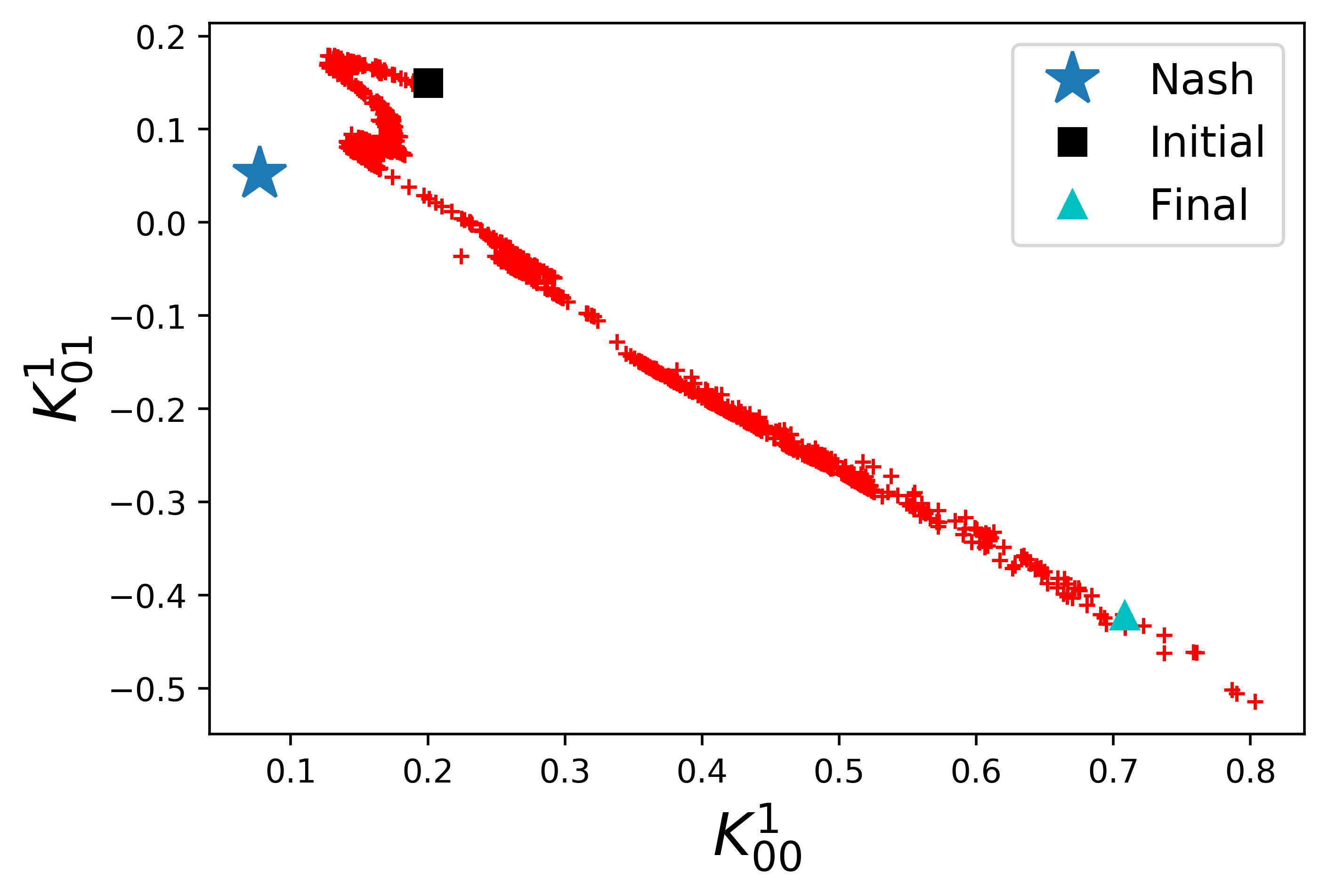}
    \caption{\label{subfig:noi_mf_traj_div}$\sigma^2=0.01$}
  \end{subfigure}
  \begin{subfigure}[h]{0.42\textwidth}
    \includegraphics[width=\textwidth]{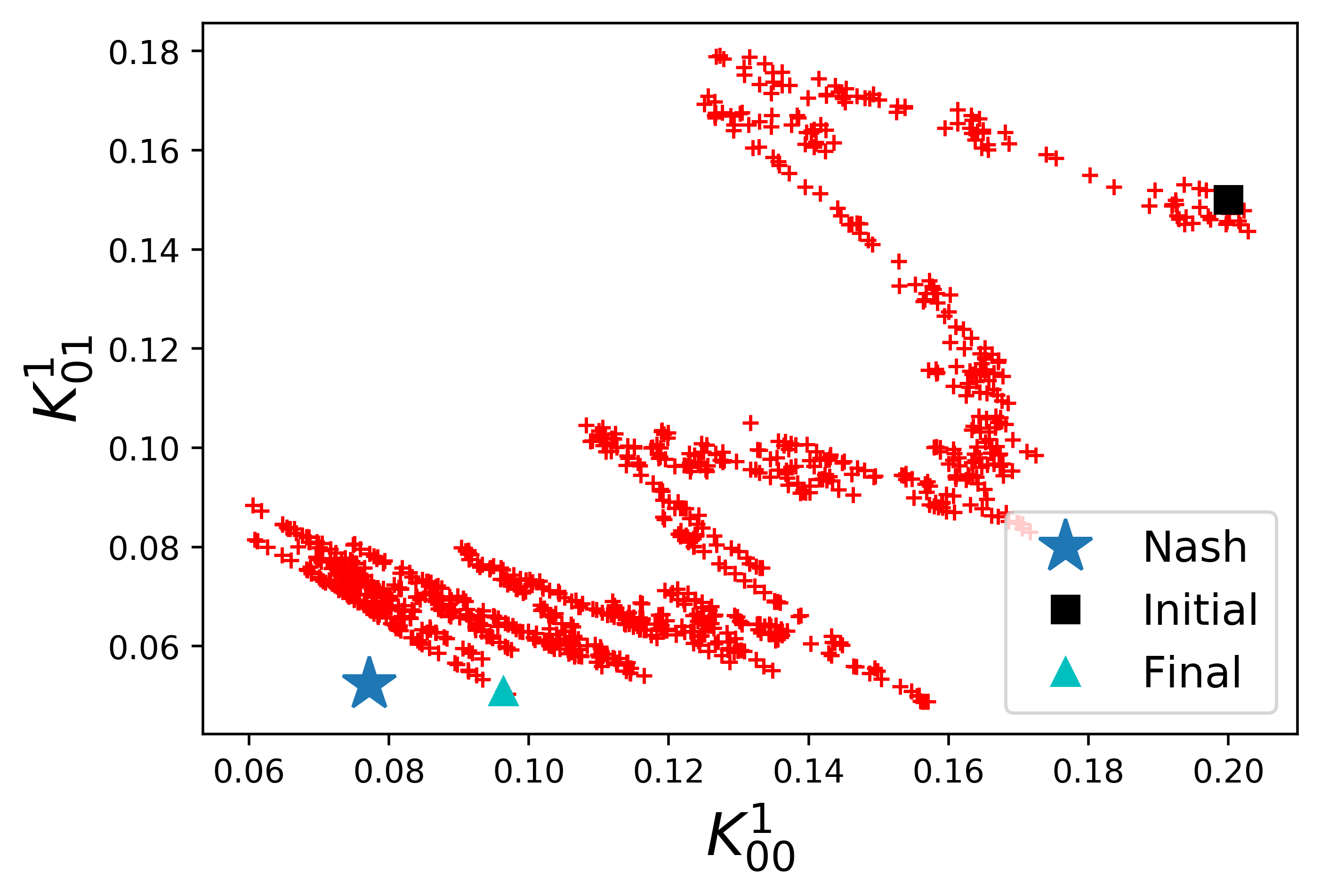}
    \caption{\label{subfig:noi_mf_traj_conv}$\sigma^2=0.1$.}
  \end{subfigure}
  \caption{\label{fig:noi_mf_traj}Trajectory (indicated in red) of learned policy $K_0^{1}=(K_{00}^1,K_{01}^1)$ ($\eta_1=\eta_2=0.001$, $M=1000$). }
\end{figure}

\paragraph{Trajectories of Learned Policies.} We show the trajectories of the learned policy at $t=0$ of player 1 by performing the natural policy gradient algorithm with unknown parameters under $\sigma^2=0.01$ and $\sigma^2=0.1$ in Figure \ref{fig:noi_mf_traj} (with 1000 iterations). In the case of a lower level of system noise ($\sigma^2=0.01$), the natural policy gradient does not converge to Nash equilibrium with $\eta_1=\eta_2=0.001$ in Figure \ref{subfig:noi_mf_traj_div}, whereas when the level of noise is increased to $\sigma^2=0.1$, the learned policy approaches the target within 1000 iterations with the same step size.

\subsection{Convergence in a Three-player Game}\label{sec:three_player}

\begin{figure}[H]
  \centering
  \begin{subfigure}[b]{0.42\textwidth}
    \includegraphics[width=\textwidth]{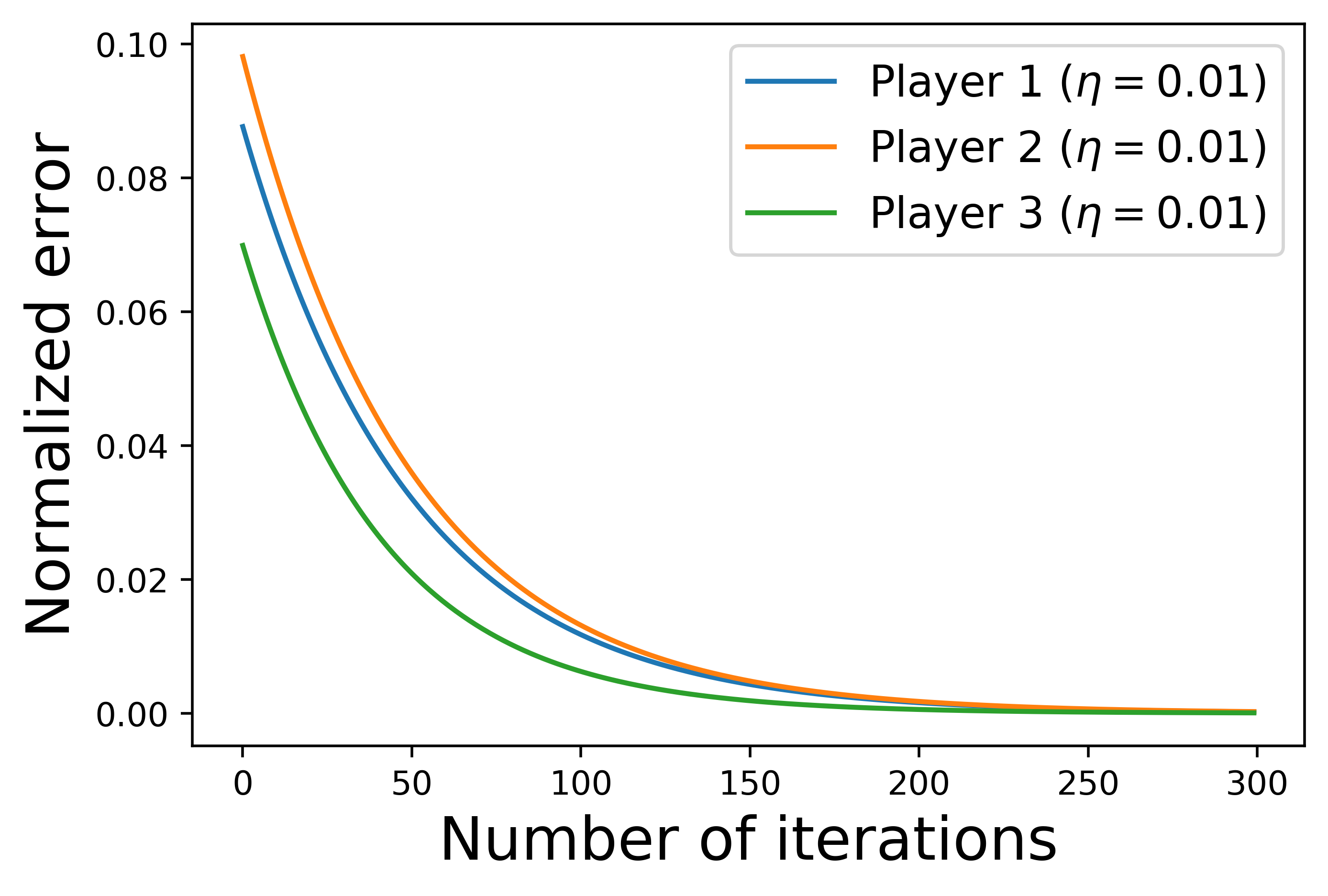}
    \caption{Known parameters.}
  \end{subfigure}
  \begin{subfigure}[b]{0.42\textwidth}
    \includegraphics[width=\textwidth]{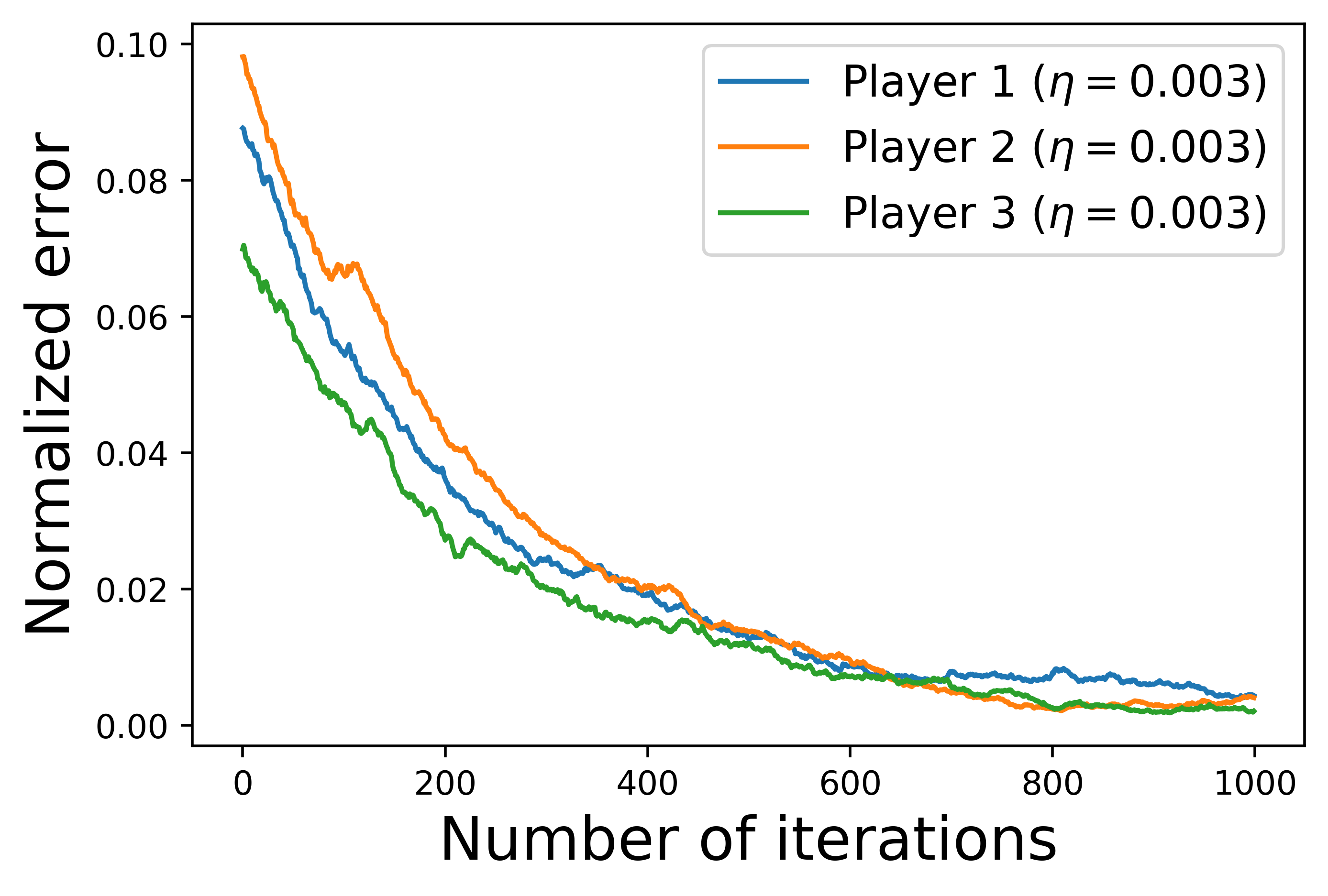}
    \caption{Unknown parameters.}
  \end{subfigure}
  \caption{\label{fig:three_player_cost}Performance of the natural policy gradient algorithm with known and unknown parameters.}
\end{figure}

In this section, we perform the the natural policy gradient method with known and unknown parameters in the following three-player general-sum game example. We show the convergence of the algorithms in Figure \ref{fig:three_player_cost}.
 
\paragraph{Set-up.} We set-up the model parameters and initial policies as follows:
\begin{enumerate}
    \item Parameters: 
    \[
A_t = 
\begin{bmatrix}
0.05 & -0.1 & 0.1\\
0.1 & 0.2 & -0.06\\
-0.02 & 0.03 &0.1
\end{bmatrix},\quad
B_t^1 = 
\begin{bmatrix}
0.05 \\
0.01 \\
-0.01
\end{bmatrix},\quad
B_t^2 = 
\begin{bmatrix}
0.01 \\
-0.05\\
-0.02
\end{bmatrix},\quad
B_t^3 = 
\begin{bmatrix}
-0.02 \\
0.01\\
0.05
\end{bmatrix},
\]
\[
W = 
\begin{bmatrix}
0.1 & 0.01 & 0.02 \\
0.01 & 0.2 & 0.01\\
0.02 & 0.01 &  0.1
\end{bmatrix},\quad
Q_T^1=Q_T^2 =Q_t^1=Q_t^2 = 
\begin{bmatrix}
0.2 & 0 & 0\\
0 & 0.2 & 0 \\
0 & 0 & 0.2 \\
\end{bmatrix},
\]
$R_t^1(t)=R_t^2(t)=0.5$, $R_t^3(t)=0.6$, and $T=5$.
    \item Initialization: Take $x_0=(x_0^1,x_0^2,x_0^3)$ where $x_0^1, x_0^2, x_0^3$ are independent and sampled from  $\mathcal{N}(0.3,0.2)$ and $\mathcal{N}(0.2,0.3)$, and $\mathcal{N}(0.3,0.2)$ respectively. The initial policies are $\pmb{K}^{1,(0)}=(0.35,0.01,0.1)$, $\pmb{K}^{2,(0)}=(-0.3,-0.2,0)$, and $\pmb{K}^{3,(0)}=(-0.3,0.1,0)$.
\end{enumerate} 

\paragraph{Convergence.} We plot the normalized error for each player in the case of known parameters and also unknown parameters in Figure \ref{fig:three_player_cost}. We can see that with three players, the algorithms still have a reasonably fast speed of convergence in practice.

\section*{Acknowledgement}
We thank Paul Barnes and Amr El Zanfally at BP plc for providing part of the motivation for this work and for their continued support.
%\newpage
%\section{Discussion}
%\paragraph{Global convergence.}

%\paragraph{Theoretical guarantee for the model-free method.}

%\newpage
\bibliographystyle{plain}
\bibliography{references}

\newpage
\appendix
\input{appendix}

\end{document}

%% file: appendix.tex
\section{The One-step Contraction Lemma for the Vanilla Policy Gradient Method}\label{app:vani}
The convergence result for the natural policy gradient method can be extended to the case of the vanilla policy gradient method. The key step is to prove the one-step contraction (Lemma \ref{lemma:one_step_loc_conv}) for the vanilla version. This can be done by modifying some parts of the current Lemma \ref{lemma:one_step_loc_conv}. We first recall the definition of $g_1$ and $g_2$ as follows:
\begin{equation*}
    g_1 : = \frac{\srmin}{\|\Sigma_{\pmb{K}^{*}}\|},
\end{equation*}
and
\begin{equation*}
    g_2 := 20(N-1)^2\,T^2\,d\, \frac{(\gamma_B)^4\max_i\{C^i(\pmb{K}^{i},\pmb{K}^{-i*})\}^4}{\sqmin^2\srmin}\left(\frac{\rho_{\pmb{K}}^{2T}-1}{\rho_{\pmb{K}}^2-1}\right)^2.
\end{equation*}
We further define $\widetilde{g}_2$ as
\begin{eqnarray}\label{defn_g2_tilde}
   \widetilde{g}_2 \hspace{-0.2cm}&:=&\hspace{-0.3cm} g_2\left(\rho_{\pmb{K}}^{2T} \|\Sigma_0\| + \big(\rho_{\pmb{K}}^{2T}+1\big)\|W\| \right)^2\nonumber\\
   &=&\hspace{-0.3cm}20(N-1)^2\,T^2\,d\, \frac{(\gamma_B)^4\max_i\{C^i(\pmb{K}^{i},\pmb{K}^{-i*})\}^4}{\sqmin^2\srmin}\left(\frac{\rho_{\pmb{K}}^{2T}-1}{\rho_{\pmb{K}}^2-1}\right)^2\hspace{-0.2cm}\left(\rho_{\pmb{K}}^{2T} \|\Sigma_0\| + \big(\rho_{\pmb{K}}^{2T}+1\big)\|W\| \right)^2,
\end{eqnarray}
and $g_3$ as
\begin{equation}\label{eqn:defn_g3}
    g_3:=\frac{5\max_i\{C^i(\pmb{K}^{i},\pmb{K}^{-i*})\}}{\sqmin}\Big(\frac{\max_i\{C^i(\pmb{K}^{i},\pmb{K}^{-i*})\}}{\sqmin}+\rho_{\pmb{K}}^{2T} \|\Sigma_0\| + \big(\rho_{\pmb{K}}^{2T}+1\big)\|W\| \Big).
\end{equation}
We also write 
$C^{i,-i*}=C^i(\pmb{K}^i,\pmb{K}^{-i*})$, $C^{i*} = C^i(\pmb{K}^{*})$ and $C^{i\prime,-i*}=C^i(\pmb{K}^{i\prime},\pmb{K}^{-i*})$ to simplify notation.

% {\color{blue}Should we give a formal Theorem for the vanilla case?}

\begin{Lemma}[One-step contraction for vanilla policy gradient]\label{lemma:one_step_cont_vani}
Assume Assumptions \ref{ass:nonzero_cost}, \ref{ass:nonzero_initial_noise}, and  \ref{ass:nonzero_exist_sol} hold, and that 
\begin{equation}\label{eqn:one_step_noise_cond_vani}
\sx^7 > \max\left\{\frac{\widetilde{g}_2}{g_1},g_3^{\,\, \frac{7}{2}}\right\}.
\end{equation}
Also assume the policy update step for player $i$ at time $t$ is given by
\begin{equation}\label{eqn:one_step_update_vani}
    K_t^{i\prime} = K_t^i - \eta\nabla_{K_t^i}C^i(\pmb{K}),
\end{equation}
where
\begin{equation}\label{eqn:step_size_vani}
    \eta \leq \min\left\{I_3\,,I_4\,,\frac{\|\Sigma_{\pmb{K}^{*}}\|}{\sx^2\srmin}\right\}
\end{equation}
with
 \begin{equation*}
 \begin{split}
 I_3 & = \left\{20 T\sx\frac{\rho_{\pmb{K}}(\rho_{\pmb{K}}^{2T}-1)}{\rho_{\pmb{K}}^2-1}\left(\sum_{i=1}^NC^{i,-i*} +\sqmin T\|W\|\right)\gamma_B\max_i\{\max_{t}\{\|\nabla_{K_t^i}C^i(\pmb{K})\|\}\}+2\sqmin\sx\right.\\
 &\quad \left.+8\big(\gamma_R+\frac{(\gamma_B)^2}{\sx}\sum_{i=1}^NC^{i,-i*}\big)\left(\rho_{\pmb{K}}^{2T} \|\Sigma_0\| + \big(\rho_{\pmb{K}}^{2T}+1\big)\|W\| \right)^2\sum_{i=1}^N\{C^{i,-i*}\}\right\}^{-1}\cdot\sqmin(\sx)^2,\\
 I_4 &= \left\{\big(\max_i\{k_i\}\big)\frac{10 T\sum_{i=1}^NC^{i,-i*}}{(10 T-1)\sqmin}\left(\gamma_R+(\gamma_B)^2\frac{\sum_{i=1}^NC^{i,-i*}}{\sx}\right)+\frac{2d}{\sx}\left(\frac{10 T\sum_{i=1}^NC^{i,-i*}}{(10 T-1)\sqmin}\right)^2\right.\cdot\\
 &\quad\quad\left.\left(\gamma_R+(\gamma_B)^2\frac{\sum_{i=1}^NC^{i,-i*}}{\sx}\right)^2\left(\rho_{\pmb{K}}^{2T} \|\Sigma_0\| + \big(\rho_{\pmb{K}}^{2T}+1\big)\|W\| \right)^2\right\}^{-1}\cdot\frac{d}{80\sx^2} \left(\frac{10 T\min_i\{C^{i*}\}}{(10 T-1)\sqmin}\right)^2.
 \end{split}
 \end{equation*}
Let $\alpha:=\sx^2 g_1- \widetilde{g}_2/\sx^5>0$. Then, we have 
\begin{enumerate}
    \item  $\eta \in (0,\frac{1}{\alpha})$; and
    \item the following inequality holds
\begin{equation}\label{eq:c_sum_onestep_final_vani} \sum_{i=1}^N \left(C^{i}(\pmb{
K}^{i\prime},\pmb{K}^{-i*})-C^{i}(\pmb{K}^*)\right)
\leq (1-\alpha\eta)\left(\sum_{i=1}^N \big(C^{i}(\pmb{K}^i,\pmb{K}^{-i*}) - C^{i}(\pmb{K}^*)\big)\right).
\end{equation}
\end{enumerate}
\end{Lemma}

\begin{Remark}\label{remark:comp_noi_cond}{\rm
We compare the noise condition \eqref{eqn:one_step_noise_cond_vani} for the vanilla policy gradient method and the condition \eqref{eqn:one_step_noise_cond} for the natural version in Lemma \ref{lemma:one_step_loc_conv}. We mainly focus on the orders of $N$, $T$, and $\rho_{\pmb{K}}$, and ignore other constants in \eqref{eqn:one_step_noise_cond_vani} and \eqref{eqn:one_step_noise_cond}. For the natural version, we need
\begin{equation}\label{noise_order_natu}
    \sx > \left(\frac{g_2}{g_1}\right)^{1/5}=\mathcal{O}\left((N-1)^{2/5}\,T^{2/5}\,\rho_{\pmb{K}}^{4T/5}\right),
\end{equation}
and for the vanilla version, we need
\begin{equation}\label{noise_order_vani}
    \sx> \max\left\{\frac{\widetilde{g}_2}{g_1},g_3^{\,\, \frac{7}{2}}\right\}^{1/7}=\mathcal{O}\left((N-1)^{2/7}\,T^{2/7}\,\rho_{\pmb{K}}^{8T/7}\right).
\end{equation}
The order of $\rho_{\pmb{K}}$ is higher in \eqref{noise_order_vani}, and the orders of $T$ and $(N-1)$ are slightly higher in \eqref{noise_order_natu}. Thus when $\rho_{\pmb{K}}^{T}$ is small, the vanilla method has a weaker noise assumption. This may happen when the time horizon $T$ is small and the policy $\pmb{K}$ is close to the Nash equilibrium. In contrast, when $\rho_{\pmb{K}}^{T}$ is large, \eqref{noise_order_natu} is weaker than \eqref{noise_order_vani} and the natural method is superior to the vanilla method. Additionally, when $N-1$ is very large (and $\rho_{\pmb{K}}^{T}$ does not blow up) ,  \eqref{noise_order_vani} leads to a weaker assumption.

We also note that other than the noise condition, there are also some slight differences between the step size conditions \eqref{eqn:step_size_vani} for the vanilla method and \eqref{eqn:step_size} for the natural method, mainly in the order of $\rho_{\pmb{K}}^{2T}$ appearing in the denominator of $I_1, I_3$, and $I_4$.}
\end{Remark}

\begin{proof}[Proof of Lemma~\ref{lemma:one_step_cont_vani}.] We break this proof up into a series of steps. \\ [.1in]
{\it Step 1:} We first consider the consequences of the condition  $\eta\leq \min\{I_3,I_4\}$. Straightforward calculations show that when condition $\eta\leq I_3$ is satisfied, the following inequalities hold:
\begin{enumerate}
    \item $\forall i=1,\cdots,N$,
    \begin{eqnarray}\label{eqn:I1_conseq_3_vani}
    \|K_t^{i\prime}-K_t^i\| = \eta \|\nabla_{K_t^i} C^i(\pmb{K})\| \leq \frac{\sqmin\sx}{20 T \gamma_B C^{i,-i*}}.    
    \end{eqnarray}
    \item $\forall i=1,\cdots,N$, \begin{eqnarray}\label{eqn:I1_conseq_1_vani}
    &&\eta \left(\frac{ \rho_{\pmb{K}}^{2T}-1}{\rho_{\pmb{K}}^2-1} \left(\frac{C^{i,-i*}}{\sqmin}+T\|W\|\right)2\rho_{\pmb{K}}\,\gamma_B\sum_{t=0}^{T-1}\|\nabla_{K_t^i}C^i(\pmb{K})\|\right)\nonumber\\
    &\leq& I_3 \left(\frac{ \rho_{\pmb{K}}^{2T}-1}{\rho_{\pmb{K}}^2-1} \left(\frac{\sum_{i=1}^N C^{i,-i*}}{\sqmin}+T\|W\|\right)2\rho_{\pmb{K}}\,\gamma_BT\max_t\{\|\nabla_{K_t^i}C^i(\pmb{K})\|\}\right)\nonumber \\
    &\leq& \frac{\sx}{10}.
\end{eqnarray}
\item $\forall i=1,\cdots,N$, \begin{eqnarray}\label{eqn:I1_conseq_2_vani}
\eta\leq I_3&\leq& \frac{\sx^2}{2\sx+4\frac{2C^{i,-i*}}{\sqmin}\left(\rho_{\pmb{K}}^{2T} \|\Sigma_0\| + \big(\rho_{\pmb{K}}^{2T}+1\big)\|W\| \right)^2(\gamma_R+\gamma_B^2\frac{C^{i,-i*}}{\sx})}.
\end{eqnarray}
\end{enumerate}
In the case where  $\eta\leq I_4$, we have  $\forall i=1,\cdots,N$
\begin{eqnarray}\label{eqn:I2_conseq_vani}
&&4\eta\,k_i\frac{10 TC^{i,-i*}}{(10 T-1)\sqmin}\left(\gamma_R+(\gamma_B)^2\frac{C^{i,-i*}}{\sx}\right)\nonumber\\
&& +8\eta\frac{d}{\sx}\left(\frac{10 TC^{i,-i*}}{(10 T-1)\sqmin}\right)^2\left(\gamma_R+(\gamma_B)^2\frac{C^{i,-i*}}{\sx}\right)^2\left(\rho_{\pmb{K}}^{2T} \|\Sigma_0\| + \big(\rho_{\pmb{K}}^{2T}+1\big)\|W\| \right)^2\nonumber\\ 
&\leq& 4 I_4 \left(\,k_i\frac{10 TC^{i,-i*}}{(10 T-1)\sqmin}\left(\gamma_R+(\gamma_B)^2\frac{C^{i,-i*}}{\sx}\right)\right.\nonumber\\
&&\left.\qquad+2\frac{d}{\sx}\left(\frac{10 TC^{i,-i*}}{(10 T-1)\sqmin}\right)^2\left(\gamma_R+(\gamma_B)^2\frac{C^{i,-i*}}{\sx}\right)^2\left(\rho_{\pmb{K}}^{2T} \|\Sigma_0\| + \big(\rho_{\pmb{K}}^{2T}+1\big)\|W\| \right)^2\right)\nonumber\\
&\leq& \frac{d}{20\sx^2} \left(\frac{10 T\min_i\{C^{i*}\}}{(10 T-1)\sqmin}\right)^2 \leq \frac{d}{20\sx^2} \left(\frac{10 T\max_i\{C^{i,-i*}\}}{(10 T-1)\sqmin}\right)^2.
\end{eqnarray}
By \eqref{eqn:I1_conseq_3_vani} and Lemma~\ref{lemma:nonzero_bds_P_Sigma} we have
\begin{equation*}
    \gamma_B\|K_t^{i\prime}-K_t^i\|  \leq \frac{\sqmin\sx}{20 TC^{i,-i*}} \leq \frac{1}{20 T^2}.
\end{equation*}
Therefore, we have $\rho_{\pmb{K},\pmb{K}^{\prime}}\leq \rho_{\pmb{K}}$ by Lemma \ref{lemma:rho_upper_bd}. \\[.1in]

{\it Steps 2 and 3:} The results in Steps 2 and 3 in Lemma \ref{lemma:one_step_loc_conv} for the natural policy gradient method still hold for the vanilla policy gradient method by the consequences \eqref{eqn:I1_conseq_3_vani} and \eqref{eqn:I1_conseq_1_vani}. Here we omit the proof and state the following results which will be used in Steps 4 and 5:
\begin{equation}\label{eqn:bd_pertSigma_vani}
   \big\|\Sigma_{\pmb{K}^{i\prime},\pmb{K}^{-i*}}\big\| \leq\frac{10 T\,C^{i,-i*}}{(10 T-1)\sqmin},
\end{equation}
and
\begin{eqnarray}
     \left\|E_{t,i}^{\pmb{K}}-E_{t,i}^{\pmb{K}^i,\pmb{K}^{-i*}}\right\|\leq \frac{(\gamma_B)^2C^{i,-i*}}{\sx}\left(\frac{2(\rho_{\pmb{K}}^{2T}-1)}{\rho_{\pmb{K}}^2-1}\sum_{j=1,j\neq i}^N\vertiii{\pmb{K}^{j}-\pmb{K}^{j*}}\right)\label{eqn:diff_E_bd_vani}.
\end{eqnarray}
\\[.1in]

{\it Step 4:} We can now estimate the cost difference between using $\pmb{K}^i$ and the update $\pmb{K}^{i\prime}$. 
By By Lemma~\ref{lemma:nonzero_almost_smoothness} we have
\begin{eqnarray}
     &&C^{i\prime,-i*}-C^{i,-i*}\nonumber\\
     && = \sum_{t=0}^{T-1}\Big[\Tr\big(\Sigma_t^{\pmb{K}^{i\prime},\pmb{K}^{-i*}} (K_t^{i\prime}-K_t^i)^{\top}(R_t^{i}+(B_t^i)^{\top}P_{t +1,i}^{\pmb{K}^{i},\pmb{K}^{-i*}}B_t^i)(K_t^{i\prime}-K_t^i)\big)\nonumber\\
     &&+2\Tr\big(\Sigma_t^{\pmb{K}^{i\prime},\pmb{K}^{-i*}} (K_t^{i\prime}-K_t^i)^{\top}E_{t,i}^{\pmb{K}^i,\pmb{K}^{-i*}}\big)\Big]. \label{eqn:Cdiff_inte1_vani}
\end{eqnarray}     
For the vanilla policy gradient method, we have the following update rule
\[
K_t^{i'}=K_t^i - \eta\nabla_{K_t^i}C^i(\pmb{K})=K_t^i - 2\eta\,E_{t,i}^{\pmb{K}}\,\Sigma_t^{\pmb{K}}
\]
by Lemma~\ref{lemma:nonzero_policygrad}. Then plugging in $K_t^{i'}-K_t^i=- 2\eta\,E_{t,i}^{\pmb{K}}\,\Sigma_t^{\pmb{K}}$ into \eqref{eqn:Cdiff_inte1_vani} leads to
\begin{eqnarray}
&& C^{i\prime,-i*}-C^{i,-i*} \nonumber \\
&& = \sum_{t=0}^{T-1}\Big[4\eta^2\Tr\big(\Sigma_t^{\pmb{K}^{i\prime},\pmb{K}^{-i*}}\Sigma_t^{\pmb{K}} (E_{t,i}^{\pmb{K}})^\top(R_t^{i}+(B_t^i)^{\top}P_{t +1,i}^{\pmb{K}^{i},\pmb{K}^{-i*}}B_t^i)E_{t,i}^{\pmb{K}}\Sigma_t^{\pmb{K}}\big)\nonumber\\
&&\quad\quad\quad-4\eta\Tr\big(\Sigma_t^{\pmb{K}^{i\prime},\pmb{K}^{-i*}} \Sigma_t^{\pmb{K}}(E_{t,i}^{\pmb{K}})^\top E_{t,i}^{\pmb{K}^i,\pmb{K}^{-i*}}\big)\Big]\nonumber\\
&& = \sum_{t=0}^{T-1}\Big[4\eta^2\Tr\Big(\Sigma_t^{\pmb{K}^{i\prime},\pmb{K}^{-i*}}\Sigma_t^{\pmb{K}} (E_{t,i}^{\pmb{K}}-E_{t,i}^{\pmb{K}^i,\pmb{K}^{-i*}}+E_{t,i}^{\pmb{K}^i,\pmb{K}^{-i*}})^\top(R_t^{i}+(B_t^i)^{\top}P_{t +1,i}^{\pmb{K}^{i},\pmb{K}^{-i*}}B_t^i)\cdot\nonumber\\
&&\quad\quad\quad(E_{t,i}^{\pmb{K}}-E_{t,i}^{\pmb{K}^i,\pmb{K}^{-i*}}+E_{t,i}^{\pmb{K}^i,\pmb{K}^{-i*}})\Sigma_t^{\pmb{K}}\Big)-4\eta\Tr\Big(\Sigma_t^{\pmb{K}^{i\prime},\pmb{K}^{-i*}} (E_{t,i}^{\pmb{K}}\Sigma_t^{\pmb{K}}-E_{t,i}^{\pmb{K}^i,\pmb{K}^{-i*}}\Sigma_t^{\pmb{K}^{i\prime},\pmb{K}^{-i*}}\nonumber\\
&& \quad\quad\quad+E_{t,i}^{\pmb{K}^i,\pmb{K}^{-i*}}\Sigma_t^{\pmb{K}^{i\prime},\pmb{K}^{-i*}})^\top E_{t,i}^{\pmb{K}^i,\pmb{K}^{-i*}}\Big)\Big]\nonumber\\
&& = \sum_{t=0}^{T-1}\Big[4\eta^2\Tr\big(\Sigma_t^{\pmb{K}^{i\prime},\pmb{K}^{-i*}}\Sigma_t^{\pmb{K}} (E_{t,i}^{\pmb{K}}-E_{t,i}^{\pmb{K}^i,\pmb{K}^{-i*}})^\top(R_t^{i}+(B_t^i)^{\top}P_{t +1,i}^{\pmb{K}^{i},\pmb{K}^{-i*}}B_t^i)(E_{t,i}^{\pmb{K}}-E_{t,i}^{\pmb{K}^i,\pmb{K}^{-i*}})\Sigma_t^{\pmb{K}}\big)\nonumber\\
&&\quad\quad+ 8\eta^2\Tr\big(\Sigma_t^{\pmb{K}^{i\prime},\pmb{K}^{-i*}}\Sigma_t^{\pmb{K}} (E_{t,i}^{\pmb{K}}-E_{t,i}^{\pmb{K}^i,\pmb{K}^{-i*}})^\top(R_t^{i}+(B_t^i)^{\top}P_{t +1,i}^{\pmb{K}^{i},\pmb{K}^{-i*}}B_t^i)E_{t,i}^{\pmb{K}^i,\pmb{K}^{-i*}}\Sigma_t^{\pmb{K}}\big) \nonumber\\
&& \quad\quad+4\eta^2\Tr\big(\Sigma_t^{\pmb{K}^{i\prime},\pmb{K}^{-i*}}\Sigma_t^{\pmb{K}} (E_{t,i}^{\pmb{K}^i,\pmb{K}^{-i*}})^\top(R_t^{i}+(B_t^i)^{\top}P_{t +1,i}^{\pmb{K}^{i},\pmb{K}^{-i*}}B_t^i)E_{t,i}^{\pmb{K}^i,\pmb{K}^{-i*}}\Sigma_t^{\pmb{K}}\big)\nonumber\\
&&\quad\quad-4\eta\Tr\big(\Sigma_t^{\pmb{K}^{i\prime},\pmb{K}^{-i*}} (E_{t,i}^{\pmb{K}}\Sigma_t^{\pmb{K}}-E_{t,i}^{\pmb{K}^i,\pmb{K}^{-i*}}\Sigma_t^{\pmb{K}^{i\prime},\pmb{K}^{-i*}})^\top E_{t,i}^{\pmb{K}^i,\pmb{K}^{-i*}}\big)\nonumber\\
&&\quad\quad-4\eta\Tr\big(\Sigma_t^{\pmb{K}^{i\prime},\pmb{K}^{-i*}} (E_{t,i}^{\pmb{K}^i,\pmb{K}^{-i*}})^\top E_{t,i}^{\pmb{K}^i,\pmb{K}^{-i*}}\Sigma_t^{\pmb{K}^{i\prime},\pmb{K}^{-i*}}\big)\Big],\label{eqn:vanilla_inte1}
\end{eqnarray}
where the first equation holds by the updating rule, the second equation holds by adding and subtracting $E_{t,i}^{\pmb{K}^i,\pmb{K}^{-i*}}$ and $E_{t,i}^{\pmb{K}^i,\pmb{K}^{-i*}}\Sigma_t^{\pmb{K}^{i\prime},\pmb{K}^{-i*}}$ terms, and the third equation holds by expanding terms.
Now, letting $\omega^2=\frac{2}{\sx}$ in
\begin{equation}\label{eqn:trace_upper_bound}
    2\Tr(A^\top B)=\Tr(A^\top B + B^\top A)\leq \omega^2 \Tr(A^\top A)+\frac{1}{\omega^2}\Tr(B^\top B),
\end{equation}
(which holds for any matrices $A$ and $B$ of the same dimension), the second term in \eqref{eqn:vanilla_inte1} can be bounded by
\begin{eqnarray}
&&8\eta^2\Tr\big(\Sigma_t^{\pmb{K}^{i\prime},\pmb{K}^{-i*}}\Sigma_t^{\pmb{K}} (E_{t,i}^{\pmb{K}}-E_{t,i}^{\pmb{K}^i,\pmb{K}^{-i*}})^\top(R_t^{i}+(B_t^i)^{\top}P_{t +1,i}^{\pmb{K}^{i},\pmb{K}^{-i*}}B_t^i)E_{t,i}^{\pmb{K}^i,\pmb{K}^{-i*}}\Sigma_t^{\pmb{K}}\big) \nonumber\\
&\leq& 8\eta^2\frac{\sx}{4}\Tr\big((E_{t,i}^{\pmb{K}^i,\pmb{K}^{-i*}})^\top E_{t,i}^{\pmb{K}^i,\pmb{K}^{-i*}}\big)+ 8\eta^2\frac{1}{\sx}\Tr\big(\Sigma_t^{\pmb{K}}\Sigma_t^{\pmb{K}^{i\prime},\pmb{K}^{-i*}}\Sigma_t^{\pmb{K}} (E_{t,i}^{\pmb{K}}-E_{t,i}^{\pmb{K}^i,\pmb{K}^{-i*}})^\top\cdot \nonumber\\
&&(R_t^{i}+(B_t^i)^{\top}P_{t +1,i}^{\pmb{K}^{i},\pmb{K}^{-i*}}B_t^i)(R_t^{i}+(B_t^i)^{\top}P_{t +1,i}^{\pmb{K}^{i},\pmb{K}^{-i*}}B_t^i)(E_{t,i}^{\pmb{K}}-E_{t,i}^{\pmb{K}^i,\pmb{K}^{-i*}})\Sigma_t^{\pmb{K}}\Sigma_t^{\pmb{K}^{i\prime},\pmb{K}^{-i*}}\Sigma_t^{\pmb{K}}\big)\label{eqn:vanilla_inte2}.
\end{eqnarray}
Now using the fact \eqref{eqn:trace_upper_bound} again with $\omega^2=\frac{2}{\sx^2}$, we can also bound  the second last term in \eqref{eqn:vanilla_inte1} as follows
\begin{eqnarray}
&& -4\eta\Tr\big(\Sigma_t^{\pmb{K}^{i\prime},\pmb{K}^{-i*}} (E_{t,i}^{\pmb{K}}\Sigma_t^{\pmb{K}}-E_{t,i}^{\pmb{K}^i,\pmb{K}^{-i*}}\Sigma_t^{\pmb{K}^{i\prime},\pmb{K}^{-i*}})^\top E_{t,i}^{\pmb{K}^i,\pmb{K}^{-i*}}\big)\nonumber\\
&=& -4\eta\Tr\big(\Sigma_t^{\pmb{K}^{i\prime},\pmb{K}^{-i*}} \big((E_{t,i}^{\pmb{K}}-E_{t,i}^{\pmb{K}^i,\pmb{K}^{-i*}})\Sigma_t^{\pmb{K}}+E_{t,i}^{\pmb{K}^i,\pmb{K}^{-i*}}(\Sigma_t^{\pmb{K}}-\Sigma_t^{\pmb{K}^{i\prime},\pmb{K}^{-i*}})\big)^\top E_{t,i}^{\pmb{K}^i,\pmb{K}^{-i*}}\big)\nonumber\\
&=&-4\eta\Tr\big(\Sigma_t^{\pmb{K}^{i\prime},\pmb{K}^{-i*}}\Sigma_t^{\pmb{K}}(E_{t,i}^{\pmb{K}}-E_{t,i}^{\pmb{K}^i,\pmb{K}^{-i*}})^\top E_{t,i}^{\pmb{K}^i,\pmb{K}^{-i*}}\big)\nonumber\\
&&-4\eta\Tr\big(\Sigma_t^{\pmb{K}^{i\prime},\pmb{K}^{-i*}}(\Sigma_t^{\pmb{K}}-\Sigma_t^{\pmb{K}^{i\prime},\pmb{K}^{-i*}}) (E_{t,i}^{\pmb{K}^i,\pmb{K}^{-i*}})^\top E_{t,i}^{\pmb{K}^i,\pmb{K}^{-i*}}\big)\nonumber\\
&\leq& 4\eta\frac{\sx^2}{4}\Tr\big((E_{t,i}^{\pmb{K}^i,\pmb{K}^{-i*}})^\top E_{t,i}^{\pmb{K}^i,\pmb{K}^{-i*}}\big)+\frac{4\eta}{\sx^2}\Tr\big(\Sigma_t^{\pmb{K}^{i\prime},\pmb{K}^{-i*}}\Sigma_t^{\pmb{K}}(E_{t,i}^{\pmb{K}}-E_{t,i}^{\pmb{K}^i,\pmb{K}^{-i*}})^\top(E_{t,i}^{\pmb{K}}-E_{t,i}^{\pmb{K}^i,\pmb{K}^{-i*}})\cdot\nonumber\\
&&\Sigma_t^{\pmb{K}}\Sigma_t^{\pmb{K}^{i\prime},\pmb{K}^{-i*}}\big)+4\eta\|\Sigma_t^{\pmb{K}^{i\prime},\pmb{K}^{-i*}}\|\,\|\Sigma_t^{\pmb{K}}-\Sigma_t^{\pmb{K}^{i\prime},\pmb{K}^{-i*}}\|\Tr\big( (E_{t,i}^{\pmb{K}^i,\pmb{K}^{-i*}})^\top E_{t,i}^{\pmb{K}^i,\pmb{K}^{-i*}}\big)\nonumber\\
&\leq& \eta\sx^2\Tr\big((E_{t,i}^{\pmb{K}^i,\pmb{K}^{-i*}})^\top E_{t,i}^{\pmb{K}^i,\pmb{K}^{-i*}}\big)+\frac{4\eta}{\sx^2}\|E_{t,i}^{\pmb{K}}-E_{t,i}^{\pmb{K}^i,\pmb{K}^{-i*}}\|^2\Tr\big(\Sigma_t^{\pmb{K}^{i\prime},\pmb{K}^{-i*}}\Sigma_t^{\pmb{K}}\Sigma_t^{\pmb{K}}\Sigma_t^{\pmb{K}^{i\prime},\pmb{K}^{-i*}}\big)\nonumber\\
&&+4\eta\|\Sigma_t^{\pmb{K}^{i\prime},\pmb{K}^{-i*}}\|\,\|\Sigma_t^{\pmb{K}}-\Sigma_t^{\pmb{K}^{i\prime},\pmb{K}^{-i*}}\|\Tr\big( (E_{t,i}^{\pmb{K}^i,\pmb{K}^{-i*}})^\top E_{t,i}^{\pmb{K}^i,\pmb{K}^{-i*}}\big).\label{eqn:inte5_vani}
\end{eqnarray}
Then plugging the above bounds into \eqref{eqn:vanilla_inte1} gives
\begin{eqnarray}     
&& C^{i\prime,-i*}-C^{i,-i*} \nonumber \\    
&& \leq \sum_{t=0}^{T-1}\Big[4\eta^2\Tr\big(\Sigma_t^{\pmb{K}^{i\prime},\pmb{K}^{-i*}}\Sigma_t^{\pmb{K}} (E_{t,i}^{\pmb{K}}-E_{t,i}^{\pmb{K}^i,\pmb{K}^{-i*}})^\top(R_t^{i}+(B_t^i)^{\top}P_{t +1,i}^{\pmb{K}^{i},\pmb{K}^{-i*}}B_t^i)(E_{t,i}^{\pmb{K}}-E_{t,i}^{\pmb{K}^i,\pmb{K}^{-i*}})\Sigma_t^{\pmb{K}}\big) \nonumber \\
&&\quad\quad+ 8\eta^2\frac{\sx}{4}\Tr\big((E_{t,i}^{\pmb{K}^i,\pmb{K}^{-i*}})^\top E_{t,i}^{\pmb{K}^i,\pmb{K}^{-i*}}\big)+ 8\eta^2\frac{1}{\sx}\Tr\big(\Sigma_t^{\pmb{K}}\Sigma_t^{\pmb{K}^{i\prime},\pmb{K}^{-i*}}\Sigma_t^{\pmb{K}} (E_{t,i}^{\pmb{K}}-E_{t,i}^{\pmb{K}^i,\pmb{K}^{-i*}})^\top\cdot \nonumber\\
&&\qquad\,\,\quad(R_t^{i}+(B_t^i)^{\top}P_{t +1,i}^{\pmb{K}^{i},\pmb{K}^{-i*}}B_t^i)(R_t^{i}+(B_t^i)^{\top}P_{t +1,i}^{\pmb{K}^{i},\pmb{K}^{-i*}}B_t^i)(E_{t,i}^{\pmb{K}}-E_{t,i}^{\pmb{K}^i,\pmb{K}^{-i*}})\Sigma_t^{\pmb{K}}\Sigma_t^{\pmb{K}^{i\prime},\pmb{K}^{-i*}}\Sigma_t^{\pmb{K}}\big) \nonumber\\
&& \quad\quad+4\eta^2\Tr\big(\Sigma_t^{\pmb{K}^{i\prime},\pmb{K}^{-i*}}\Sigma_t^{\pmb{K}} (E_{t,i}^{\pmb{K}^i,\pmb{K}^{-i*}})^\top(R_t^{i}+(B_t^i)^{\top}P_{t +1,i}^{\pmb{K}^{i},\pmb{K}^{-i*}}B_t^i)E_{t,i}^{\pmb{K}^i,\pmb{K}^{-i*}}\Sigma_t^{\pmb{K}}\big)\nonumber\\
&&\quad\quad+ \eta\sx^2\Tr\big((E_{t,i}^{\pmb{K}^i,\pmb{K}^{-i*}})^\top E_{t,i}^{\pmb{K}^i,\pmb{K}^{-i*}}\big)+\frac{4\eta}{\sx^2}\|E_{t,i}^{\pmb{K}}-E_{t,i}^{\pmb{K}^i,\pmb{K}^{-i*}}\|^2\Tr\big(\Sigma_t^{\pmb{K}^{i\prime},\pmb{K}^{-i*}}\Sigma_t^{\pmb{K}}\Sigma_t^{\pmb{K}}\Sigma_t^{\pmb{K}^{i\prime},\pmb{K}^{-i*}}\big)\nonumber\\
&&\quad\quad+4\eta\|\Sigma_t^{\pmb{K}^{i\prime},\pmb{K}^{-i*}}\|\,\|\Sigma_t^{\pmb{K}}-\Sigma_t^{\pmb{K}^{i\prime},\pmb{K}^{-i*}}\|\Tr\big( (E_{t,i}^{\pmb{K}^i,\pmb{K}^{-i*}})^\top E_{t,i}^{\pmb{K}^i,\pmb{K}^{-i*}}\big)\nonumber\\
&&\quad\quad-4\eta\sx^2\Tr\big( (E_{t,i}^{\pmb{K}^i,\pmb{K}^{-i*}})^\top E_{t,i}^{\pmb{K}^i,\pmb{K}^{-i*}}\big)\Big]\nonumber \\
&& \leq \sum_{t=0}^{T-1}\Big[\Big(4\eta^2\|\Sigma_t^{\pmb{K}^{i\prime},\pmb{K}^{-i*}}\|\,\|\Sigma_t^{\pmb{K}}\|^2\Tr\big(R_t^{i}+(B_t^i)^{\top}P_{t +1,i}^{\pmb{K}^{i},\pmb{K}^{-i*}}B_t^i\big)+ \frac{8\eta^2}{\sx}\|R_t^{i}+(B_t^i)^{\top}P_{t +1,i}^{\pmb{K}^{i},\pmb{K}^{-i*}}B_t^i\|^2\cdot\nonumber\\
&&\quad\quad\Tr\big(\Sigma_t^{\pmb{K}}\Sigma_t^{\pmb{K}^{i\prime},\pmb{K}^{-i*}}\Sigma_t^{\pmb{K}}\Sigma_t^{\pmb{K}} \Sigma_t^{\pmb{K}^{i\prime},\pmb{K}^{-i*}}\Sigma_t^{\pmb{K}}\big)+\frac{4\eta}{\sx^2} \Tr\big(\Sigma_t^{\pmb{K}^{i\prime},\pmb{K}^{-i*}}\Sigma_t^{\pmb{K}}\Sigma_t^{\pmb{K}}\Sigma_t^{\pmb{K}^{i\prime},\pmb{K}^{-i*}}\big)\Big)\cdot \nonumber\\
&& \quad\quad \|E_{t,i}^{\pmb{K}}-E_{t,i}^{\pmb{K}^i,\pmb{K}^{-i*}}\|^2 +\big(2\eta^2\sx+4\eta^2\|\Sigma_t^{\pmb{K}^{i\prime},\pmb{K}^{-i*}}\|\|\Sigma_t^{\pmb{K}}\|^2\|R_t^{i}+(B_t^i)^{\top}P_{t +1,i}^{\pmb{K}^{i},\pmb{K}^{-i*}}B_t^i\|\nonumber\\
&&\quad\quad\quad+\eta\sx^2 +4\eta\|\Sigma_t^{\pmb{K}^{i\prime},\pmb{K}^{-i*}}\|\,\|\Sigma_t^{\pmb{K}}-\Sigma_t^{\pmb{K}^{i\prime},\pmb{K}^{-i*}}\|-4\eta\sx^2\big)\Tr\big( (E_{t,i}^{\pmb{K}^i,\pmb{K}^{-i*}})^\top E_{t,i}^{\pmb{K}^i,\pmb{K}^{-i*}}\big)\Big],\label{eqn:Cdiff_inte3_vani}
\end{eqnarray}
where the first inequality holds by \eqref{eqn:vanilla_inte2} and \eqref{eqn:inte5_vani}, and the second inequality holds by the trace inequality \eqref{eqn:trace_inequality} and rearranging terms.
Now we bound the term $\|\Sigma_t^{\pmb{K}}\|$. By \eqref{eqn:intermedidate1} and \eqref{eqn:G_t_bound} we have
\begin{eqnarray}
\|\Sigma_t^{\pmb{K}}\| &=& \|\mathcal{G}_{t-1}(\Sigma_0)\| + \Big\|\sum_{s=1}^{t-1} D_{t-1,s} W D_{t-1,s}^{\top}\Big\| +\|W\|\nonumber\\
&\leq& \rho_{\pmb{K}}^{2T} \|\Sigma_0\| + \big(\rho_{\pmb{K}}^{2T}+1\big)\|W\|. \label{eqn:sigma_bd_vani}
\end{eqnarray}
By \eqref{eqn:one_step_noise_cond_vani} we have $\sx^2\geq g_3$, which implies
\begin{eqnarray}\label{eqn:vanilla_inte3}
    \sx^2 &\geq& 4\frac{10 TC^{i,-i*}}{(10 T-1)\sqmin}\Big(\frac{10 TC^{i,-i*}}{(10T-1)\sqmin}+\rho_{\pmb{K}}^{2T} \|\Sigma_0\| + \big(\rho_{\pmb{K}}^{2T}+1\big)\|W\| \Big)\nonumber\\
    &\geq&4 \|\Sigma_t^{\pmb{K}^{i\prime},\pmb{K}^{-i*}}\|\,\left(\|\Sigma_t^{\pmb{K}}\|+\|\Sigma_t^{\pmb{K}^{i\prime},\pmb{K}^{-i*}}\|\right)\nonumber\\
    &\geq& 4 \|\Sigma_t^{\pmb{K}^{i\prime},\pmb{K}^{-i*}}\|\,\|\Sigma_t^{\pmb{K}}-\Sigma_t^{\pmb{K}^{i\prime},\pmb{K}^{-i*}}\|.
\end{eqnarray}
Now, since $\big\|\Sigma_{\pmb{K}^{i\prime},\pmb{K}^{-i*}}\big\| \leq\frac{2\,C^{i,-i*}}{\sqmin}$ by \eqref{eqn:bd_pertSigma_vani}, we can bound the step size condition in \eqref{eqn:I1_conseq_2_vani} by
\begin{equation}\label{eqn:vanilla_inte4}
     \eta \leq \frac{\sx^2}{2\sx+4\|\Sigma_t^{\pmb{K}^{i\prime},\pmb{K}^{-i*}}\|\|\Sigma_t^{\pmb{K}}\|^2\|R_t^{i}+(B_t^i)^{\top}P_{t +1,i}^{\pmb{K}^{i},\pmb{K}^{-i*}}B_t^i\|}.
\end{equation}    
Combining \eqref{eqn:vanilla_inte3} and \eqref{eqn:vanilla_inte4} gives 
\begin{eqnarray}
&&2\eta^2\sx+4\eta^2\|\Sigma_t^{\pmb{K}^{i\prime},\pmb{K}^{-i*}}\|\|\Sigma_t^{\pmb{K}}\|^2\|R_t^{i}+(B_t^i)^{\top}P_{t +1,i}^{\pmb{K}^{i},\pmb{K}^{-i*}}B_t^i\|\nonumber\\
&&+\eta\sx^2 +4\eta\|\Sigma_t^{\pmb{K}^{i\prime},\pmb{K}^{-i*}}\|\,\|\Sigma_t^{\pmb{K}}-\Sigma_t^{\pmb{K}^i,\pmb{K}^{-i*}}\|-4\eta\sx^2\nonumber \\
&\leq & -\eta\sx^2.\label{eqn:inte_4_vani}
\end{eqnarray}
Hence, using this in \eqref{eqn:Cdiff_inte3_vani}, we have
     \begin{equation}
C^{i\prime,-i*}-C^{i,-i*} \leq \eta\,h_{\rm diff}^i\sum_{t=0}^{T-1}\|E_{t,i}^{\pmb{K}}-E_{t,i}^{\pmb{K}^i,\pmb{K}^{-i*}}\|^2  -\eta\sx^2 \sum_{t=0}^{T-1}\Tr\big( (E_{t,i}^{\pmb{K}^i,\pmb{K}^{-i*}})^\top E_{t,i}^{\pmb{K}^i,\pmb{K}^{-i*}}\big)\label{eqn:onestep_inte_vani},
\end{equation}
where 
\begin{equation*}
\begin{split}
 h_{\rm diff}^i &:= 4\eta\,k_i\frac{10 TC^{i,-i*}}{(10 T-1)\sqmin}\left(\rho_{\pmb{K}}^{2T} \|\Sigma_0\| + \big(\rho_{\pmb{K}}^{2T}+1\big)\|W\| \right)^2\left(\gamma_R+(\gamma_B)^2\frac{C^{i,-i*}}{\sx}\right)\\
 &\qquad+8\eta\frac{d}{\sx}\left(\frac{10 TC^{i,-i*}}{(10 T-1)\sqmin}\right)^2\left(\rho_{\pmb{K}}^{2T} \|\Sigma_0\| + \big(\rho_{\pmb{K}}^{2T}+1\big)\|W\| \right)^4\left(\gamma_R+(\gamma_B)^2\frac{C^{i,-i*}}{\sx}\right)^2\\
 &\qquad +4\frac{d}{\sx^2} \left(\frac{10 TC^{i,-i*}}{(10 T-1)\sqmin}\right)^2\left(\rho_{\pmb{K}}^{2T} \|\Sigma_0\| + \big(\rho_{\pmb{K}}^{2T}+1\big)\|W\| \right)^2.
 \end{split}
\end{equation*}
%Here \eqref{eqn:Cdiff_inte1} holds by Lemma ; \eqref{eqn:Cdiff_inte2} holds by \eqref{eqn:one_step_update} and Lemma ;
%\eqref{eqn:Cdiff_inte3} holds by  \eqref{eqn:onestep_inte} holds since, by step size condition \eqref{eqn:step_size},we have 
%by , which leads to

Therefore, by \eqref{eqn:diff_E_bd_vani}, \eqref{eqn:bd_pertSigma_vani}, Lemma \ref{lemma:grad_domi}, and Lemma \ref{lemma:nonzero_bds_P_Sigma},
\begin{eqnarray}
 C^{i\prime,-i*}-C^{i,-i*} &\leq & \eta\,h_{\rm diff}^i\,T\,\left[\frac{(\gamma_B)^2C^{i,-i*}}{\sx}\frac{2(\rho_{\pmb{K}}^{2T}-1)}{\rho_{\pmb{K}}^2-1}\sum_{j=1,j\neq i}^N\vertiii{\pmb{K}^{j}-\pmb{K}^{j*}}\right]^2 \nonumber \\ &&\quad-\eta\sx^2\sum_{t=0}^{T-1}\Tr\big( (E_{t,i}^{\pmb{K}^i,\pmb{K}^{-i*}})^\top E_{t,i}^{\pmb{K}^i,\pmb{K}^{-i*}}\big) \nonumber \\
     &\leq & \eta\,h_{\rm glob}\left(\sum_{j=1,j\neq i}^N\vertiii{\pmb{K}^{j}-\pmb{K}^{j*}}\right)^2-\eta\frac{\sx^2\srmin}{\|\Sigma_{\pmb{K}^{*}}\|}\big(C^{i,-i*}-C^{i*}\big), \label{eqn:ind_cost_diff_vani}
\end{eqnarray}
where
\begin{equation}\label{eqn:defn_hloc_vani}
 \begin{split}
      h_{\rm glob} &= 4T\left[\eta\,k_i\frac{10 T\max_i\{C^{i,-i*}\}}{(10 T-1)\sqmin}\left(\rho_{\pmb{K}}^{2T} \|\Sigma_0\| + \big(\rho_{\pmb{K}}^{2T}+1\big)\|W\| \right)^2\left(\gamma_R+(\gamma_B)^2\frac{\max_i\{C^{i,-i*}\}}{\sx}\right)\right.\\
      &\qquad+2\eta\frac{d}{\sx}\left(\frac{10 T\max_i\{C^{i,-i*}\}}{(10 T-1)\sqmin}\right)^2\left(\rho_{\pmb{K}}^{2T} \|\Sigma_0\| + \big(\rho_{\pmb{K}}^{2T}+1\big)\|W\| \right)^4\left(\gamma_R+(\gamma_B)^2\frac{\max_i\{C^{i,-i*}\}}{\sx}\right)^2\\
      &\qquad\left.+\,\frac{d}{\sx^2} \left(\frac{10 T\max_i\{C^{i,-i*}\}}{(10 T-1)\sqmin}\right)^2\left(\rho_{\pmb{K}}^{2T} \|\Sigma_0\| + \big(\rho_{\pmb{K}}^{2T}+1\big)\|W\| \right)^2\right]\\
      &\qquad\cdot\left[\frac{(\gamma_B)^2\max_i\{C^{i,-i*}\}}{\sx}\frac{2(\rho_{\pmb{K}}^{2T}-1)}{\rho_{\pmb{K}}^2-1}\right]^2.
 \end{split}
\end{equation}\\[.1in]

{\it Step 5:}
Finally we can establish the one step contraction. Using \eqref{eqn:ind_cost_diff_vani}, we have
\begin{eqnarray}
 C^{i\prime,-i*}-C^{i*} &=& C^{i\prime,-i*}-C^{i,-i*}+C^{i,-i*}-C^{i*}\nonumber\\
 &\leq& \left(1-\eta\frac{\sx^2\srmin}{\|\Sigma_{\pmb{K}^{*}}\|}\right)\big(C^{i,-i*}-C^{i*}\big) + \eta\,h_{\rm glob}\left(\sum_{j=1,j\neq i}^N\vertiii{\pmb{K}^{j}-\pmb{K}^{j*}}\right)^2.\label{eq:C_diff_vani}
\end{eqnarray}
Hence by Lemma \ref{lemma:nonzero_almost_smoothness} and \eqref{eqn:Ci_diff_Ki}, we have
\begin{eqnarray*}
\sum_{j=1,j\neq i}^N\left(C^{j,-j*} - C^{j*}\right) 
 \geq \frac{\sx\srmin}{T} \sum_{j=1,j\neq i}^N\vertiii{\pmb{K}^j-\pmb{K}^{j*}}^2 \geq \frac{\sx\srmin}{T(N-1)} \left(\sum_{j=1,j\neq i}^N\vertiii{\pmb{K}^j-\pmb{K}^{j*}}\right)^2,
\end{eqnarray*}
and thus
\begin{equation}
    \begin{split}
        C^{i\prime,-i*}-C^{i*}& \leq \big(1-\eta\frac{\sx^2\srmin}{\|\Sigma_{\pmb{K}^{*}}\|}\big)\big(C^{i,-i*}-C^{i*}\big) + \,\eta\,\, h_{\rm glob}\frac{T(N-1)}{\sx \srmin}\left(\sum_{j=1,j\neq i}^N (C^{j,-j*} - C^{j*}) \right).\label{eq:C_diff1_vani}
    \end{split}
\end{equation}
Summing up \eqref{eq:C_diff1_vani} for $i=1,\cdots,N$, we have
\begin{equation}\label{eq:c_sum_onestep_vani}
\sum_{i=1}^N \left(C^{i\prime,-i*}-C^{i*}\right)
\leq \left(1-\eta\frac{\sx^2\srmin}{\|\Sigma_{\pmb{K}^{*}}\|}+\eta(N-1)h_{\rm glob}\frac{T(N-1)}{\sx \srmin}\right)\left(\sum_{i=1}^N (C^{i,-i*} - C^{i*})\right).
\end{equation}
Since $\eta\leq I_4$, we have \eqref{eqn:I2_conseq_vani} and
then
 \begin{eqnarray*}
  h_{\rm diff}^i &\leq& (4+\frac{1}{20})\frac{d}{\sx^2} \left(\frac{10 T \max_i\{C^{i,-i*}\}}{(10 T-1)\sqmin}\right)^2\left(\rho_{\pmb{K}}^{2T} \|\Sigma_0\| + \big(\rho_{\pmb{K}}^{2T}+1\big)\|W\| \right)^2\\
  &\leq& 5\frac{d}{\sx^2} \left(\frac{ \max_i\{C^{i,-i*}\}}{\sqmin}\right)^2\left(\rho_{\pmb{K}}^{2T} \|\Sigma_0\| + \big(\rho_{\pmb{K}}^{2T}+1\big)\|W\| \right)^2,\quad\text{and}\quad h_{\rm glob} \leq \bar{h}_{\rm glob},
 \end{eqnarray*}
where $\bar{h}_{\rm glob}$ is given by
\begin{equation}\label{eqn:defn_hloc_bar_vani}
  \begin{split}
      \bar{h}_{\rm glob} &= 5\,T\,d\, \frac{\max_i\{C^{i,-i*}\}^4(\gamma_B)^4}{\sqmin^2\sx^4}\left[\frac{2(\rho_{\pmb{K}}^{2T}-1)}{\rho_{\pmb{K}}^2-1}\right]^2\left(\rho_{\pmb{K}}^{2T} \|\Sigma_0\| + \big(\rho_{\pmb{K}}^{2T}+1\big)\|W\| \right)^2.
  \end{split}
 \end{equation}
Under condition \eqref{eqn:one_step_noise_cond_vani}, we have
\begin{eqnarray*}
\sx^2 g_1 - \frac{\widetilde{g}_2}{\sx^5} = \frac{\sx^2\srmin}{\|\Sigma_{\pmb{K}^{*}}\|}-(N-1)^2\bar{h}_{\rm glob}\frac{T}{\sx\srmin}>0,
 \end{eqnarray*}
which indicates that $\alpha\eta>0$. Since $\eta\leq \frac{\|\Sigma_{\pmb{K}^{*}}\|}{\sx^2\srmin}$ by assumption, we have
 \begin{eqnarray*}
 \eta \leq \frac{\|\Sigma_{\pmb{K}^{*}}\|}{\sx^2\srmin}<\left( \frac{\sx^2\srmin}{\|\Sigma_{\pmb{K}^{*}}\|}-(N-1)^2\bar{h}_{\rm glob}\frac{T}{\sx\srmin}\right)^{-1} = \left(\sx^2 g_1 - \frac{\widetilde{g}_2}{\sx^5}\right)^{-1}.
 \end{eqnarray*}
Recall that in the statement we define $\alpha=\sx^2 g_1-\widetilde{g}_2/\sx^5$. Therefore we have $\alpha\eta<1$. Along with \eqref{eq:c_sum_onestep_vani}, we obtain the one-step contraction \eqref{eq:c_sum_onestep_final_vani}.
\end{proof}

We now explain the main differences between the above proof and the proof for the natural policy gradient method. Since for the vanilla method we have $K_t^{i'}-K_t^i=- 2\eta\,E_{t,i}^{\pmb{K}}\,\Sigma_t^{\pmb{K}}$ rather than $K_t^{i'}-K_t^i=- 2\eta\,E_{t,i}^{\pmb{K}}$ in the case of natural version, the extra $\Sigma_t^{\pmb{K}}$ term needs to be dealt with when calculating the individual cost difference $C^{i\prime,-i*}-C^{i,-i*}$. More precisely, the presence of $\Sigma_t^{\pmb{K}}$ causes more terms which need to be bounded (see equation \eqref{eqn:vanilla_inte2}) when finding an upper bound on $C^{i\prime,-i*}-C^{i,-i*}$ in equation \eqref{eqn:vanilla_inte1}, and the term $\|\Sigma_t^{\pmb{K}}\|$ also needs to be bounded above (see equation \eqref{eqn:sigma_bd_vani}). Due to these extra terms, the amount of noise needed in the system has a more complex form given in \eqref{eqn:one_step_noise_cond_vani}, where a new $g_3$ is introduced to guarantee equation \eqref{eqn:inte_4_vani} holds. This further demonstrates that the natural policy gradient method can be considered as a normalized version of the vanilla policy gradient method.

%% file: main.bbl
\begin{thebibliography}{10}

\bibitem{bagnell2003covariant}
J~Andrew Bagnell and Jeff Schneider.
\newblock Covariant policy search.
\newblock In {\em Proceedings of 18th International Joint Conference on
  Artificial Intelligence (IJCAI-03)}, pages 1019--1024, August 2003.

\bibitem{balduzzi2018mechanics}
David Balduzzi, Sebastien Racaniere, James Martens, Jakob Foerster, Karl Tuyls,
  and Thore Graepel.
\newblock The mechanics of n-player differentiable games.
\newblock In {\em International Conference on Machine Learning}, pages
  354--363. PMLR, 2018.

\bibitem{BasarOlsder1999}
Tamer Ba{\c{s}}ar and Geert~Jan Olsder.
\newblock {\em Dynamic Non-Cooperative Game Theory}.
\newblock Society for Industrial and Applied Mathematics, 1998.

\bibitem{basei2021logarithmic}
Matteo Basei, Xin Guo, Anran Hu, and Yufei Zhang.
\newblock Logarithmic regret for episodic continuous-time linear-quadratic
  reinforcement learning over a finite-time horizon.
\newblock {\em Available at SSRN 3848428}, 2021.

\bibitem{bowling2002multiagent}
Michael Bowling and Manuela Veloso.
\newblock Multiagent learning using a variable learning rate.
\newblock {\em Artificial Intelligence}, 136(2):215--250, 2002.

\bibitem{bu2019global}
Jingjing Bu, Lillian~J Ratliff, and Mehran Mesbahi.
\newblock Global convergence of policy gradient for sequential zero-sum linear
  quadratic dynamic games.
\newblock {\em arXiv preprint arXiv:1911.04672}, 2019.

\bibitem{FGKM2018}
Maryam Fazel, Rong Ge, Sham~M Kakade, and Mehran Mesbahi.
\newblock Global convergence of policy gradient methods for the linear
  quadratic regulator.
\newblock In {\em International Conference on Machine Learning}, pages
  1467--1476. PMLR, 2018.

\bibitem{fiez2020implicit}
Tanner Fiez, Benjamin Chasnov, and Lillian Ratliff.
\newblock Implicit learning dynamics in stackelberg games: Equilibria
  characterization, convergence analysis, and empirical study.
\newblock In {\em International Conference on Machine Learning}, pages
  3133--3144. PMLR, 2020.

\bibitem{foerster2017learning}
Jakob~N Foerster, Richard~Y Chen, Maruan Al-Shedivat, Shimon Whiteson, Pieter
  Abbeel, and Igor Mordatch.
\newblock Learning with opponent-learning awareness.
\newblock {\em arXiv preprint arXiv:1709.04326}, 2017.

\bibitem{hambly2020policy}
Ben~M Hambly, Renyuan Xu, and Huining Yang.
\newblock Policy gradient methods for the noisy linear quadratic regulator over
  a finite horizon.
\newblock {\em Available at SSRN}, 2020.

\bibitem{houthooft2016vime}
Rein Houthooft, Xi~Chen, Yan Duan, John Schulman, Filip De~Turck, and Pieter
  Abbeel.
\newblock Vime: Variational information maximizing exploration.
\newblock {\em arXiv preprint arXiv:1605.09674}, 2016.

\bibitem{huang2006large}
Minyi Huang, Roland~P Malham{\'e}, and Peter~E Caines.
\newblock Large population stochastic dynamic games: closed-loop mckean-vlasov
  systems and the nash certainty equivalence principle.
\newblock {\em Communications in Information \& Systems}, 6(3):221--252, 2006.

\bibitem{JSLGWZ2018}
Junqi Jin, Chengru Song, Han Li, Kun Gai, Jun Wang, and Weinan Zhang.
\newblock Real-time bidding with multi-agent reinforcement learning in display
  advertising.
\newblock In {\em Proceedings of the 27th ACM International Conference on
  Information and Knowledge Management}, pages 2193--2201, 2018.

\bibitem{kakade2001natural}
Sham~M Kakade.
\newblock A natural policy gradient.
\newblock {\em Advances in {N}eural {I}nformation {P}rocessing {S}ystems}, 14,
  2001.

\bibitem{lasry2007mean}
Jean-Michel Lasry and Pierre-Louis Lions.
\newblock Mean field games.
\newblock {\em Japanese journal of mathematics}, 2(1):229--260, 2007.

\bibitem{letcher2019differentiable}
Alistair Letcher, David Balduzzi, S{\'e}bastien Racaniere, James Martens, Jakob
  Foerster, Karl Tuyls, and Thore Graepel.
\newblock Differentiable game mechanics.
\newblock {\em The Journal of Machine Learning Research}, 20(1):3032--3071,
  2019.

\bibitem{mania2019certainty}
Horia Mania, Stephen Tu, and Benjamin Recht.
\newblock Certainty equivalence is efficient for linear quadratic control.
\newblock In {\em Advances in Neural Information Processing Systems}, pages
  10154--10164, 2019.

\bibitem{mazumdar2019policy}
Eric Mazumdar, Lillian~J Ratliff, Michael~I Jordan, and S~Shankar Sastry.
\newblock Policy-gradient algorithms have no guarantees of convergence in
  linear quadratic games.
\newblock {\em arXiv preprint arXiv:1907.03712}, 2019.

\bibitem{mazumdar2020gradient}
Eric Mazumdar, Lillian~J Ratliff, and S~Shankar Sastry.
\newblock On gradient-based learning in continuous games.
\newblock {\em SIAM Journal on Mathematics of Data Science}, 2(1):103--131,
  2020.

\bibitem{nocedal2006numerical}
Jorge Nocedal and Stephen Wright.
\newblock {\em Numerical optimization}.
\newblock Springer Science \& Business Media, 2006.

\bibitem{peters2008natural}
Jan Peters and Stefan Schaal.
\newblock Natural actor-critic.
\newblock {\em Neurocomputing}, 71(7-9):1180--1190, 2008.

\bibitem{rajeswaran2017towards}
Aravind Rajeswaran, Kendall Lowrey, Emanuel~V Todorov, and Sham~M Kakade.
\newblock Towards generalization and simplicity in continuous control.
\newblock {\em Advances in Neural Information Processing Systems}, 30, 2017.

\bibitem{roudneshin2020reinforcement}
Masoud Roudneshin, Jalal Arabneydi, and Amir~G Aghdam.
\newblock Reinforcement learning in nonzero-sum linear quadratic deep
  structured games: Global convergence of policy optimization.
\newblock In {\em 2020 59th IEEE Conference on Decision and Control (CDC)},
  pages 512--517. IEEE, 2020.

\bibitem{Saniuk1987}
J.~Saniuk and I.~Rhodes.
\newblock A matrix inequality associated with bounds on solutions of algebraic
  riccati and lyapunov equations.
\newblock {\em IEEE Transactions on Automatic Control}, 32(8):739--740, 1987.

\bibitem{schulman2015trust}
John Schulman, Sergey Levine, Pieter Abbeel, Michael Jordan, and Philipp
  Moritz.
\newblock Trust region policy optimization.
\newblock In {\em International Conference on Machine Learning}, pages
  1889--1897. PMLR, 2015.

\bibitem{SSS2016}
Shai Shalev-Shwartz, Shaked Shammah, and Amnon Shashua.
\newblock Safe, multi-agent, reinforcement learning for autonomous driving.
\newblock {\em arXiv preprint arXiv:1610.03295}, 2016.

\bibitem{singh2000nash}
Satinder Singh, Michael~J Kearns, and Yishay Mansour.
\newblock Nash convergence of gradient dynamics in general-sum games.
\newblock In {\em UAI}, pages 541--548. Citeseer, 2000.

\bibitem{song2019convergence}
Xinliang Song, Tonghan Wang, and Chongjie Zhang.
\newblock Convergence of multi-agent learning with a finite step size in
  general-sum games.
\newblock {\em arXiv preprint arXiv:1903.02868}, 2019.

\bibitem{wang2020reinforcement}
Haoran Wang, Thaleia Zariphopoulou, and Xun~Yu Zhou.
\newblock Reinforcement learning in continuous time and space: A stochastic
  control approach.
\newblock {\em J. Mach. Learn. Res.}, 21:198--1, 2020.

\bibitem{wang1986trace}
Sheng-De Wang, Te-Son Kuo, and Chen-Fa Hsu.
\newblock Trace bounds on the solution of the algebraic matrix riccati and
  lyapunov equation.
\newblock {\em IEEE Transactions on Automatic Control}, 31(7):654--656, 1986.

\bibitem{zhang2010multi}
Chongjie Zhang and Victor Lesser.
\newblock Multi-agent learning with policy prediction.
\newblock In {\em Twenty-fourth AAAI conference on artificial intelligence},
  2010.

\bibitem{zhang2020sample}
Junzi Zhang, Jongho Kim, Brendan O'Donoghue, and Stephen Boyd.
\newblock Sample efficient reinforcement learning with reinforce.
\newblock {\em arXiv preprint arXiv:2010.11364}, 2020.

\bibitem{zhang2019policy}
Kaiqing Zhang, Zhuoran Yang, and Tamer Ba{\c{s}}ar.
\newblock Policy optimization provably converges to {N}ash equilibria in
  zero-sum linear quadratic games.
\newblock In {\em Advances in Neural Information Processing Systems}, 2019.

\bibitem{zhang2021derivative}
Kaiqing Zhang, Xiangyuan Zhang, Bin Hu, and Tamer Ba{\c{s}}ar.
\newblock Derivative-free policy optimization for risk-sensitive and robust
  control design: Implicit regularization and sample complexity.
\newblock {\em arXiv preprint arXiv:2101.01041}, 2021.

\end{thebibliography}
